\documentclass{amsart}
\pdfoutput=1
\usepackage[table,usenames,dvipsnames]{xcolor}
\usepackage{amsmath}
\usepackage{amssymb}
\usepackage{enumitem}
\usepackage[firstinits=true,doi=false,url=true,backend=bibtex]{biblatex}
\usepackage{graphicx}
\usepackage{hyperref}
\hypersetup{
    colorlinks,
    linkcolor={red!50!black},
    citecolor={blue!90!black},
    urlcolor={blue!80!black}
}

\usepackage{xcolor}
\usepackage{bm}
\usepackage{enumitem}
\setlist{nolistsep}
\newcommand{\epsi}{\ensuremath{\varepsilon}}

\newcommand{\Real}{\mathbb{R}}
\newcommand{\Pcal}{\mathcal{P}}
\newcommand{\grad}{\nabla}
\newcommand{\gradb}{\bar\nabla}
\newcommand{\dive}{\grad\cdot }
\newcommand{\Fcal}{\mathcal{F}}
\newcommand{\Tcal}{\mathcal{T}}

\newcommand{\sym}{\text{sym}}
\renewcommand{\vec}[1]{\bm{#1}}
\newcommand{\abs}[1]{\left|#1\right|}
\newcommand{\ddiv}{\texttt{div}}
\newcommand{\dis}[1]{\ensuremath{\texttt{#1}}}
\newcommand{\fracpar}[2]{\frac{\partial #1}{\partial #2}}
\DeclareMathOperator{\id}{I}
\DeclareMathOperator{\trace}{trace}

\title{Comparison between cell-centered and nodal based discretization schemes for
  linear elasticity}

\author{Halvor M. Nilsen \and Jan Nordbotten \and Xavier Raynaud}

\newcounter{mycase}
\newcounter{mysubcase}[mycase]

\begin{document}

\maketitle
\begin{abstract}
  In this paper we study newly developed methods for linear elasticity on polyhedral
  meshes. Our emphasis is on applications of the methods to geological models.
  Models of subsurface, and in particular sedimentary rocks, naturally lead to
  general polyhedral meshes. Numerical methods which can directly handle such
  representation are highly desirable. Many of the numerical challenges in simulation
  of subsurface applications come from the lack of robustness and accuracy of
  numerical methods in the case of highly distorted grids. In this paper we
  investigate and compare the Multi-Point Stress Approximation (MPSA) and the Virtual
  Element Method (VEM) with regards to grid features that are frequently seen in
  geological models and likely to lead to a lack of accuracy of the methods. In
  particular we look how the methods perform near the incompressible limit. This work
  shows that both methods are promising for flexible modeling of subsurface
  mechanics.
\end{abstract}

\section{Introduction}

Modeling of sedimentary subsurface rock naturally leads to general unstructured grids
because of stratigraphic layering, erosion and faults. The industry standard for
grids in reservoir modeling is the Corner-Point grids (cp-grids). Other geometrical
grid formats have been proposed to improve on this format, but all compact
representations of the underlying geology will lead to cells with high aspect ratios,
distorted cells, large variations in cell volumes and faces areas. Methods that are
valid on general polyhedral grids and are robust for different grid types will
greatly simplify the modeling of subsurface physics for multiphase flow encountered
in the oil and gas industry. The workhorse method there is the finite volume
discretization based on a two point flux approximation (TPFA). The method is not
convergent for general grids and can introduce large grid-orientation effects, see
for example \cite[Figure 3]{NLN11:faults}, but is very robust due to its monotonicity
properties, which often result in faster computation times. The multi-point flux
approximation (MPFA) method has been developed to solve the convergence problems and
has been successfully applied to minimize grid-orientation effects
\cite{Aavatsmark02,MRST12:comg}, but due to lack of monotonicity, the method is
difficult to apply to complex grids such as those arising from real reservoir
models. Based on a mixed formulation, the mimetic finite difference method has been
proposed for incompressible flow \cite{brezzi2005family,brezzi2005convergence}, but
problems arise in the case of fully compressible black oil models, as the method
introduces non-monotonicity and significantly more degrees of freedom. In recent
years, coupling of geo-mechanical effects with subsurface flow has become more
important in many areas including oil and gas production from mature fields,
fractured tight reservoirs as well as geothermal application and risk assessment of
CO2 injection. Realistic modeling of these geological cases is hampered by
differences in the way geo-mechanics and flow models are built and discretized.

Recently, a cell-centered finite volume discretization has been proposed in
\cite{nordbotten2014cell} to specifically address problems arising in coupled
geo-mechanical and flow simulation of porous media. The method is inspired from the
MPFA discretization developed for flow problems and was thus named multi-point stress
approximation (MPSA). The MPSA method presents two appealing features for subsurface
applications. Since the method is based on the MPFA method, it shares the same data
structure which is commonly used for the flow problem where the preferred methods
remain based on finite volume discretization. Moreover, the method can operate on the
type of general polyhedral grids typically used to represent complex heterogeneous
medium. This later property is shared by the virtual element (VE) method
\cite{beirao2013basic}. The VE method builds upon the long-standing effort in the
development of mimetic finite difference (MFD) methods, see
\cite{lipnikov2014mimetic, da2014mimetic}. The MFD method reproduces at the discrete
level fundamental properties of the differentiation operators, using only the
available degrees of freedom and without explicitly constructing any finite element
basis. In this way, the method can easily handle general cell shapes. The VE method
is a reformulation of the MFD method in the finite element framework. As in the MFD
method, a complete finite element basis for a polyhedral cell is not computed, some
of the basis elements become \textit{virtual}. Both the MPSA and the VE methods for
mechanics naturally define the divergence of the displacement on cells (see
\cite{NorbottenStable2015} for MPSA), which is also the natural coupling term between
flow and mechanics, when flow is discretized with finite volume methods.  As pointed
out in \cite{nordbotten2014cell}, any attempt to extend the TPFA method to mechanics
is bound to fail as the method already fails the local patch test. The local patch
test verifies that the numerical method preserves rigid rotations, which are exact
solutions to the problem.

In this paper we will investigate the MPSA and VE methods for mechanics with special
emphasis on grid artifacts that naturally occur in geological models of
sedimentary. Even if both aspects are related, our first interest is not the
convergence properties of the methods but their performance on coarse and distorted
meshes. This paper contains the first set of tests where the MPSA method is tested in
view of applications to geosciences. In addition, we will discuss the properties of
the methods in the incompressible limit since it has practical consequences for the
short time dynamics of elasticity problems coupled with flow, as for example the
Biot's equations. We also look at the different properties of the methods for
different types of grids and how the methods can incorporate features like fractures.

\section{Presentation of the methods}

We study the methods for the standard equations of linear elasticity given by
\begin{equation}
  \label{eq:lin_elast_cont}
  \begin{aligned}
    \dive \sigma &= \vec{f},\\
    \epsi &= \frac{1}{2}(\grad +\grad^{T})\vec{ u},\\
    \sigma &= C \epsi,
  \end{aligned}
\end{equation}
where $\sigma$ is the Cauchy stress tensor, $\epsi$ the infinitesimal strain tensors
and $u$ the displacement field. The linear operator $C$ is a fourth-order stiffness
tensor. Since both $\sigma$ and $\epsilon$ are symmetric, the Voigt notation is
convenient. In Voigt notation, a three-dimensional symmetric tensor $\{\epsi_{ij}\}$
is represented as an element of $\Real^6$ with components
$[\epsi_{11},\epsi_{22},\epsi_{33},\epsi_{23},\epsi_{13},\epsi_{12}]^T$ while a
two-dimensional symmetric tensor is represented by a vector in $\Real^3$ given by
$[\epsi_{11},\epsi_{22},\epsi_{12}]^T$. For isotropic materials, we have the
constitutive equations
\begin{equation}
  \label{eq:isosigma}
  \sigma = 2\mu\epsilon + \lambda\trace(\epsilon)\id.
\end{equation}
We summarize the description of the methods given in \cite{gain2014} for the VE
method and in \cite{nordbotten2014cell} for MPSA. In the case of VE, we do not use
the nodal representations of the load and traction terms. Instead we use traction and
load terms defined on faces and cells, respectively. This is consistent with the
physical meaning of these terms in addition to the fact that the integration rules
hold exactly. The advantages of this evaluation of the volume force will be discussed
in a forthcoming paper.

\subsection{The Virtual Element Method}
\label{subsec:vem}
As the classical finite element method, the VE method starts from the linear
elasticity equations written in the weak form
\begin{equation}
  \label{eq:weakform}
  \int_{\Omega}\epsi( \vec{v}):C \epsi (\vec{u})  \,dx = \int_{\Omega} \vec{v}\cdot\vec{f} \,dx  \quad \text{for all} \quad \vec{v}.
\end{equation}
In \eqref{eq:weakform}, we use the standard scalar product for matrices defined as
\begin{equation*}
  \alpha:\beta = \trace(\alpha^t\beta) = \sum_{i,j=1}^{3}\alpha_{i,j}\beta_{i,j},
\end{equation*}
for any two matrices $\alpha,\beta \in \Real^{3\times 3}$. We have also introduced
the symmetric gradient $\epsi$ given by
\begin{equation*}
  \epsi(u) = (\grad + \grad^{T} )\vec{u},
\end{equation*}
for any displacement $\vec{u}$. The fundamental idea in the VE method is to compute
on each element an approximation $a_K^h$ of the bilinear form
\begin{equation}
  \label{eq:defak}
  a_K(\vec{u}, \vec{v}) =  \int_{K} \epsi(\vec{u}) : C \epsi(\vec{v})\,dx,
\end{equation}
that, in addition of being symmetric, positive definite and coercive (uniformly with
respect to the grid size if we want convergence), is also exact for linear
functions. Note that in this paper, we only consider first-order methods. If higher
order methods are used, the exactness must hold for polynomials of a given degree
where the degree determines the order of the method. These methods were first
introduced as mimetic finite element methods but later developed further under the
name of virtual element methods (see \cite{da2014mimetic} for discussions). The
degrees of freedom are chosen as in the standard finite element methods to ensure the
continuity at the boundaries and an element-wise assembly of the bilinear forms
$a_K^h$. We have followed the implementation described in \cite{gain2014}. In a
first-order VE method, the projection operator $\Pcal$ into the space of linear
displacement has to be computed locally for each cell. The VE approach ensures that
it can be computed exactly for each basis element. The projection operator is defined
with respect to the metric induced by the bilinear form $a_K$. The projection is
self-adjoint so that we have the following Pythagoras identity,
\begin{equation}
  \label{eq:pythagoras}
  a_K(\vec{u}, \vec{v}) = a_K(\Pcal\vec{u}, \Pcal\vec{v}) + a_K((\id - \Pcal)\vec{u}, (\id - \Pcal)\vec{v})
\end{equation}
for all displacement field $\vec{u}$ and $\vec{v}$ (In order to keep this
introduction simple, we do not state the requirements on regularity which is needed
for the displacement fields). In \cite{gain2014}, an explicit expression for $\Pcal$
is given so that we do not even have to compute the projection. Indeed, we have
$\Pcal = \Pcal_{R} + \Pcal_{C}$ where $\Pcal_{R}$ is the projection on the space $R$
of pure rotations and $\Pcal_{C}$ the projection on the space $C$ of constant shear
strain. The spaces $R$ and $C$ are defined as
\begin{align*}
  R &= \left\{ \vec{a} + B(\vec{x} - \bar{\vec{x}})\ |\ \vec{a}\in\Real^3, B\in\Real^{3\times 3}, B^T = -B \right\},\\
  C &= \left\{ B(\vec{x} - \bar{\vec{x}})\ |\ B\in\Real^{3\times 3}, B^T = B \right\}.
\end{align*}
Then, the discrete bilinear form $a_K^h$ is defined as
\begin{equation}
  \label{eq:disenervem}
  a_K^h(\vec{u}, \vec{v}) = a_K(\Pcal\vec{u}, \Pcal\vec{v}) + s_K((\id - \Pcal)\vec{u},(\id - \Pcal)\vec{v})
\end{equation}
where $s_K$ is a symmetric positive matrix which is chosen such that $a_K^h$ remains
coercive. Note the similarities between \eqref{eq:disenervem} and
\eqref{eq:pythagoras}. Since $\Pcal_{R}$ and $\Pcal_{C}$ are orthogonal and
$\Pcal_{R}$ maps into the null space of $a_K$ (rotations do not produce any change in
the energy), we have that the first term on the right-hand side of
\eqref{eq:pythagoras} and \eqref{eq:disenervem} can be simplified to
\begin{equation*}
  a_K(\Pcal\vec{u}, \Pcal\vec{v}) = a_K(\Pcal_C\vec{u}, \Pcal_C\vec{v}).
\end{equation*}
The expression \eqref{eq:disenervem} immediately guarantees the consistency of the
method, as we get from \eqref{eq:disenervem} that, for linear displacements, the
discrete energy coincides with the exact energy. Since the projection operator can be
computed exactly for all elements in the basis - and in particular for the
\textit{virtual} basis elements for which we do not have explicit expressions - the
local matrix can be written only in terms of the degrees of freedom of the method. In
our case the degrees of freedom of the method are the value of displacement at the
node. Let us denote $\vec{\varphi}_i$ a basis for these degrees of freedom. The
matrix $(A_{K})_{i,j}=a_K(\vec{\varphi}_i, \vec{\varphi}_j)$ is given by
\begin{equation}
  \label{eq:assembvem}
  A_K =  |K| \ \ W_C^T D W_C + (I-\Pcal)^T  S_K (I-\Pcal).
\end{equation}
In \eqref{eq:assembvem}, $W_C$ is the projection operator from the values of node
displacements to the space of constant shear strain and $S_K$, which corresponds to a
discretization of $s_K$ in \eqref{eq:disenervem}, is a symmetric positive matrix
which guarantees the positivity of $A_{K}$. There is a large amount of freedom in the
choice of $S_K$ but it has to be scaled correctly. We choose the same $S_K$ as in
\cite{gain2014}. The matrix $D$ in \eqref{eq:assembvem} corresponds to the tensor $C$
rewritten in Voigt notations so that, in three dimensions, we have
\begin{equation*}
  D_{ij} = \epsi_i : C \epsi_j,\quad \text{ for }i,j = 1, \ldots, 6.
\end{equation*}
Finally, the matrices $A_K$ are used to assemble the global matrix $A$ corresponding
to $a^h$.

\subsection{Multi-Point Stress Approximation}
\label{subsec:mpsaexpl}
\begin{figure}[h]
  \centering
  \includegraphics[width = \textwidth]{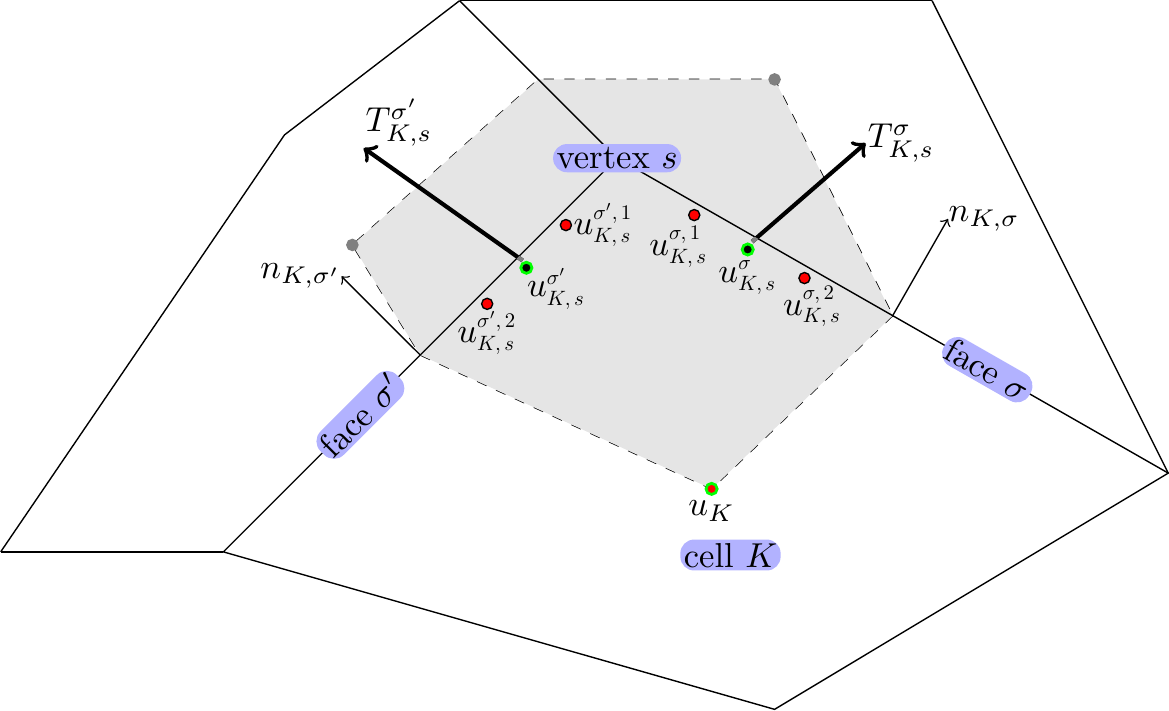}
  \caption{Illustration picture for the MPSA method. The interaction region for the
    vertex $s$ is represented in gray.  The degrees of freedom for the interaction
    sub-region $(K, s)$ are filled in red. We have represented the points that are
    used to define the discrete gradient with a green contour. We hope the other
    notations are self-explanatory, see \cite{nordbotten2015convergence} for a
    complete description.}
  \label{fig:mpsanotation}
\end{figure}

The Multi-Point Stress Approximation (MPSA) has its origin in the MPFA method
\cite{Aavatsmark02} which is a finite volume method for fluid flow. Its derivation is
based on discrete principles for the conservation of momentum and the continuity of
the forces. We use the same notations as in \cite{nordbotten2015convergence}, which
are also summarized in Figure \ref{fig:mpsanotation} (Note that $\sigma$ in this
section no longer denotes the Cauchy stress tensor but a face of a cell). On each
interaction sub-region, say $(K,s)$, we consider the degrees of freedom that are
given by a cell-value, $u_K$, and values at the two quadrature points on each
sub-face,
\begin{equation}
  \label{eq:defsubfacevar}
\left\{u_{K,s}^{\sigma,\beta}\right\}_{\sigma\in\Fcal_s\cap\Fcal_K,\beta=\{1,2\}}.
\end{equation}
Here, $\Fcal_s$ and $\Fcal_K$ denote the set of faces which have non empty
intersection with vertex $s$ and cell $K$, respectively. To simplify the
presentation, we consider only the 2D problem so that the displacement $u_K$, for
example, belongs to $\Real^2$. Note that the degrees of freedom at the quadrature
points are useful for deriving the method but will be removed in the assembly
process. On each outer face of an interaction sub-region, say $(K,s,\sigma)$, we can
define the average value
\begin{equation}
  \label{eq:defavervals}
  u_{K,s}^\sigma = \sum_{\beta=1}^{2}\omega_\beta u_{K,s}^{\sigma,\beta}
\end{equation}
using the Gauss quadrature weights $\omega_\beta$. From the average values computed
in \eqref{eq:defavervals}, we can define, uniquely in 2D and with some restrictions
in 3D (see \cite{nordbotten2015convergence}), a gradient operator which corresponds
to the linear approximation that takes the values $u_K$ at the cell center and
$u_{K,s}^\sigma$ at the center of each sub-face. We denote this gradient
operator by $\gradb_{K,s}$ and it is a mapping which, from the degrees of freedom of
the interaction sub-region, yields a two-dimensional tensor,
\begin{equation}
  \label{eq:defdisgrad}
  u_K,\left\{u_{K,s}^{\sigma,\beta}\right\}_{\substack{\sigma\in\Fcal_s\cap\Fcal_s\\\beta=\{1,2\}}} \to \gradb_{K,s}u.
\end{equation}
Here, the arrow means that the values at the right are computed using the quantities
at the left. We use the same convention below. Now that the discrete gradient
$\gradb_{K,s}$ has been defined, we approximate the forces on the sub-faces as
\begin{equation}
  \label{eq:defTkssigma}
  T_{K,s}^{\sigma} = m_\sigma^{s}\left(2\mu_{K}\gradb_{K, s}^{\sym}(u) + \lambda_K(\gradb_{K, s}\cdot u) I\right)\cdot n_{K,\sigma}
\end{equation}
where $\gradb_{K, s}^{\sym}(u)=\frac12(\gradb_{K, s}(u) + \gradb_{K, s}(u)^T)$ is the
discrete symmetric gradient operator and $\gradb_{K,s}\cdot u = \trace(\gradb_{K,s}(u))$. Equation \eqref{eq:defTkssigma} is a direct discrete transcription of the
constitutive equation \eqref{eq:isosigma}. The force acting on a cell-face is
naturally defined as the sum of the forces acting on all the corresponding sub-faces,
that is
\begin{equation}
  \label{eq:defTJsigma}
  T_{K}^{\sigma} = \sum_{\{s\,|\,\sigma \in \Fcal_s\}} T_{K,s}^{\sigma}.
\end{equation}
We get the first part of the discrete system of equations by imposing conservation of
momentum: For each cell, the sum of the forces applied to all faces is equal to the
external force applied to the cell, that is
\begin{equation}
  \label{eq:consmomentum}
  \sum_{\sigma\in\Fcal_K} T_{K}^{\sigma} = \int_{K}f(x)\,dx.
\end{equation}
The second part of the system of equations is obtained by defining the linear
interpolation operator $I_{FV,s}$ which from cell values yields all the remaining
degrees of freedom in the interaction region,
\begin{equation}
  \label{eq:defInterpoper}
  \left\{u_{K}\right\}_{K\in\Tcal_s} \overset{I_{FV,s}}{\longrightarrow} \left\{u_{K}\right\}_{K\in\Tcal_s}, \left\{u_{K,s}^{\sigma,\beta}\right\}_{\substack{\sigma\in\Fcal_s,K\in\Tcal_s\\\beta=\{1,2\}}}.
\end{equation}
Here, $\Tcal_s$ denotes the set of cells which contains the vertex $s$. We will see
shortly how $I_{FV,s}$ is determined. Assuming that it is defined, we can use
\eqref{eq:defTkssigma} and \eqref{eq:defdisgrad} to compute $T_{K,s}^{\sigma}$ by
only using the cell values $u_{K'}$ for $K'\in\Tcal_s$. Schematically, we have
\begin{equation}
  \label{eq:Tkssigmasimp}
   \left\{u_{K'}\right\}_{K'\in\Tcal_s} \overset{I_{FV,s}}{\longrightarrow}
  \begin{array}[h]{l}
    \left\{u_{K}\right\}_{K\in\Tcal_s},\\ \left\{u_{K,s}^{\sigma,\beta}\right\}_{\substack{\sigma\in\Fcal_s,K\in\Tcal_s\\\beta=\{1,2\}}}
  \end{array}
  \overset{\gradb_{K,s}}{\longrightarrow} (\gradb u)_{K,s} \overset{\text{by \eqref{eq:defTkssigma}}}{\longrightarrow}
  T_{K,s}^{\sigma}.
\end{equation}
It is important to note that the sequence of operations given by
\eqref{eq:Tkssigmasimp} is finally local in the sense that it only involves
cell-values in the interaction region $\Tcal_s$ of the node $s$. The operator $I_{FV,
  s}$ takes essentially care of this local reduction and we can now explain how the
coefficients of this linear mapping are determined. First, we require that
the forces are continuous at each face,
\begin{equation}
  \label{eq:contforce}
  T_{K,s}^{\sigma}  = - T_{K',s}^{\sigma}
\end{equation}
whenever $\sigma = K\cap K'$. The remaining degrees of freedom to define
$I_{FV, s}$ are not enough to impose the continuity of the displacement. Instead, we
use them to minimize the jump of the displacement at the interface. Thus, the
coefficients of $I_{FV,s}$ are determined by solving the least square problem
\begin{equation}
  \label{eq:leastsqpb}
  \min \sum_{\sigma \in \Fcal_s}\sum_{\beta\in\{1,2\}}\sum_{K,K'\in\Tcal_\sigma}w_{K'.K}\left|u_{K',s}^{\sigma,\beta} - u_{K,s}^{\sigma,\beta}\right|^2
\end{equation}
with the constraints given by \eqref{eq:contforce}. The weights $w_{K,K'}$ can be
chosen as the harmonic mean of the largest eigenvalue of the stiffness tensor $C$ of
the adjacent cells $K$ and $K'$. Once this is done, the result of the assembly
process leads us to a linear mapping of the form
\begin{equation}
  \label{eq:TKssigmatrans}
  T_{K,s}^{\sigma} = \sum_{K'\in\Tcal s} t_{K, K', s, \sigma}u_{K'}.
\end{equation}
The local coefficient tensors $t_{K, K', s, \sigma}$ are referred to as sub-face
stress weight tensors, and generalize the notion of transmissibilities from the
scalar diffusion equation \cite{nordbotten2014cell}. The stress continuity condition
\eqref{eq:contforce} implies that $t_{K, K', s, \sigma}$ whenever $\sigma = K\cap
K'$. The system of equations for linear elasticity are then given by the discrete
conservation of momentum \eqref{eq:consmomentum}, the definition of the force on
faces \eqref{eq:defTJsigma} and the multi-point approximation of sub-face forces
given by \eqref{eq:TKssigmatrans}.

\clearpage

\subsection{Fundamental differences between the methods}

There are fundamental differences between the VE and finite element (FE) methods in
the assembly process. In the VE framework, the matrix $A$ defining the discrete
equation $A u = F$ is computed element for element, based on rock parameters and the
geometry of the cell. In the MPSA method, we first calculate fluxes from cell center
displacement. This calculation requires to solve the singular value problem which
corresponds to \eqref{eq:leastsqpb}. Then the contribution to a matrix element is
calculated by summing up the contributions from each sub-face. So VE operates on the
element, while MPSA operates on interaction regions. Interaction regions can be
associated with the dual grid. The MPSA method considered here also needs to solve a
constrained least-square problem on each interaction region.  Recently, a new variant
of the MPSA method has been developed based on ideas from weakly symmetric mixed
finite element discretizations. This new variant circumvents the least squares
problem by enforcing displacement continuity at a single point for each face within
each interaction region, see \cite{Keilegavlen2016}.

\subsubsection{Comparison of the method with respect to grid features}

In the context of geosciences, the MPSA method has the advantage to allow for an easy
treatment of fractures. A fracture appearing at the interface between two cells can
be modeled by decoupling the corresponding face. If we denote this face by $\sigma$,
Equation \eqref{eq:contforce} is replaced by $T_{K,s}^{\sigma} = 0 = -
T_{K',s}^{\sigma}$ and the displacement values at the Gauss points,
$u_{K,s}^{\sigma,\beta}$ and $u_{K',s}^{\sigma,\beta}$ for $K,K'\in\Tcal_\sigma$ are
removed from the sum in \eqref{eq:leastsqpb}. The method, before removing the degrees
of freedom, is similar to a mixed method. This can be seen more explicitly in
\cite{Keilegavlen2016}, where the MPSA method which is presented there is very
similar to the PEERS elements for triangles \cite{arnold1984peers}. The difference is
that the PEERS elements have one set of forces on edges and an addition bubble
function to obtain stability for the incompressible limit, while the MPSA method has
two sets of forces on the edges. To reduce the degrees of freedom of the global
system, the MPSA method sacrifices the symmetry and positive definiteness of the
local systems in order to make a block diagonal inner-product which can be
reduced. In the case of triangles, there exists a symmetric block diagonal
inner-product as noticed in \cite{lipnikov2009local}, which makes the formulation
\cite{Keilegavlen2016} attractive.

The disadvantages of the MPSA method is that it is not symmetric and only
conditionally stable, a property which is also encountered in the MPFA method
\cite{lipnikov2009local,klausen2006robust}. It may result in failure or poor results
for severely distorted grids, and strongly anisotropic media. However, the stability
of the method can be verified locally \cite{nordbotten2015convergence}. Still, this
may prevent it from being used on specific grids without extra griding. Generally the
MPSA will suffer from the same grid restrictions as MPFA.  The method requires the
inversion of local matrices, which may induce an extra cost, but this only affects
the assembly process. In the case of VEM, we can expect the same structure as for FEM
so that the same solvers can be used. The system matrix is by construction symmetric,
positive and definite. If not modified, the VE method suffers from the same
limitations with respect to locking and accuracy of stresses as a linear FE
method. In addition, forces are not as naturally defined on faces as they are in MPSA
and methods of mixed type. Some of these problems can be avoided by using techniques
developed for the FEM solutions \cite{zienkiewicz92SPR} since for simple grids these
the VE and FE methods are equivalent.

\subsubsection{Limit of incompressible elasticity}
\label{subsubsec:locking}

In the limit of incompressible elasticity, the displacement field tends to be
solenoidal, that is, divergence free. Numerical locking occurs when solenoidal fields
are poorly approximated at the discrete level. In the context of subsurface
application, numerical locking will be an issue when considering the coupling of
linear elasticity for the rock matrix with flow. To see that, let us consider the
Biot's equations \cite{biot1941general} which are commonly used in the simulation of
hydromechanical problems. The Biot's equations consist of the following linear
equations,
\begin{equation}
  \label{eq:poroelast}
  \begin{aligned}
    \dive \sigma + \grad p&= \vec{f},\\
    \fracpar{}{t}(c_0 p + \alpha\dive \vec{u}) + \dive(K\grad p)  &= 0 ,\\
    \epsi &= \frac{1}{2}(\grad +\grad^{T})\vec{ u},\\
    \sigma &= C \epsi,
  \end{aligned}
\end{equation}
where $c_0 p + \alpha\dive \vec{u}$ denotes the fluid content. The fluid content
depends on the fluid pressure $p$ and on the rock volume change given by
$\dive\vec{u}$ which is weighted by the Biot-Willis constant $\alpha$. We discretize
in time the equations \eqref{eq:poroelast} and use a superscript $n$ to denote the
values corresponding to the time step $n$. From \eqref{eq:poroelast}, we get
\begin{subequations}
  \label{eq:poroelastdisc}
  \begin{align}
    \dive \sigma^{n+1} + \grad p^{n+1}&= \vec{f}^{n+1},\\
    \label{eq:poroelastdisc2}
    c_0 p^{n+1} + \alpha\dive\vec{u}^{n+1} + \Delta t\,\dive(K\grad (p^{n+1}))  &=  c_0p^{n} + \alpha\dive\vec{u}^{n} ,\\
    \epsi^{n+1} &= \frac{1}{2}(\grad +\grad^{T})\vec{u}^{n+1},\\
    \sigma^{n+1} &= C \epsi^{n+1}.
  \end{align}
\end{subequations}
In equation \eqref{eq:poroelastdisc2}, we can see that, in the limit where the fluid
becomes incompressible, that is $c_0\approx 0$, when the time-step $\Delta t$ tends
to zero, the change in the divergence of $\vec{u}$ becomes very small. In this case,
we are computing a displacement field which is close to solenoidal and numerical
locking will potentially become an issue.

In the case of VEM and any finite element method, the material parameters enters the
discrete equations \emph{cell-wise} in the assembly of the bilinear form $a$, see
\eqref{eq:disenervem}. Letting $\lambda$ be very large compared to $\mu$ therefore
imposes a near solenoidality condition on each cell. To evaluate the level of
locking, we can compare the number of degree of freedom of the whole system with the
number of local solenoidal equations, that is, the equations which locally impose the
solenoidal condition for large $\lambda$. The heuristic is then the following: If the
number of local solenoidal equations is small with respect to the number of degree of
freedoms, then we increase the chance that the global discrete divergence operator
becomes surjective. In this case, we avoid the appearance of spurious mode which is
also responsible for locking, see \cite[section 8.3]{da2014mimetic} where this aspect
is discussed for the Stokes equation.

In the case of VEM, the ratio between the number of cells and the number of nodes
will give an indication of the sensitivity of the grid with respect to locking. The
higher this ratio is, the more likely it is that locking appears. Hence, triangular
grid, where this ratio is high, are likely to lock. One has to introduce extra face
degrees of freedom and the restriction is the same as for the case of the linear
Stokes equations. A sufficient condition for avoiding locking in 2D is that each
corner have only three faces without extra degrees of freedom. On the other side,
PEBI grid (also called Voronoi meshes) where this ratio is low will not be likely to
lock. We refer to \cite{da2014mimetic} for a detailed analysis of the necessary
conditions to avoid locking. In the case of the MPSA method, the situation is the
opposite. The discrete representation of the stress tensor is done at each node so
that the solenoidality condition will be imposed there. Therefore, the ratio between
the number of interaction regions (which is also equal to the number of nodes) with
the number of cells (which corresponds to the number of degree of freedom for MPSA)
will determine the sensitivity of the grid to locking. A PEBI grid will then much
more likely lead to locking than a triangular grid. More generally, we can conclude
that the grids where the MPSA method and the VE method lock are dual grids of each
other (not true for quadrilateral). As proven in \cite{NorbottenStable2015}, a
practical advantage of the MPSA method is that when it is coupled with a finite
volume discretization, which is the most common choice of discretization for the flow
equation, the method will be stable independently of the time-step size, even in the
limit of incompressible fluid. In the case of geological applications, the
compressibility of water is about the same as of the rock, which means that locking
is not happening. However, for mud and shale, it may be important.

It seems that the limitations are a bit less severe for MPSA although care has to be
taken in order to require the local inversion of matrices to be sufficiently accurate
so that it does not perturb the div-free part of the solution, since this part will
be multiplied with a large parameter. Finite-volume based method for flow defines a
natural divergence operator into cells. Moreover, the coupling term of the mechanical
system with a finite volume method is due to the volume changes of the cells, which
means that it will require for the mechanical system a divergence operator into
cells. Hence, in the limit of incompressible fluid and small time steps, it will lead
to the same constraint as the near-incompressibility constraint for the mechanical
system. This means that for MPSA the ratio between the number of degrees of freedom
and the near div-free condition, imposed by small time-steps in the Biot's
formulation, is independent of the grid, while for VEM it depends on the ratio
between nodes and cells. In the discretization of the Biot's equations in the
framework of MPSA, we naturally introduce a regularization by making the numerical
divergence for displacement depends on pressure. This can be seen from the
discretization in the reduction to cell-centered variables. The essential ingredients
are that the divergence operator is defined in terms of displacement at cell
boundaries and that the continuity of forces is required by using the effective
stress, that is $\sigma - \alpha p \id$, and not the continuity of the forces due to
stress with an additional force from the pressure gradient. This requires that the
discretization of the mechanics and the pressure system is done together. We also
note that the discretization of the coupling is independent of the discretization of
the gradient in the Darcy equation.

\subsubsection{Regularization methods for the near-incompressible limit}
\label{subsubsec:lockingstrateg}

We have implemented different regularization strategies to handle materials close to
the near incompressible limit. For VEM, our approach follows \cite{da2013virtual}
even if the results there hold for elements of order $k\geq 2$ while we only consider
linear elements, that is $k=1$. We will comment on that later. We use the
constitutive equation given by \eqref{eq:isosigma} and the energy in a cell $K$ is
given by
\begin{equation}
  \label{eq:akhook}
  \frac12a_K(\vec{u}, \vec{u}) =  \mu \int_{K} \epsi(\vec{u}): \epsi(\vec{u})\,dx +
  \frac{\lambda}2\int_K\abs{\dive\vec{u}}^2\,dx.
\end{equation}
We introduce
\begin{equation*}
  a_{\mu,K}(\vec{u}, \vec{v}) = \int_{K} \epsi(\vec{u}): \epsi(\vec{u})\,dx
\end{equation*}
so that \eqref{eq:akhook} can be rewritten as
\begin{equation}
  \label{eq:akhook2}
  a_K(\vec{u}, \vec{u}) =  2\mu\,a_{\mu,K}(\vec{u}, \vec{u}) +
  \lambda\int_K\abs{\dive\vec{u}}^2\,dx
\end{equation}
The coercivity of $a_K$ follows from the coercivity of $a_{\mu,K}$ but it
deteriorates when $\lambda$ get very large compared to $\mu$. In terms of the
Poisson ratio $\nu$, we have
\begin{equation*}
  \frac{\mu}{\lambda} = \frac{1 - 2\nu}{2\nu}
\end{equation*}
so that the deterioration of the coercivity occurs when $\nu$ tends to $\frac12$. In
this case, the exact solution will be very close to a solenoidal field. As mentioned
in section \ref{subsubsec:locking}, numerical locking occurs when too many degrees of
freedom are used to satisfy the solenoidal constraint so that too few are left to
approximate close-to-solenoidal solution. Instead of \eqref{eq:akhook}, let us
consider the following approximation of $a_{K}$,
\begin{equation}
  \label{eq:akhookapprox}
  a_{K,\text{app}}(\vec{u}, \vec{v}) =  2\mu\, a_{\mu,K}(\vec{u}, \vec{v}) +
  \lambda\int_K\abs{\Pi_{0,K}(\dive\vec{u})}^2\,dx.
\end{equation}
Here $\Pi_{0,K}$ is the $L^2$ projection. When $\lambda$ becomes very large, the
strong penalization of the term following $\lambda$ in \eqref{eq:akhookapprox} will
impose on the solution the constraint
\begin{equation}
  \label{eq:piokdivconst}
  \Pi_{0,K}(\dive\vec{u}) = 0
\end{equation}
while, for \eqref{eq:akhook}, it gave $\dive\vec{u} = 0$. We have in this way relaxed
the system as the constraint \eqref{eq:piokdivconst} is easier to fulfill than the
solenoidal constraint $\dive\vec{u}=0$.  More degrees of freedom are therefore left
to resolve the rest of the displacement field. At the same time, we commit a
\textit{variational crime} meaning that we base our formulation on a non-exact form
of the energy. However, in the VE method, the energy we consider is already an
approximation because of the stabilization term and we are going to see that, for the
relaxed version, we retain exactness for linear displacement. To approximate $a_{\mu,
  K}$, we use the projection $\Pcal$ and introduce a stabilization term as described
in Section \ref{subsec:vem}, that is,
\begin{equation}
  \label{eq:dismuenervem}
  a_{\mu, K}^h(\vec{u}, \vec{v}) = a_{\mu, K}(\Pcal\vec{u}, \Pcal\vec{v}) + s_{\mu, K}((\id - \Pcal)\vec{u}, (\id - \Pcal)\vec{v}).
\end{equation}
Then, the discrete approximation of the energy is given by
\begin{equation}
  \label{eq:disstabenervem}
  a_{K}^h(\vec{u}, \vec{u}) = 2\mu\, a_{\mu, K}^h(\vec{u}, \vec{u}) +   \lambda\int_K\abs{\Pi_0\dive\vec{u}}^2\,dx.
\end{equation}
As usual the total energy will obtained by summing up the cell contributions,
\begin{equation}
  \label{eq:totdiscenergreg}
  a^h(\vec{u}, \vec{v}) = \sum_{K\in\Tcal} a_{K}^{h}(\vec{u}, \vec{v}).
\end{equation}
We can check that $\Pi_0(\dive\vec{u})$ can be computed exactly for all elements of
the virtual basis. Indeed to compute the $L^2$ projection of $\dive\vec{u}$, we only
need to evaluate its zero moment, that is, the integral of $\dive\vec{u}$. A
straightforward integration by parts give us
\begin{equation}
  \label{eq:zeromomdiveu}
  \int_K \dive\vec{u}\,dx = \int_{\partial K} \vec{u}\cdot\vec{n}\,dx
\end{equation}
and, by construction, for any $\vec{u}$ which belongs to the virtual basis, $\vec{u}$
is linear on the edges so that the integral on the right-hand side above can be
computed exactly. Thus, the bilinear form $a_K^h$ in \eqref{eq:disstabenervem} can be
assembled and the corresponding system inherits the consistency property of the VE
method. Namely, if $\vec{u}$ is linear and $\vec{v}$ is one of the virtual basis
element, then
\begin{equation}
  \label{eq:constrelax}
  a_{K}^h(\vec{u},\vec{v}) = a_{K,\text{app}}(\vec{u},\vec{v}).
\end{equation}
We define a discrete divergence operator from node to cell variables as
\begin{equation*}
  \sum_{K\in\Tcal}\dis{q}_K\ddiv(\dis{u})_K = \sum_{K\in\Tcal}\dis{q}_K\int_K\dive \vec{u}\,dx
\end{equation*}
for all $\dis{q}_K$. Here $\vec{u}$ is the function in the virtual finite element
space corresponding to the nodal displacement coefficients given by $\dis{u}$. We
assemble the matrix $A$ corresponding to the bilinear form $a_{\mu}^h$ in the same
way as in section \ref{subsec:vem}. We obtain that, for any discrete nodal
displacement vector $\texttt{u}$, the discrete bilinear form $a^h$ takes the form
\begin{equation*}
  a^h(\vec{u}, \vec{u}) = \mu\,\dis{u}^TA\dis{u} + \frac{\lambda}2 \abs{\ddiv(\dis{u})}^2.
\end{equation*}
The discrete system of equations is obtained by taking the variation of $a^h$ and we
get
\begin{equation*}
  2\mu A \dis{u} + \lambda\ddiv^T\ddiv(\dis{u}) = \dis{f},
\end{equation*}
for a given right-hand side $\dis{f}$. We can rewrite this system as
\begin{subequations}
  \label{eq:vemrelaximpl}
  \begin{align}
    2\mu A \dis{u} + \dis{grad}(\dis{p}) &= \dis{f}\\
     \ddiv(\dis{u}) - \frac{1}{\lambda}\dis{p}  &= 0,
  \end{align}
\end{subequations}
where $\dis{grad} = \ddiv^T$. Then, $\dis{p}$ can then be interpreted as a
pressure. This strategy where the solenoidal constraint is relaxed using a projection
operator can be successfully applied when considering higher order method virtual
finite element methods. Indeed, in \cite{da2013virtual}, it is shown that for a
method of order $k$ is the projection operator $\Pi_{k - 2}$ is used for relaxation
then the method is unconditionally convergent with respect to the parameter
$\lambda$. Since we consider linear elements, that is $k=1$, such operator is not
available. Therefore, we need extra degree of freedom, see \cite{da2014mimetic},
where it is shown that it is only necessary to introduce an extra edge degree of
freedom for edges which connect to nodes that have more than three edges. The
following three VE methods have been implemented,
\begin{enumerate}[topsep = 2mm, , itemsep = 1mm, labelsep = 1mm, labelwidth = 3.2cm,
  itemindent = !, labelindent = !, leftmargin = !, align = parleft]
\item[\textit{VEM}:] Standard implementation, as described in Section \ref{subsec:vem},
\item[\textit{VEM-relax}:] Implementation using the relaxed version coming from
  \eqref{eq:disstabenervem}, see \eqref{eq:vemrelaximpl},
\item[\textit{VEM-relax-extra}:] Same as \textit{VEM relax} but we introduce an extra
  degree of freedom on each face so that the stability condition given in
  \cite{da2014mimetic} is fulfilled.
\end{enumerate}

For the MPSA method, a regularization of a similar nature is presented in
\cite{nordbotten2014cell} in the case of the poro-elastic equation. We detail its
application to the incompressible limit. We use the same framework and notations as
given in Section \ref{subsec:mpsaexpl}. First, we add to each cell $K$ an extra
degree of freedom $p_K$ to approximate the divergence term $\lambda_K\dive u$ in the
cell $K$. We replace the definition \eqref{eq:defTkssigma} of the forces on sub-faces
as
\begin{equation}
  \label{eq:defTkssigmaincomp}
  T_{K,s}^{\sigma} = m_\sigma^{s}\left(2\mu_{K}\gradb_{K, s}^{\sym}(u) + p_K I\right)\cdot n_{K,\sigma}
\end{equation}
The sequence of operations given by \eqref{eq:Tkssigmasimp} is essentially the same
except that it uses $p_K$ in the last step,
\begin{equation}
  \label{eq:Tkssigmasimpincomp}
   \left\{u_{K'}\right\}_{K'\in\Tcal_s} \overset{I_{FV,s}}{\longrightarrow}
  \begin{array}[h]{l}
    \left\{u_{K}\right\}_{K\in\Tcal_s},\\ \left\{u_{K,s}^{\sigma,\beta}\right\}_{\substack{\sigma\in\Fcal_s,K\in\Tcal_s\\\beta=\{1,2\}}}
  \end{array}
  \overset{\gradb_{K,s}}{\longrightarrow} (\gradb u)_{K,s} \underset{\text{using }p_K}{\overset{\text{by \eqref{eq:defTkssigmaincomp}}}{\longrightarrow}}
  T_{K,s}^{\sigma}.
\end{equation}
The linear mapping $I_{FV, s}$ is defined using the same principle as before: Given
$u_{K'}$ and $p_{K'}$ for $K'$ in the interaction region $\Tcal_{s}$, choose the
coefficients of $I_{FV, s}$ such that the forces are continuous at each sub-faces,
that is \eqref{eq:contforce} holds, and the measure of the jumps in displacement
values given by \eqref{eq:leastsqpb} is minimized. To summarize, using local
reduction, we obtain at each interaction regions $\Tcal_s$,
\begin{equation*}
  \left\{u_{K}, p_{K}\right\}_{K\in\Tcal_s} \longrightarrow \left\{T_{K,s}^{\sigma}, u_{K,s}^{\sigma} \right\}_{K\in\Tcal_s, \sigma\in \Fcal_{K,s}}
\end{equation*}
The global system of equation is then given by the equation of conservation of
momentum \eqref{eq:consmomentum} and the
following equation for the pressure
\begin{equation}
  \label{eq:defpkincomp}
  p_K = \frac{\lambda_K}{|K|}\sum_{\{s| K\in\Tcal_s\}}\sum_{\{\sigma\in\Fcal_K\cap\Fcal_s\}} (|(K,s,\sigma)|u_{K,s}^\sigma\cdot n_{K,\sigma}),
\end{equation}
where $|(K,s,\sigma)|$ denotes the length (or surface in 3D) of the the sub-face
$(K,s,\sigma)$ and $|K|$ the volume of the cell $K$. Equation \eqref{eq:defpkincomp}
is the discrete counterpart of the identity
\begin{equation*}
  \int_{K}p\,dx = \lambda_K\int_{K}\dive u\,dx= \sum_{\sigma\in\Fcal_{K}}\int_\sigma u\cdot n_{K,s}\,d\sigma.
\end{equation*}
Using the same arguments as in \cite{nordbotten2014cell}, one can prove that with
essentially the same grid restrictions that apply for the elastic and pressure
discretizations independently, the method is convergent uniformly with respect to
$\lambda$. The method introduces the extra degrees of freedom $p_K$ and it also
introduces relaxation. Indeed, the divergence term in the definition of the force is
imposed through $p_K$, that is, from the condition \eqref{eq:defpkincomp}, which is
imposed cell-wise. This represents a relaxation in comparison with the original MPSA
method, where different values of the divergence are used for each sub-interaction
region (of the type $(K,s)$). The following two MPSA methods have been implemented
\begin{enumerate}[topsep = 2mm, , itemsep = 1mm, labelsep = 1mm, labelwidth = 3.2cm,
  itemindent = !, labelindent = !, leftmargin = !, align = parleft]
\item[\textit{MPSA}:] Standard implementation, as described in Section \ref{subsec:mpsaexpl},
\item[\textit{MPSA-relax-extra}:] Implementation using the relaxed version where an extra
  pressure degree of freedom is used, see \eqref{eq:defTkssigmaincomp}.
\end{enumerate}

\section{Numerical test cases}
\label{sec:numericaltest}
The test cases are designed to study the robustness of the methods with respect to
grid features that are specific to subsurface applications. All of the code has been
written and run using the framework of MRST, \cite{MRST:2015bb}. We consider only
two-dimensional configurations and plan to study three-dimensional configurations in
subsequent works. At the moment, only full Dirichlet boundary conditions have been
implemented for MPSA but the extension to other general boundary conditions (rolling
conditions, that is, component-wise Dirichlet conditions) is not difficult but
requires some careful work. We summarize the different test cases in Table
\ref{tab:testcases}.

We pay particular attention to the error in the divergence field, because of its
central role in the coupling with poro-elasticity, and to the local behavior of the
stress fields, due to its importance in the simulation of the development of
fractures and faults. When comparing the error estimates that are presented below, it
is important to understand how the discrete $L^2$ and $L^\infty$ norms are computed
for both methods. The displacement values are defined on nodes for VEM and at cell
centers for MPSA, so that the discrete $L^\infty$ norms for the displacement, even if
not equivalenţ provide comparable estimates. For both methods, the divergence is
defined on cells, so that the discrete $L^2$ and $L^\infty$ norms are directly
comparable. The stress is piecewise constant for VEM, while MPSA defines forces on
faces. We define the discrete $L^2$-norm for stress for the MPSA method as the
summation of the discrete $L^2$ norms of the force over the edges. In this way, we
obtain an averaged quantity but it is not directly comparable with the discrete $L^2$
norm used for stress in VEM, which is the standard volume integral over the domain.

In all test cases, we use the same reference solution which is computed as follows.
We consider the displacement field $\vec{u}=[u_1, u_2]^T$ given by
\begin{align}
  \label{eq:exact_disp_field}
  u_1 & = x(1 - x)\sin(2\pi y)& u_2 & = \sin(2\pi x)\sin(2\pi y),
\end{align}
for $x$ and $y$ belonging to $[0,1]$. Using \eqref{eq:lin_elast_cont}, we compute the
force $\vec{f}$ for which $\vec{u}$, given by \eqref{eq:exact_disp_field} is the
solution. In this way, we have obtained an exact solution of
\eqref{eq:lin_elast_cont} that we will use to compare our results in all the examples
below. The boundary conditions are zero displacement on all sides. A summary of all
the numerical tests that are presented is included at the end of the article, in
Table \ref{tab:testcases}.

\subsection{Case \ref{case:fullcart}: Twisted grid with random perturbation}
\refstepcounter{mycase}\label{case:fullcart}

\begin{figure}[h]
  \centering
  \includegraphics[width=0.6\textwidth]{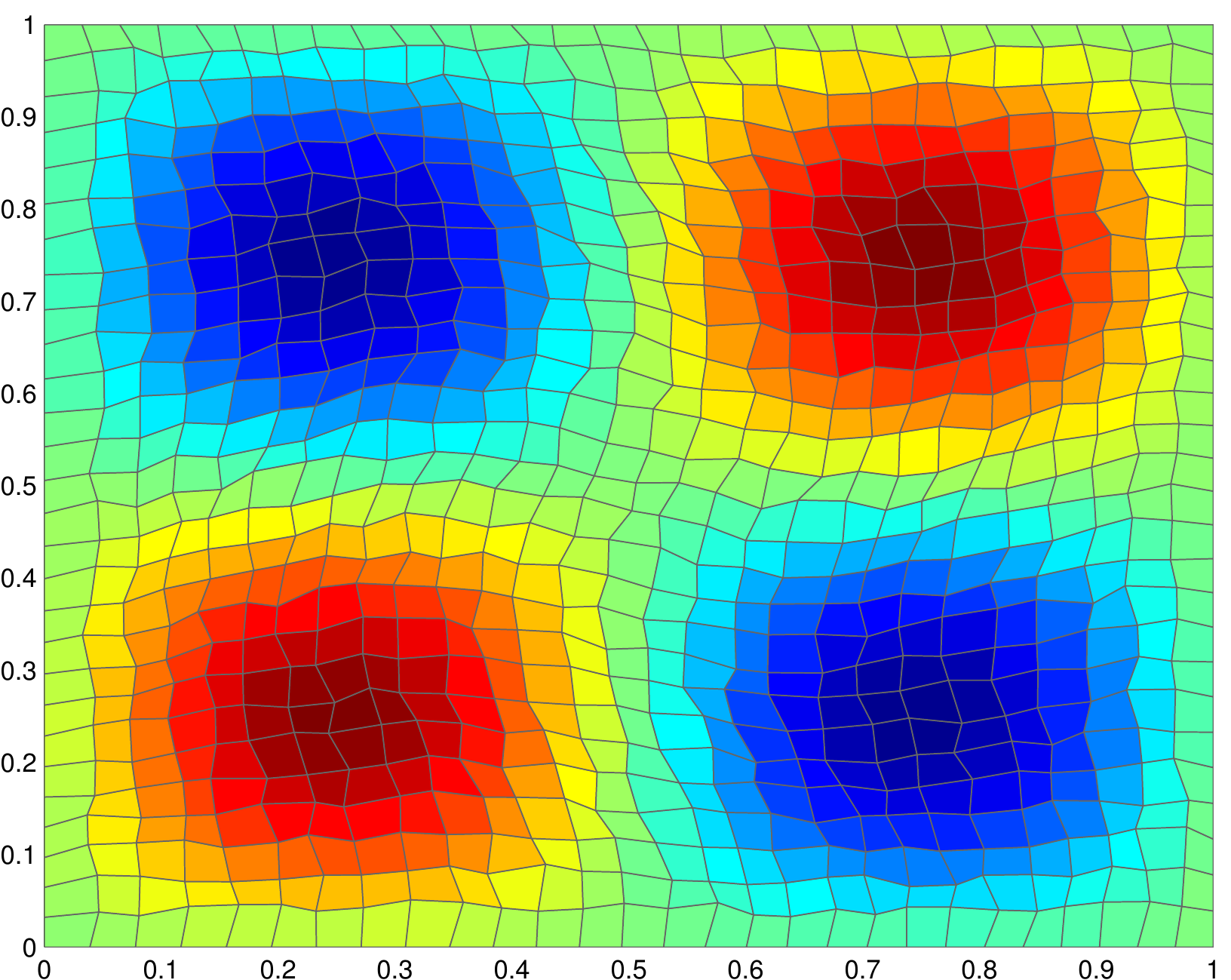}
  \caption{Example of a deformed grid, with refinement factor equal to 7 (The
    colors represents the displacement field in the $y$-direction computed using the
    MSPA method)}
  \label{fig:exdefgrid}
\end{figure}

\begin{figure}[h]
  \label{fig:conf_rough}
  \includegraphics[width=\textwidth]{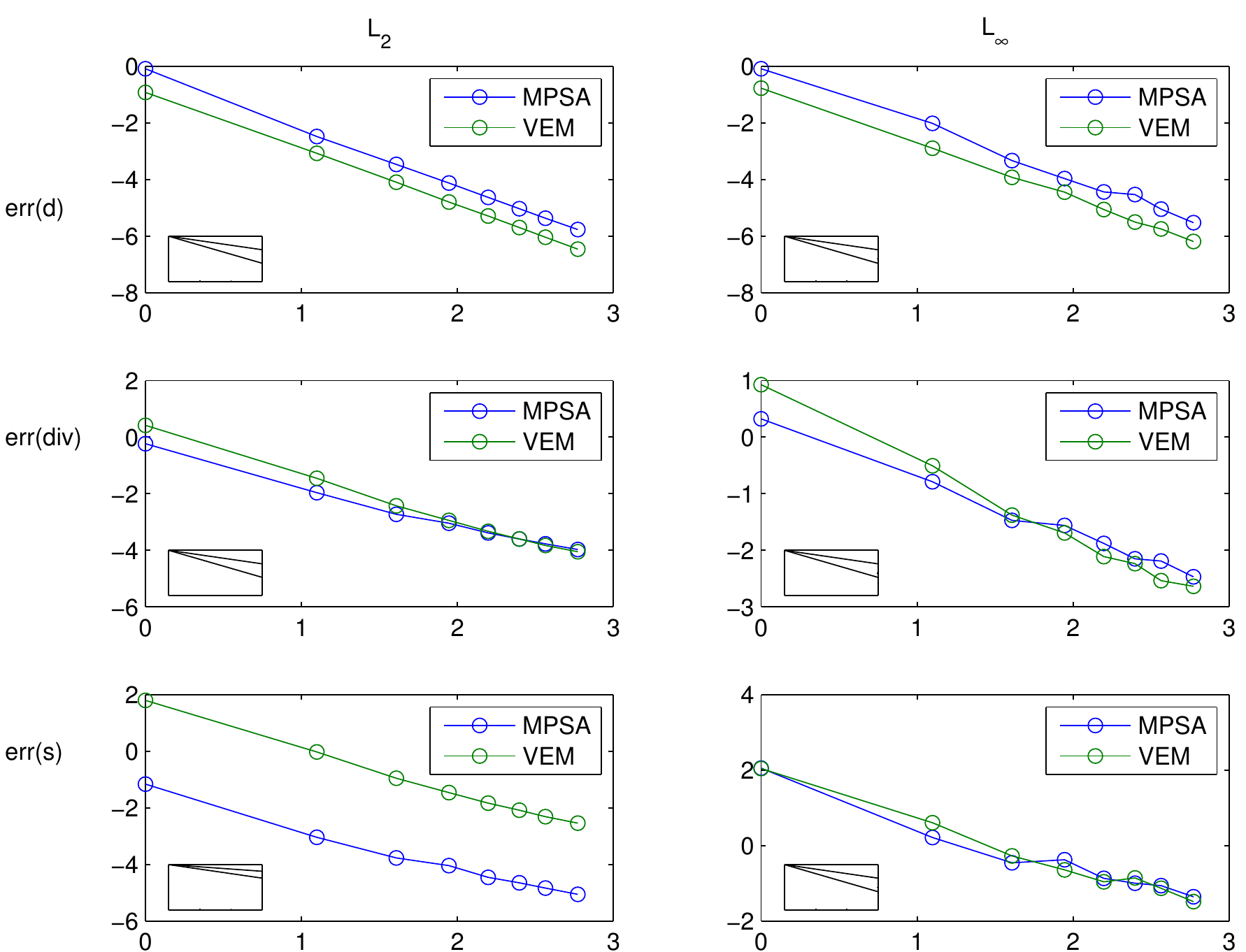}
  \caption{Convergence plot for a twisted Cartesian grid. The $L^2$-norm (left) and
    $L^\infty$-norm (right) of the error are plotted for the displacement (upper
    row), the divergence (middle row) and the stress (lower row). Logarithmic scales
    are used and the values on the x-axis give the logarithms of the refinement
    factor. In the small boxes, the slopes corresponding to a convergence factor for
    one and two are represented.}
\end{figure}

In this test case, we study the convergence properties of both methods. The VE method
is in general first order convergent, as shown in \cite{da2013virtual,da2014mimetic},
but numerical tests show second order convergence under general conditions
\cite{da2013virtual,gain2014}. For the VE method used in the near-incompressible
case, see Section \ref{sec:stabnearincomp}, where pressure is considered as an
independent variable, then the pressure converges at first order, see \cite[theorem
9.1]{da2014mimetic}. Convergence estimates for the MPSA method are not available in
the established literature but numerical tests show the same features as VEM, see
\cite{nordbotten2014cell}. Accordingly, Figure \ref{fig:conf_rough} shows convergence
rates of two and one for the displacement and the divergence, respectively. The most
refined grid is obtained by refining 16 times the initial grid. The grids which are
considered in this test case are \textit{non-regular quadrilateral} grids, see Figure
\ref{fig:exdefgrid} for an example. To generate such grids, the starting point is a
uniform Cartesian grid with a given refinement. Then, a deformation field which is
independent of the refinement factor and which lets the boundaries invariant is
applied to each node.  However such approach leads to the generation of cells that
are close to parallelogram for small refinement, meaning that the grid is strongly
regularized in the refinement process. Such regularity for the grid cannot be
expected in a realistic context and that is why add a random perturbation to each
node. To cancel out the random part in the generation of the grids, we have produced
the same error plots several times and we observe that the convergence rates keep the
same characteristics. The $L^\infty$-norm is computed in a straightforward way by
taking the maximum value over all nodes, for the VE method, and over all cell
centers, for the MPSA method. As far as the $L^2$ norm is concerned, in the case of
the VE method, we approximate it using the same quadrature rule as in \cite[Section
3.3]{gain2014}. For the MPSA method, we weight the cell-values by the volume.

\clearpage

\subsection{Case \ref{case:mixedgrid}: Mixed cell types}
\refstepcounter{mycase}\label{case:mixedgrid}

\begin{figure}[h]
  \centering
  \includegraphics[width=0.8\textwidth]{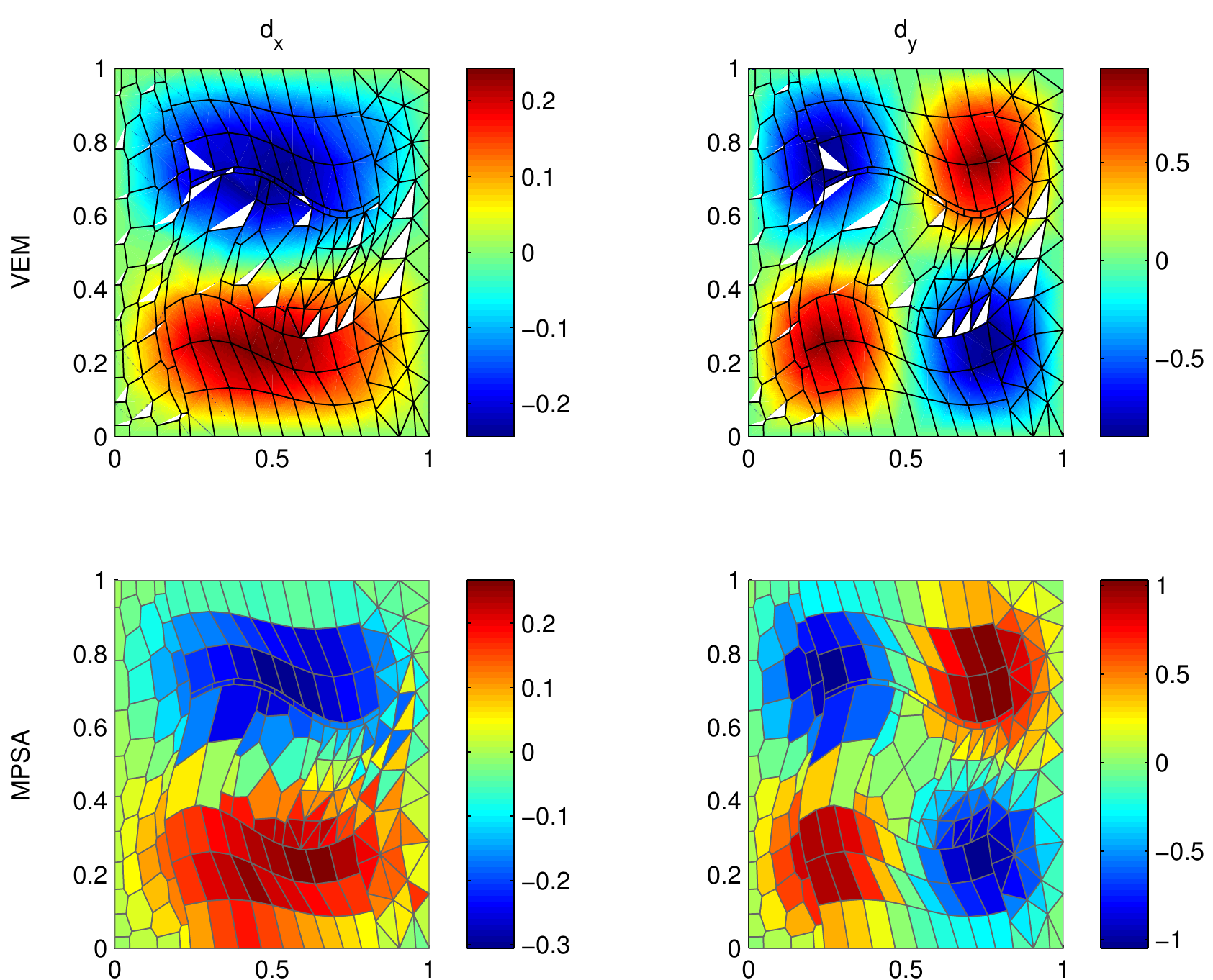}
  \caption{Displacement field $\vec{u}$ (Case \ref{case:mixedgridplain}) for the MPSA and
    VE methods.}
  \label{fig:mixed1}
\end{figure}

\begin{figure}[h]
  \centering
  \includegraphics[width=0.8\textwidth]{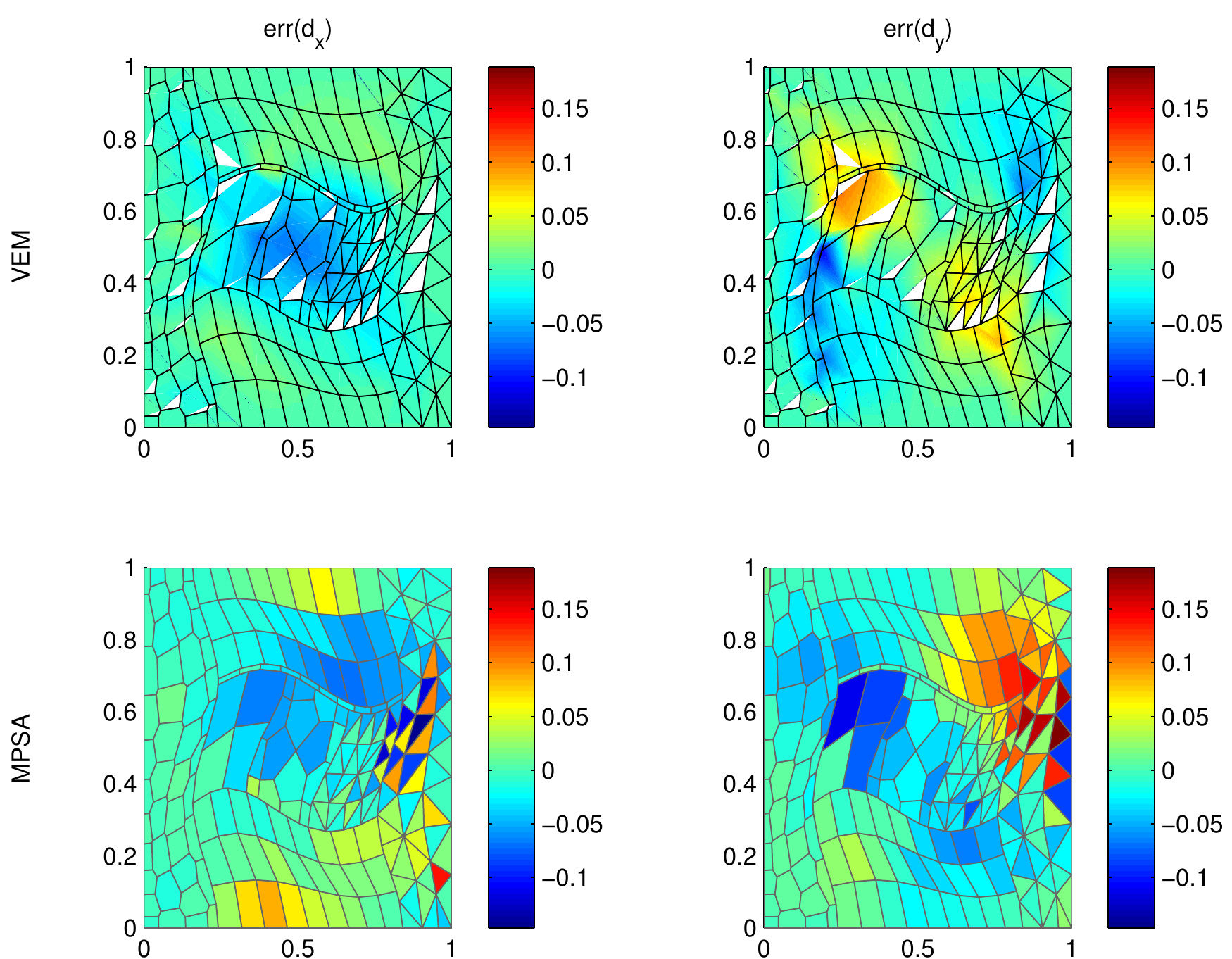}
  \caption{Error in the displacement field $\vec{u}$ for both methods (Case
    \ref{case:mixedgridplain}).}
  \label{fig:mixed2}
\end{figure}

\begin{figure}[h]
  \centering
  \includegraphics[width=0.8\textwidth]{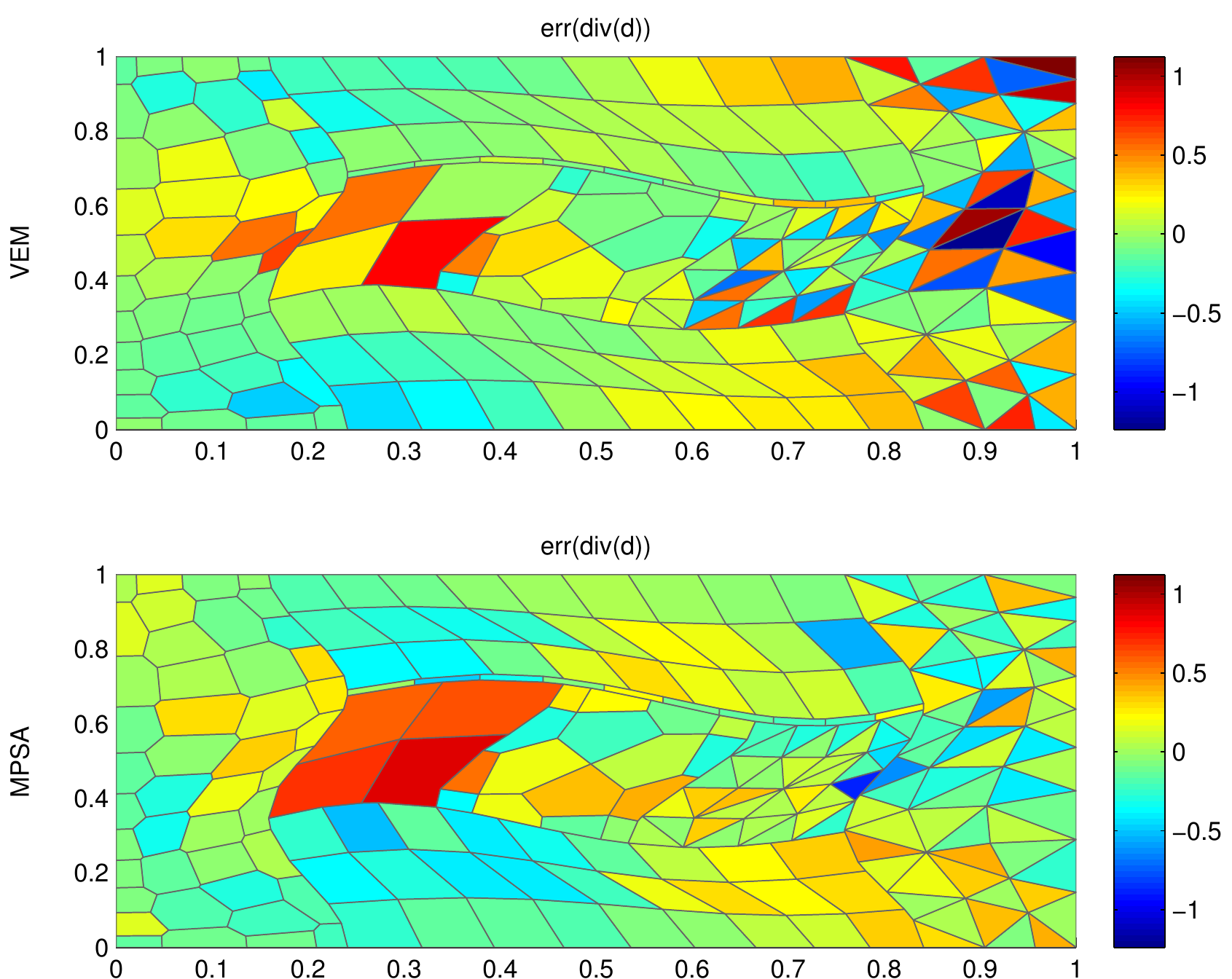}
  \caption{Error in the divergence field $\dive \vec{u}$ for both methods (Case \ref{case:mixedgridplain}).}
  \label{fig:mixed3}
\end{figure}

We set up a case (Case
\ref{case:mixedgridplain}\refstepcounter{mysubcase}\label{case:mixedgridplain}) with a
grid which mixes several difficulties. The grid is made up by, first, assembling
regions with different cell types (triangles, quadrilaterals, hexahedral) and, then,
twisting the grid. Many cells have unfavorable aspect-ratio. There are also hanging
nodes. However, as it can be seen from Figures \ref{fig:mixed1}, \ref{fig:mixed2} and
\ref{fig:mixed3}, the methods manage to capture rather well the exact solution. Note
that MPSA has problem to handle triangles where not all the angles are smaller than
$90$ degrees. Nevertheless, the error that can be made on these cells does not spread
to the whole domain.

We investigate further the case of large aspect ratio for hexahedral (Case
\ref{case:stretchhexa}\refstepcounter{mysubcase}\label{case:stretchhexa}) and
triangular grid (Case
\ref{case:stretchtriangle}\refstepcounter{mysubcase}\label{case:stretchtriangle}). Both
grids are obtained by stretching uniform grids in the $x$-direction with a given
factor. In both cases, we use an aspect ratio of $7$. In Figure \ref{fig:stretchhexa}
and \ref{fig:stretchtriangle}, we observe that the VEM method manages to handle the
hexahedral case correctly while, for the grid made by triangles, it produces reliable
values for the displacement field but oscillatory values for the divergence. As far
as the MPSA method is concerned, we exceed the restrictions on the grid that the
method requires and the method fails, see Figure \ref{fig:stretchhexa}. The case of
the triangle grid is not plotted for the MPSA method.

\begingroup
\def\imagewidth{0.3\textwidth}
\def\insertimage #1{\parbox[c][][c]{\imagewidth}{\centering\includegraphics[width=\imagewidth]{#1}}}
\def\insertcolorbar #1{\parbox[c][][c]{\imagewidth}{\centering\hspace*{1mm}\rotatebox{90}{\includegraphics[height=\imagewidth]{#1}}}}
\def\insertlegend #1{\parbox[c][][c]{\imagewidth}{\centering #1}}
\def\insertcase #1{\parbox[c][][c]{0.13\textwidth}{\centering #1}}

\begin{figure}[h]
  \centering
  \begin{tabular}{@{}c@{}c@{}c@{}c@{}}

    &
    \insertlegend{$d_x$}&
    \insertlegend{$d_y$}&
    \insertlegend{$\dive\vec{u}$}\\

    \insertcase{VEM}&
    \insertimage{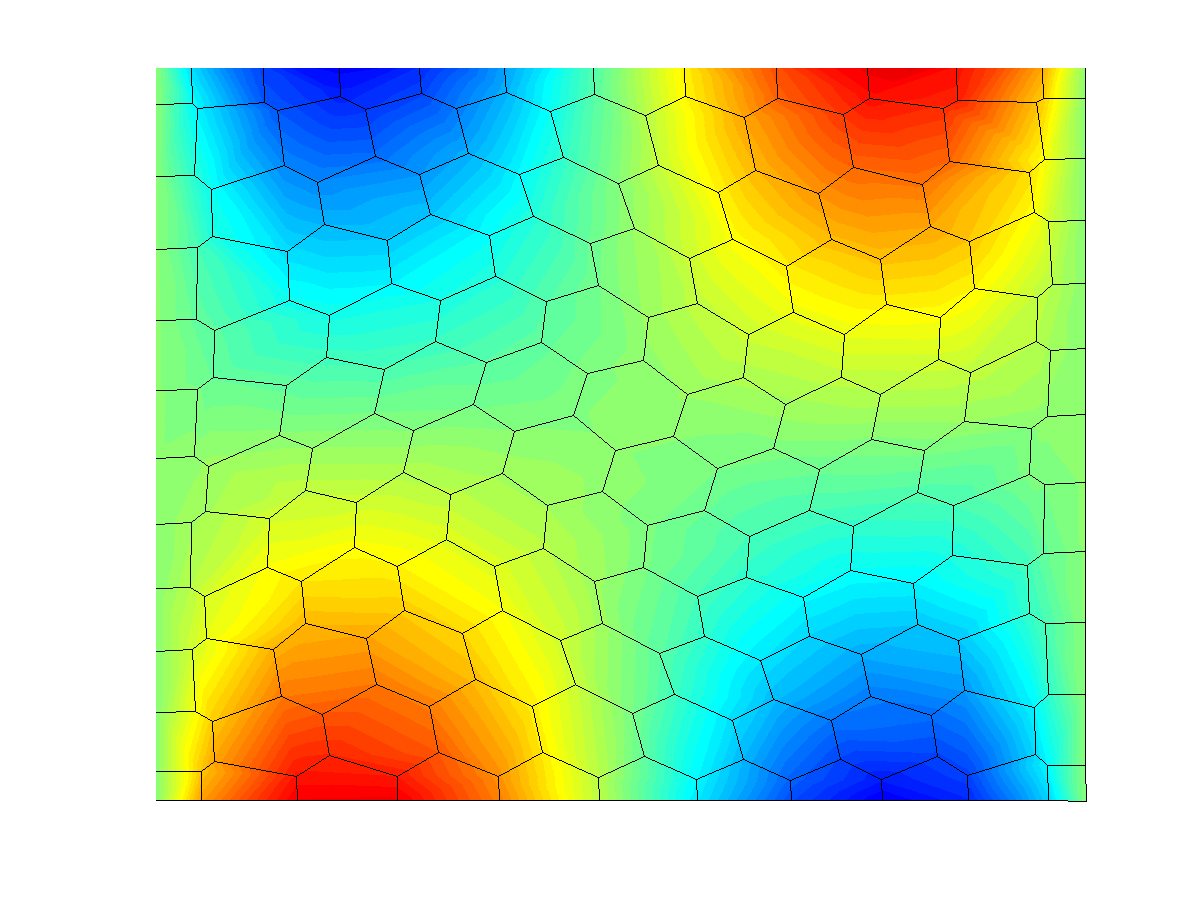}&
    \insertimage{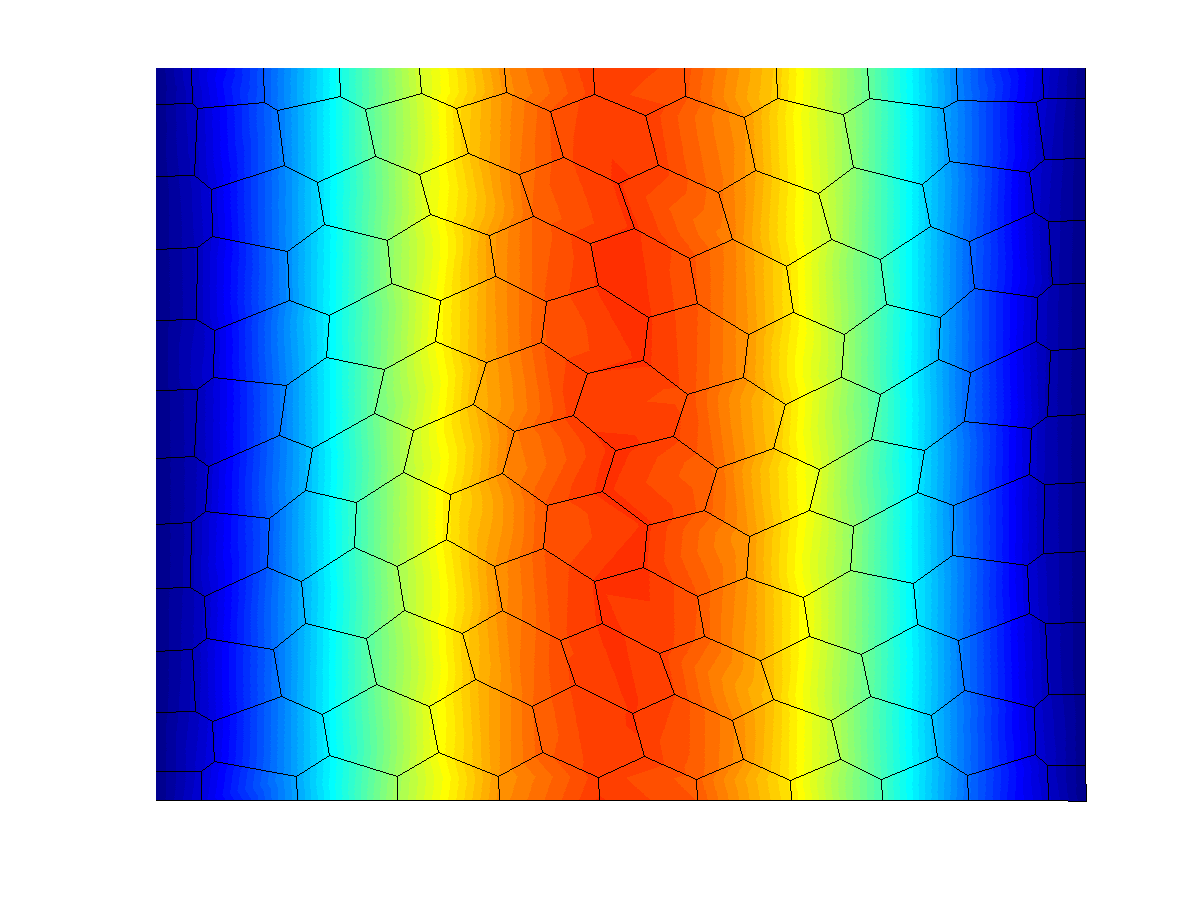}&
    \insertimage{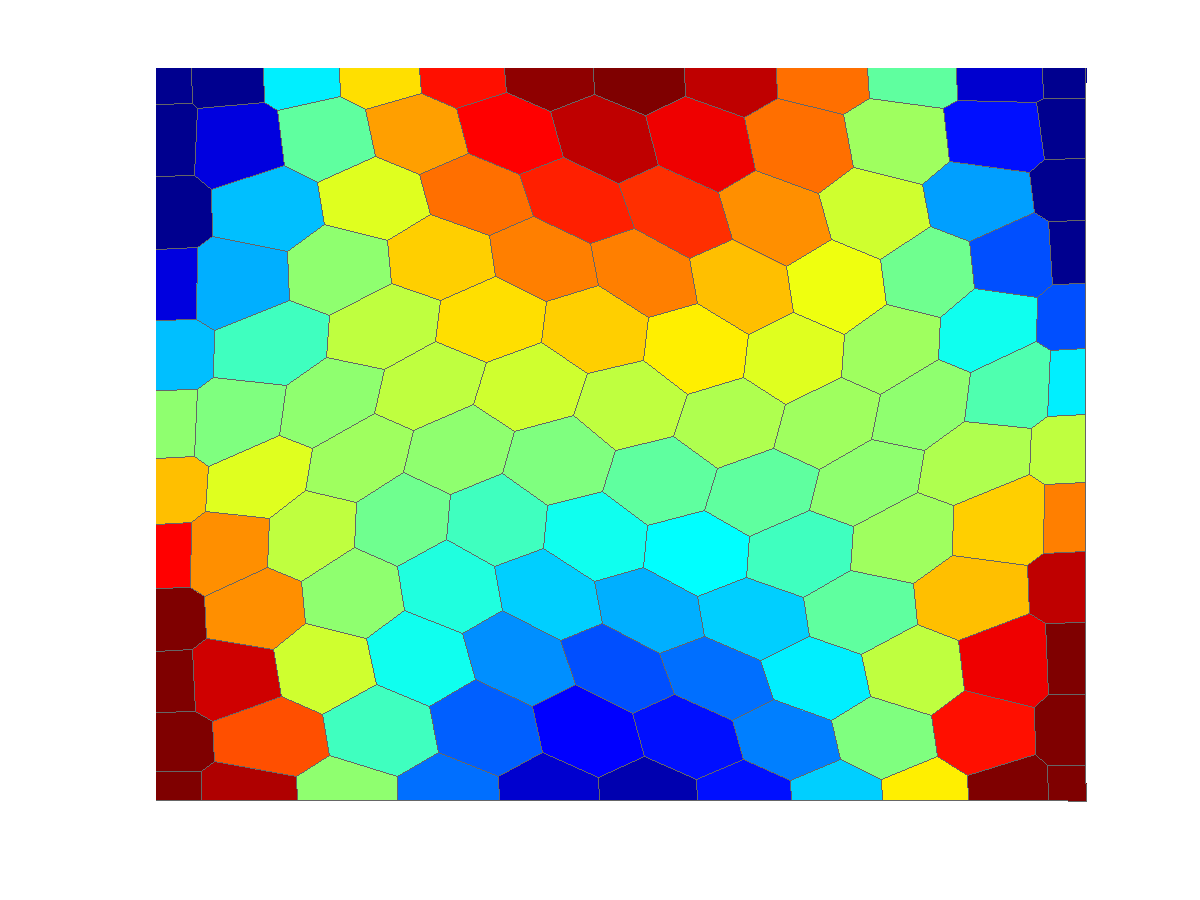}\\

    \insertcase{MPSA}&
    \insertimage{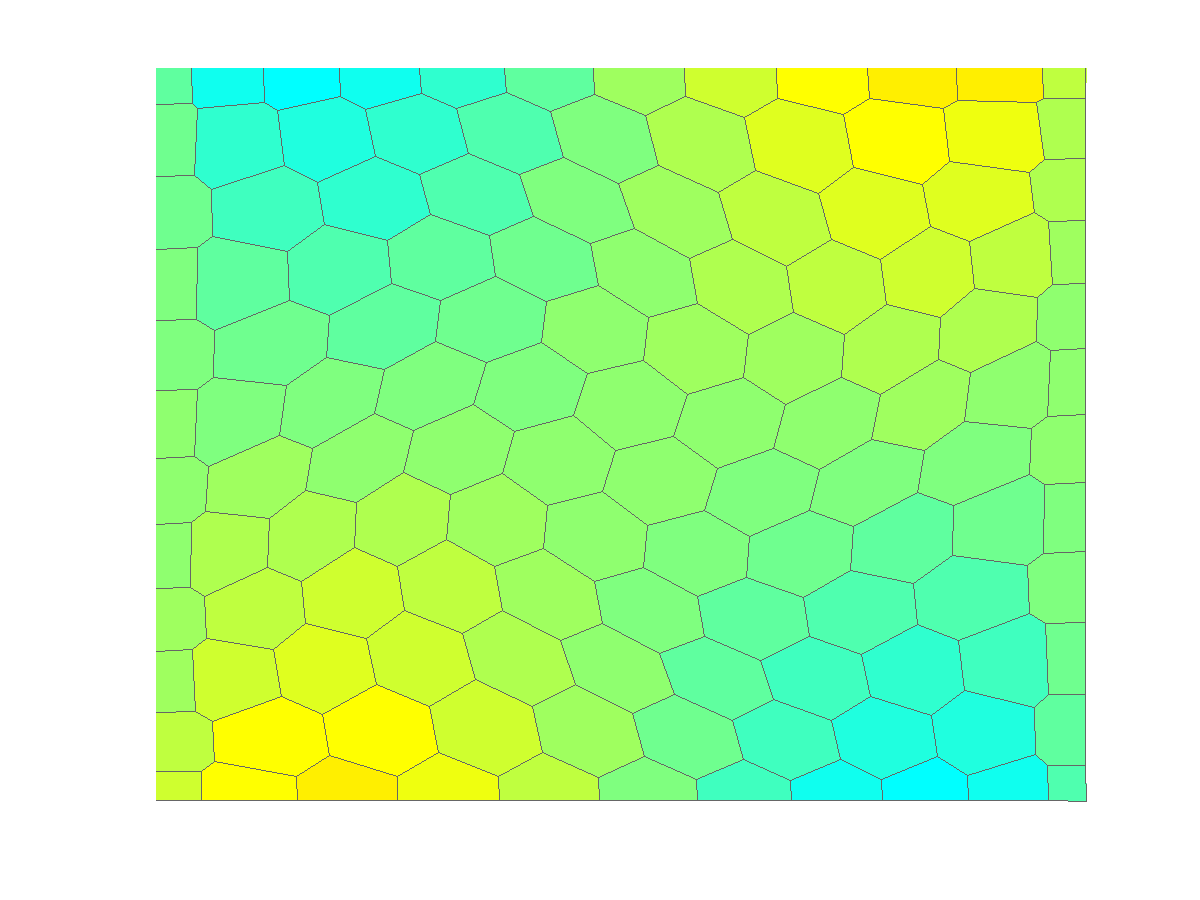}&
    \insertimage{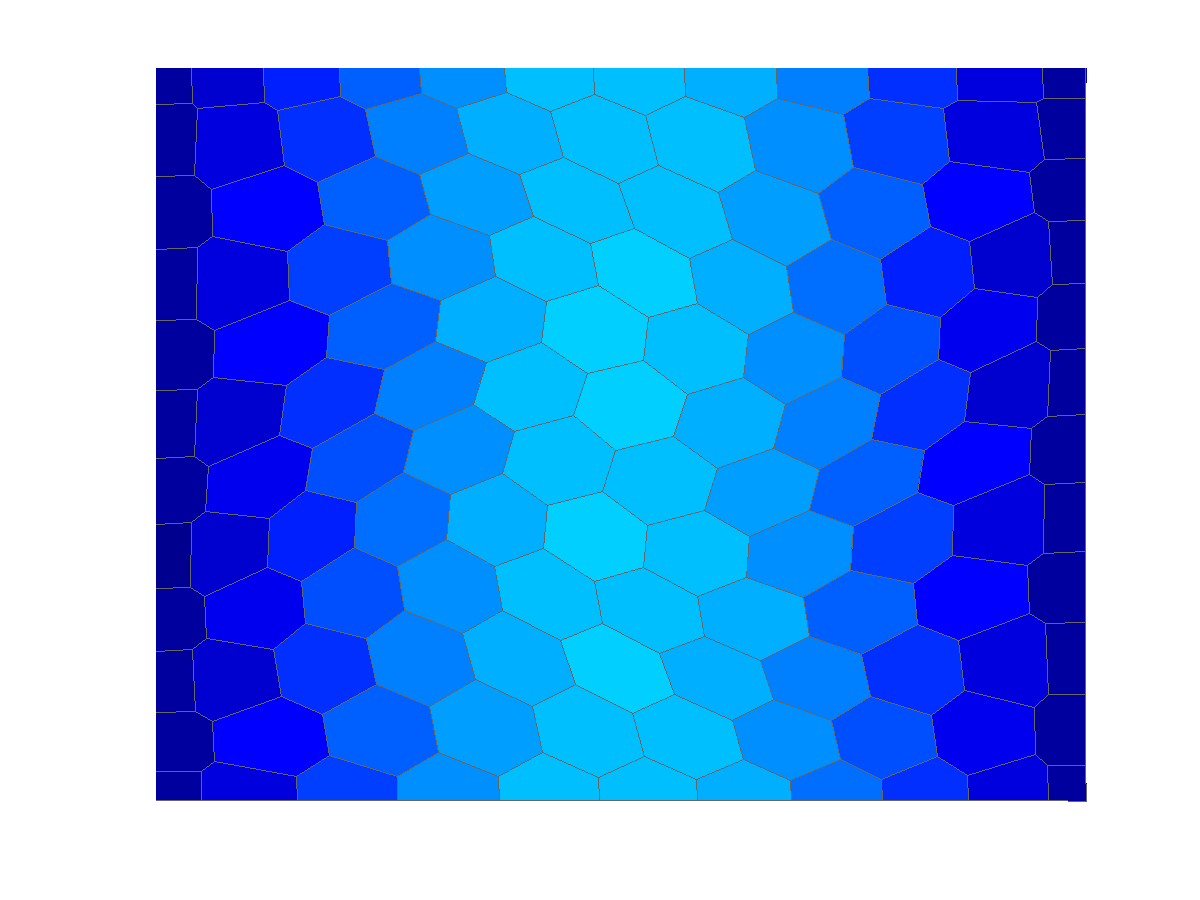}&
    \insertimage{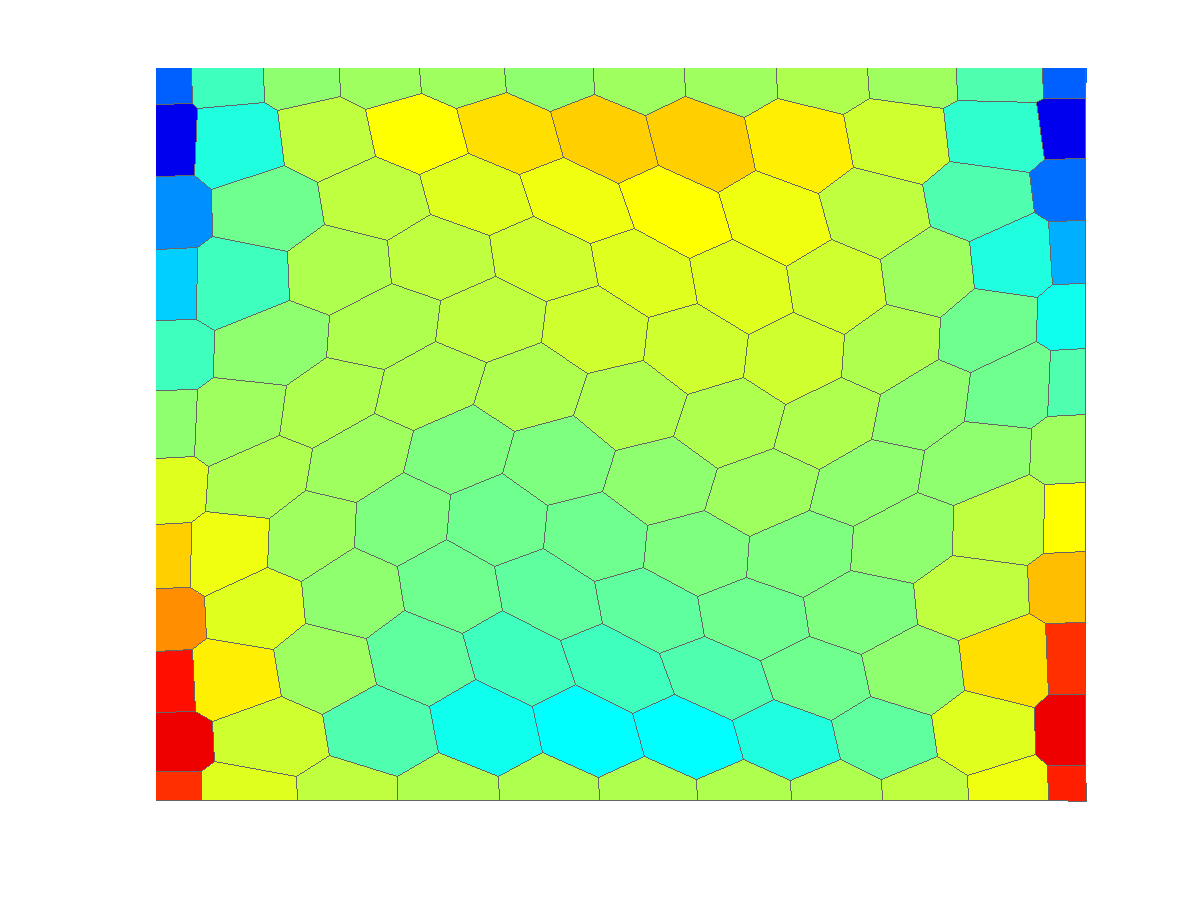}\\

    &
    \insertcolorbar{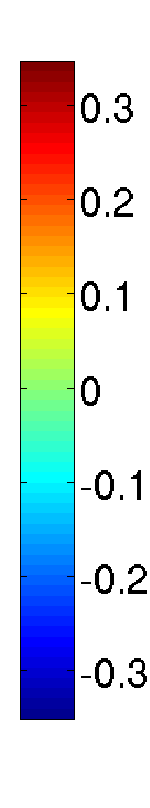}&
    \insertcolorbar{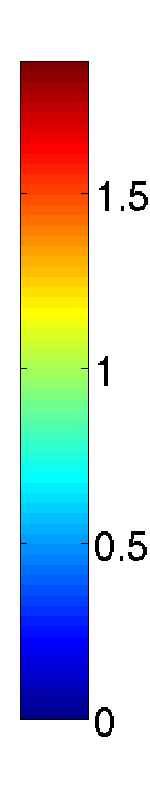}&
    \insertcolorbar{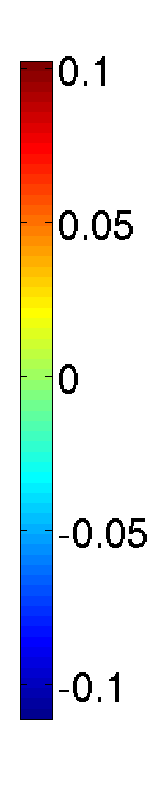}\\

  \end{tabular}
  \caption{Aspect ratio $7$ using an hexahedral grids (Case
    \ref{case:stretchhexa}). For the purpose of a better visualization, we plot the
    values of the displacement and divergence fields on a grid which is stretched
    back to a uniform grid (with aspect ratio $1$).}
  \label{fig:stretchhexa}
\end{figure}

\begin{figure}[h]
  \centering
  \begin{tabular}{@{}c@{}c@{}c@{}c@{}}

    &
    \insertlegend{$d_x$}&
    \insertlegend{$d_y$}&
    \insertlegend{$\dive\vec{u}$}\\

    \insertcase{VEM}&
    \insertimage{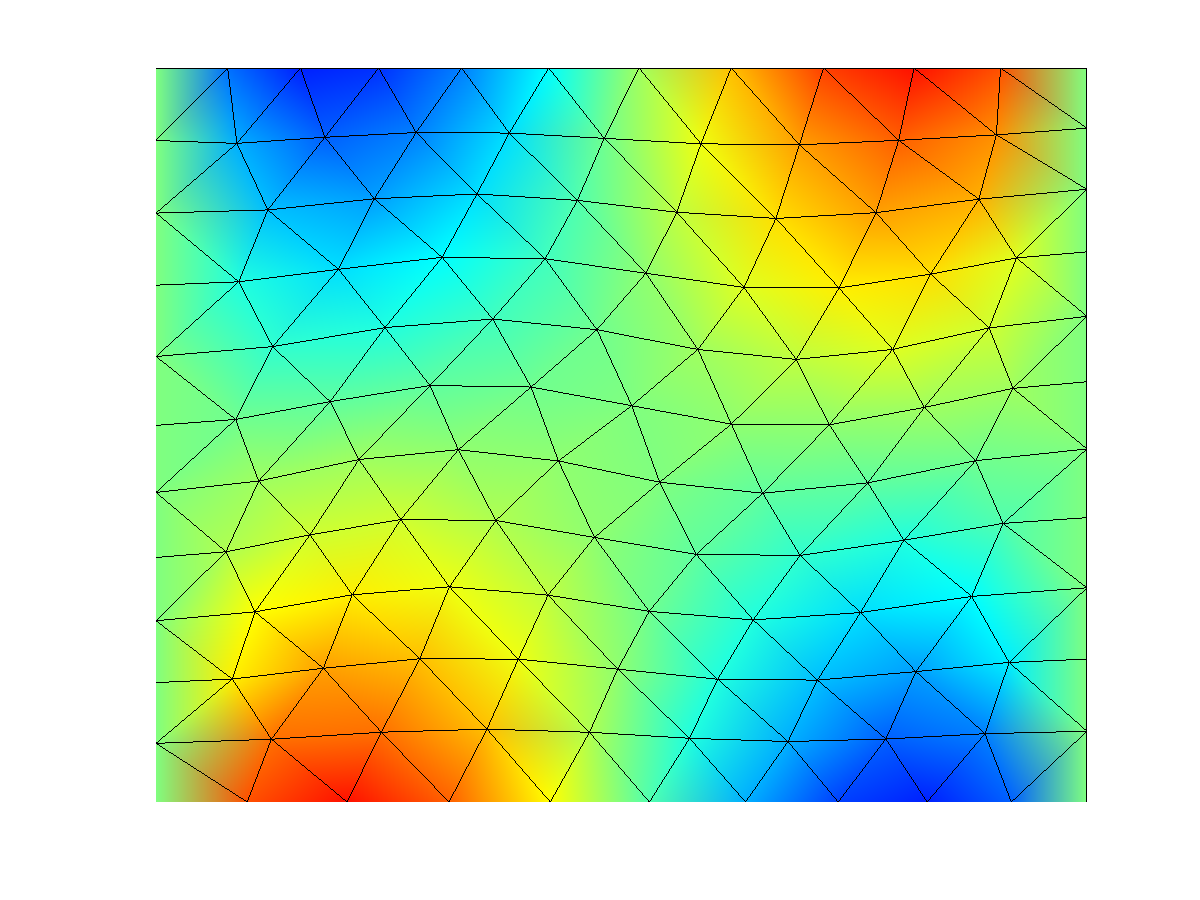}&
    \insertimage{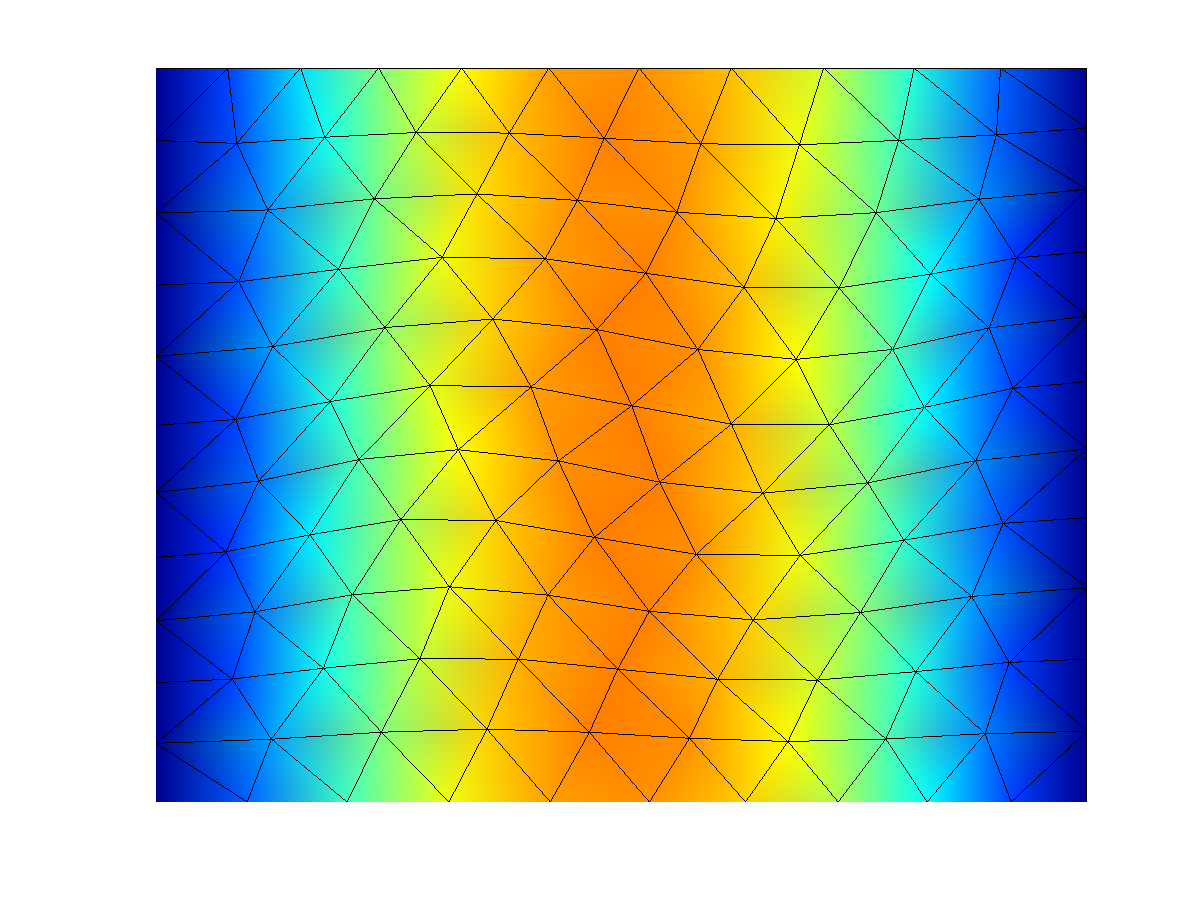}&
    \insertimage{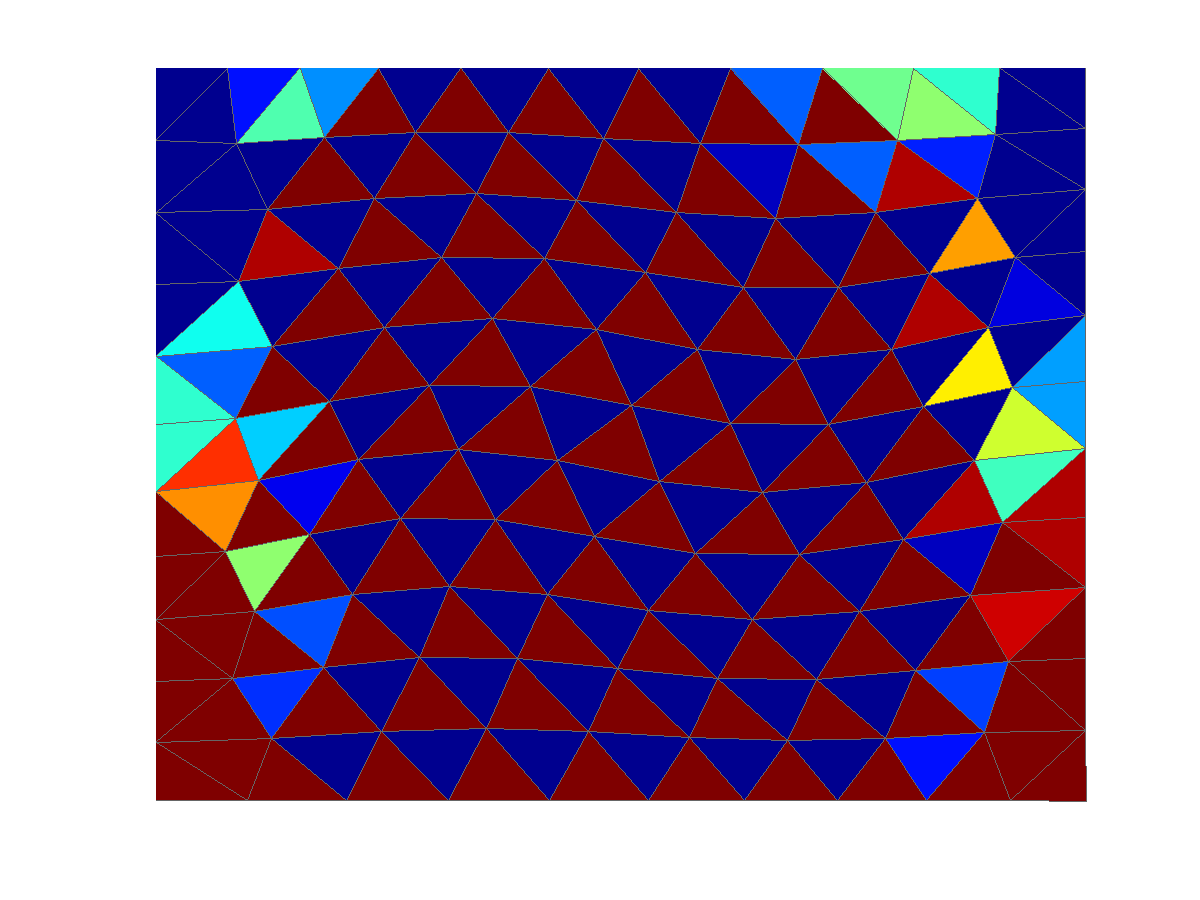}\\

    &
    \insertcolorbar{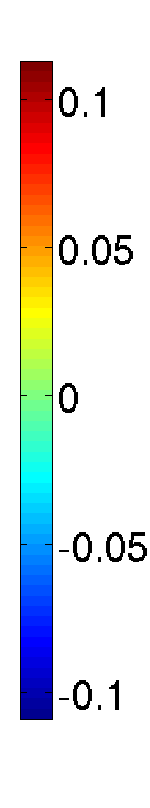}&
    \insertcolorbar{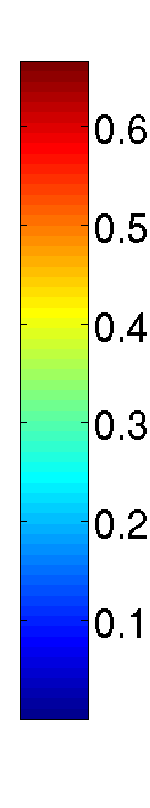}&
    \insertcolorbar{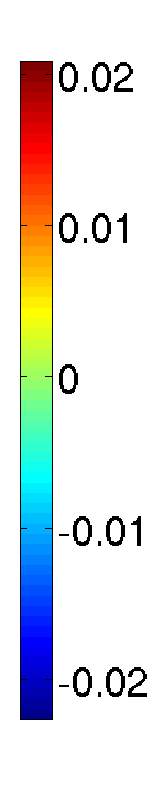}

  \end{tabular}
  \caption{Aspect ratio $7$ using a triangular grid (Case
    \ref{case:stretchtriangle}). The same visualization procedure as in Figure
    \ref{fig:stretchhexa} is used here.}
  \label{fig:stretchtriangle}
\end{figure}
\endgroup

\clearpage

\subsection{Case \ref{case:uniformVerticalRefinement}: Stability for refinement in one direction}
\refstepcounter{mycase}\label{case:uniformVerticalRefinement}

\begin{figure}[h]
  \centering
  \includegraphics[width=0.8\textwidth]{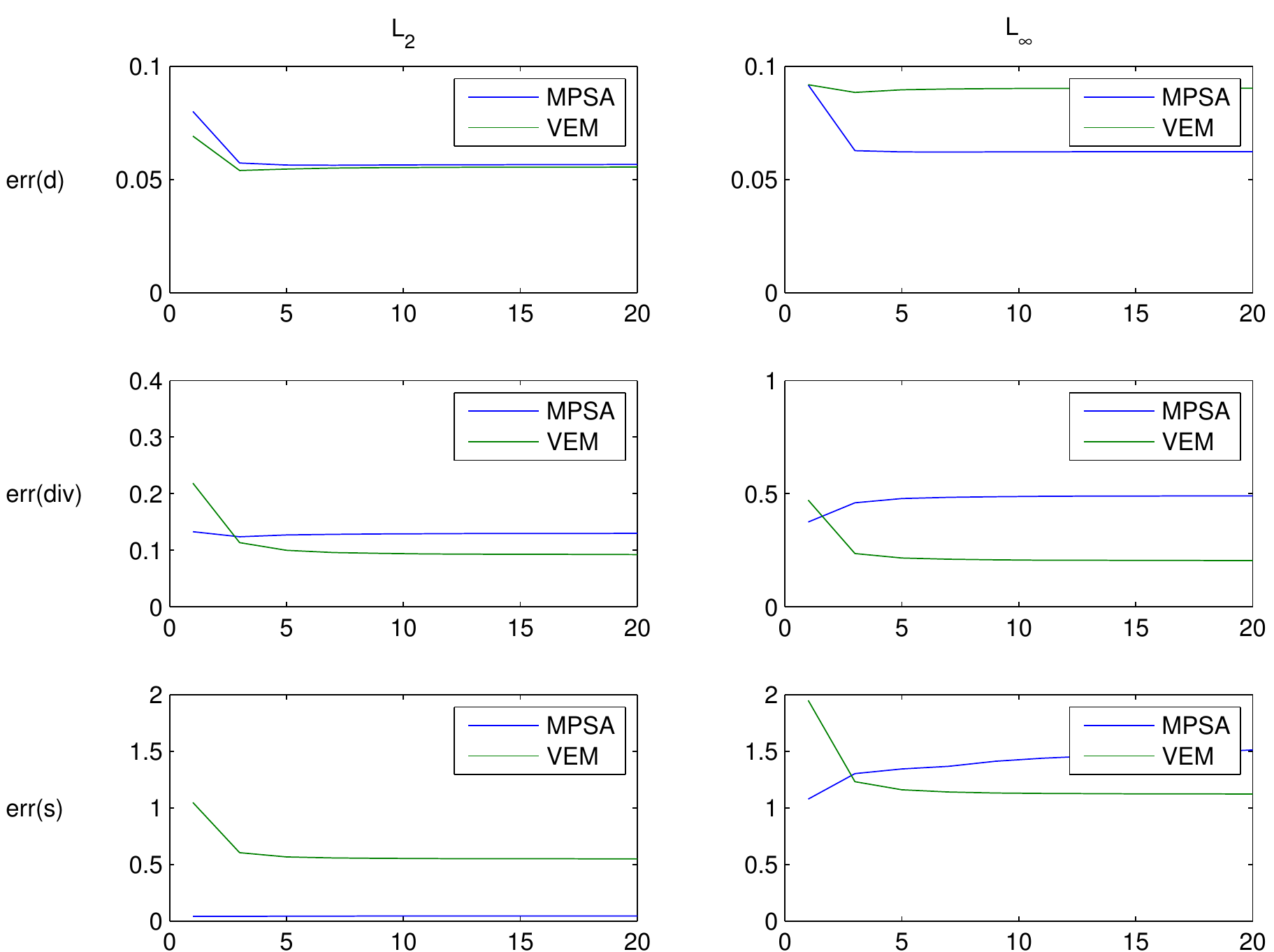}
  \caption{Plot of the errors when the grid is uniformly refined in the $y$-direction
    before being twisted (Case \ref{case:uniformVerticalRefinement}). We consider the
    The $L^2$ norm (left) and $L^\infty$-norm (right) of the errors for the
    displacement (upper row), the divergence (middle row) and the stress (lower
    row). The error remains under control for all methods even when the aspect ratio
    deteriorates. The $x$-axis indicates the refinement ratio in the
    $y$-direction.}
  \label{fig:refineonedirectionconv}
\end{figure}

\begin{figure}[h]
  \centering
  \includegraphics[width=0.8\textwidth]{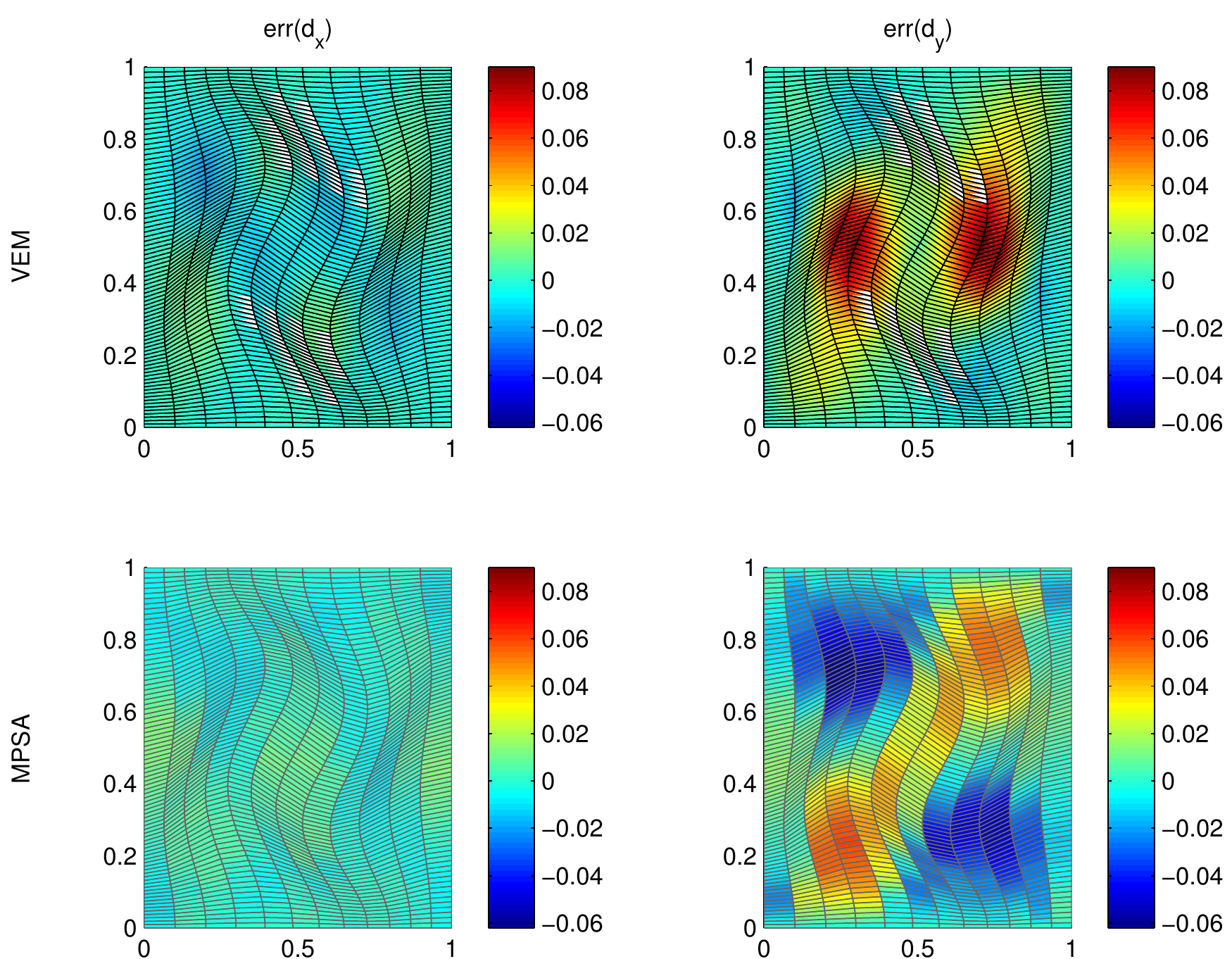}
  \caption{Plot of the displacement in both the vertical and horizontal directions
    for Case \ref{case:uniformVerticalRefinement}. In this figure, the refinement ratio is 10.}
  \label{fig:refineonedirectionerr}
\end{figure}

Grids of stratigraphic subsurface models are often designed with long and flat
cell-blocks which reflect the layered structure of the rock. Such cell-blocks have
deteriorated aspect ratio. We set up an example to test the robustness of the methods
with respect to such deterioration of the grid. We use a Cartesian grid which is
refined in the $y$ direction. Moreover, the grid is twisted to break symmetry effects
which may improve artificially the results. In Figures
\ref{fig:refineonedirectionerr} and \ref{fig:refineonedirectionconv}, we plot the
error. Of course, even if we increase the number of cells, we cannot expect any
improvement of the solution in this case but we can see that the solutions are not
significantly impaired by the refinement and the deterioration of the grid.  The
$L^2$-norms of the error for the stress are substantially different for the two
methods, but we recall that these values are not directly comparable, see the
comments at the beginning of Section \ref{sec:numericaltest}.

\clearpage

\subsection{Case \ref{case:twoRegions}: Stability with respect to decomposition of the grid in regions with
  different refinement}
\refstepcounter{mycase}\label{case:twoRegions}

Grids in subsurface simulations are typically heterogeneous, mixing cells of
different sizes and shapes. We consider two test cases where two regions of equal
size but with different refinements are set side by side. For the first case (Case
\ref{case:twoRegionsIsotropicRefinement}\refstepcounter{mysubcase}\label{case:twoRegionsIsotropicRefinement}),
the refinement in the region on the right-hand side is done in both $x$ and $y$
direction. For the second case (Case
\ref{case:twoRegionsVerticalDirectionRefinement}\refstepcounter{mysubcase}\label{case:twoRegionsVerticalDirectionRefinement}),
it is done only in the $y$ direction. In both cases, we have a coarse domain on the
left-hand side.

\begin{figure}[h]
  \centering
  \includegraphics[width=\textwidth]{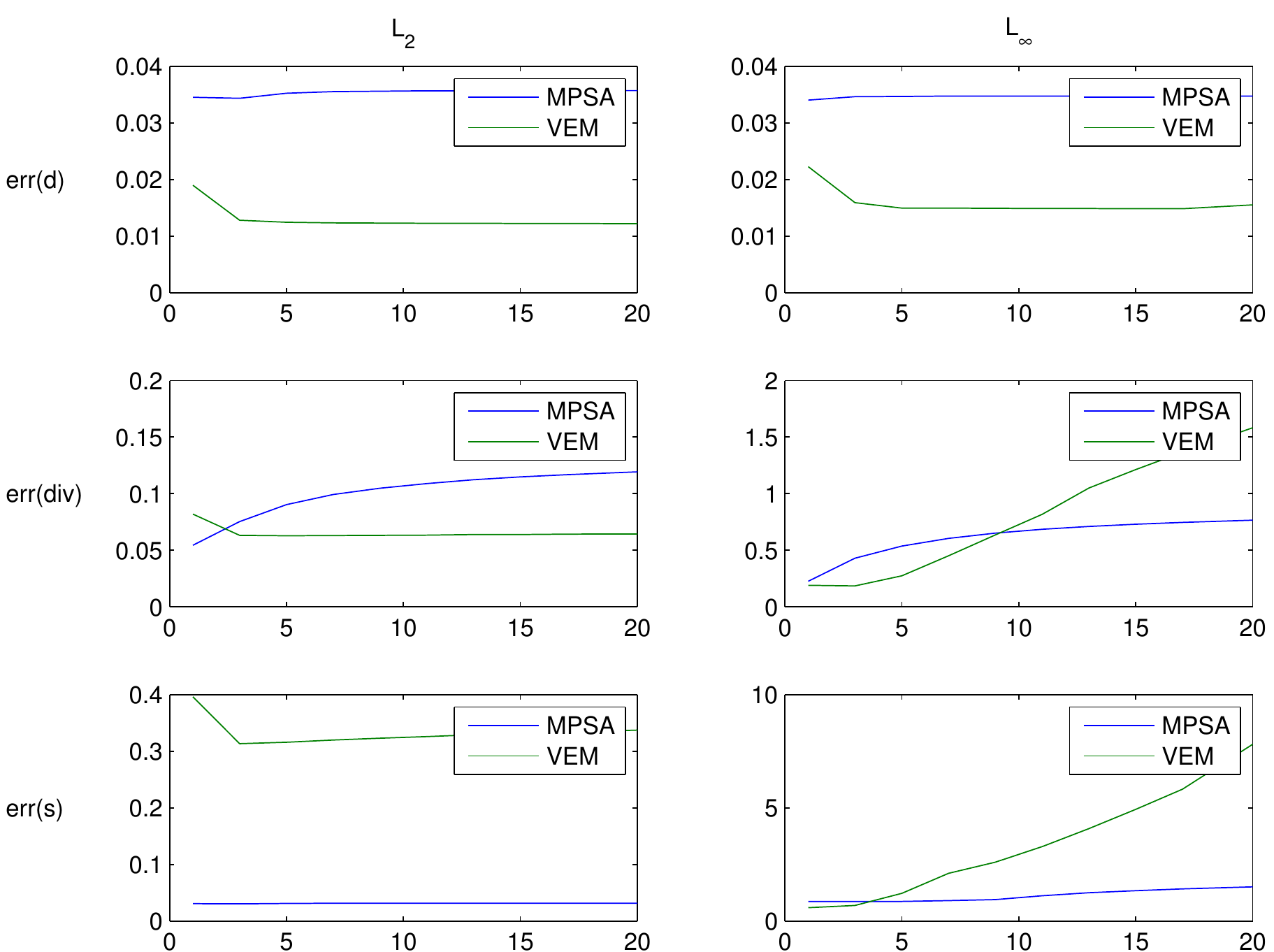}
  \caption{The region at the right-hand side is refined equally in both the $x$ and
    $y$ direction (Case \ref{case:twoRegionsIsotropicRefinement}). The $L^2$-norm (left) and
    $L^\infty$-norm (right) of the error are plotted for the displacement (upper
    row), the stress (middle row) and the divergence (lower row). The $x$-axis
    indicates the refinement factor in the right-hand side region.}
  \label{fig:convboxav10}
\end{figure}

\begin{figure}[h]
  \centering
  \includegraphics[width=0.9\textwidth]{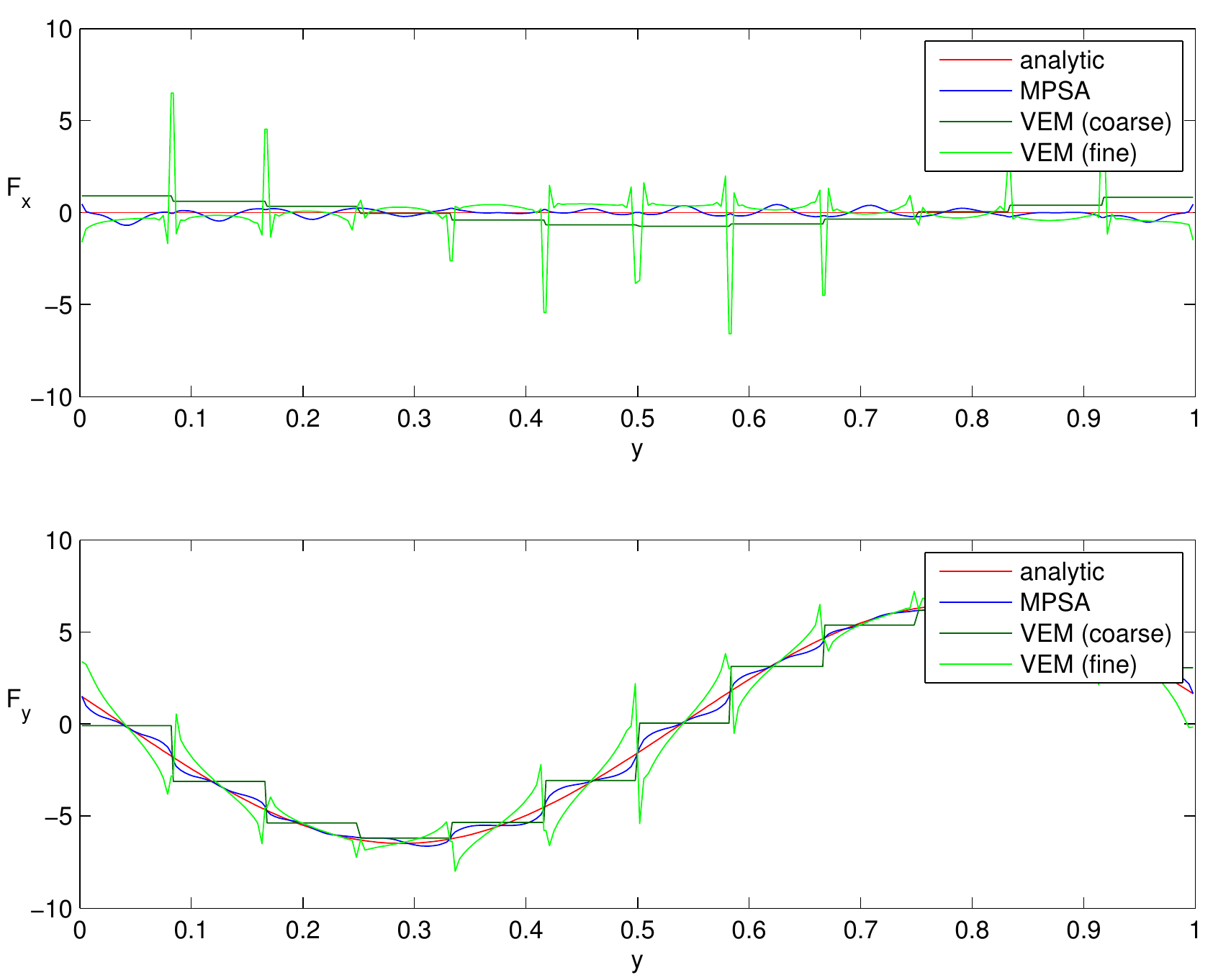}
  \caption{Plot of the forces at the interface between the two regions
    for Case \ref{case:twoRegionsIsotropicRefinement}. For the MPSA, the stress is computed on
    the faces so that the values of the stress at the interface are directly
    available. For the VE method, the stress is computed in the cells so that two
    values, one from the coarse and the other from the fine region, can be used to
    define the value of the stress at the interface.}
  \label{fig:stressfacesboxes10}
\end{figure}

\begin{figure}[h]
  \centering
  \includegraphics[width=0.7\textwidth]{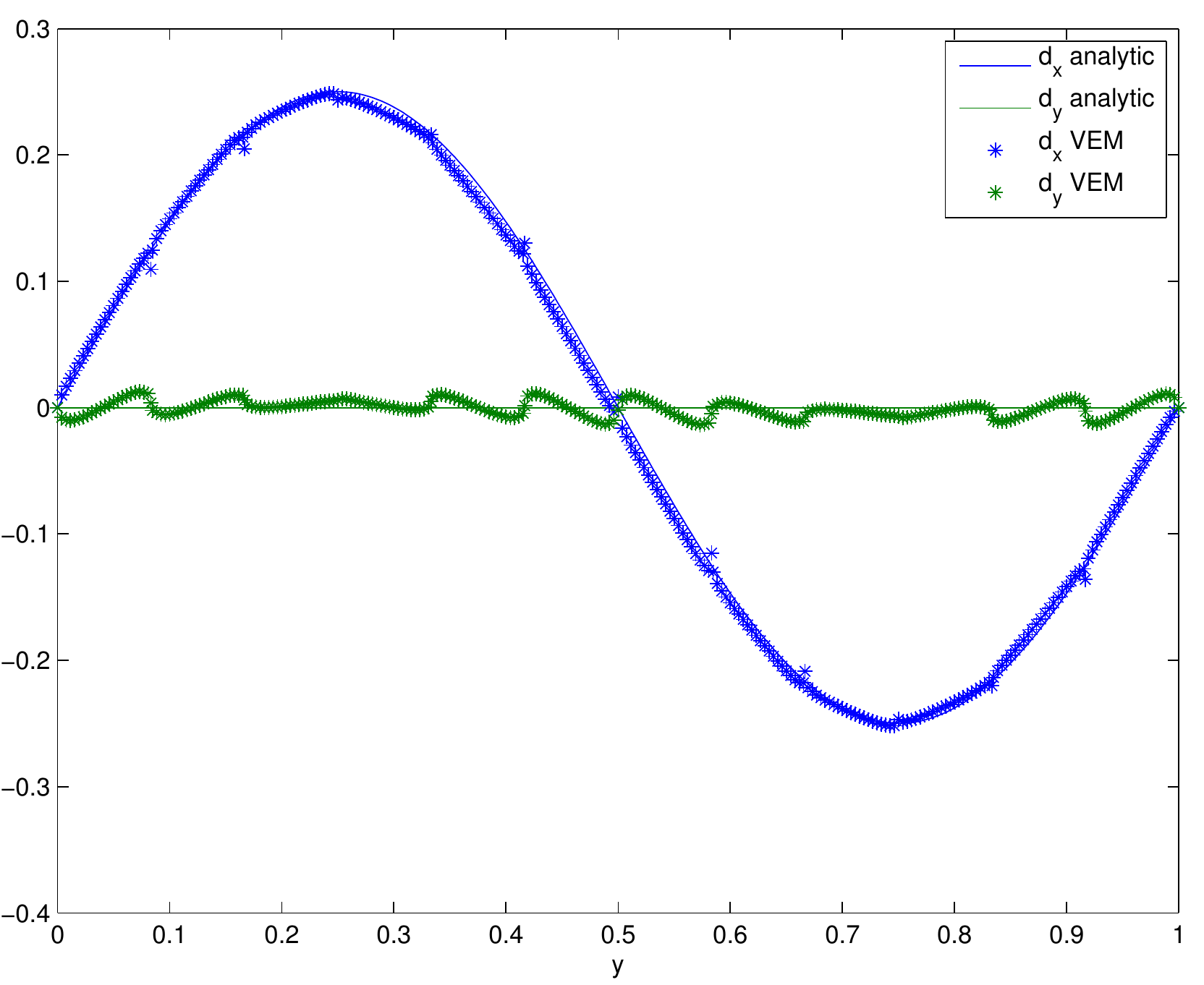}
  \caption{Plot of the error in displacement at the interface, only for the VE
    method (Case \ref{case:twoRegionsIsotropicRefinement})}
  \label{fig:vemdfacesboxesa10}
\end{figure}

\begin{figure}[h]
  \centering
  \includegraphics[width=\textwidth]{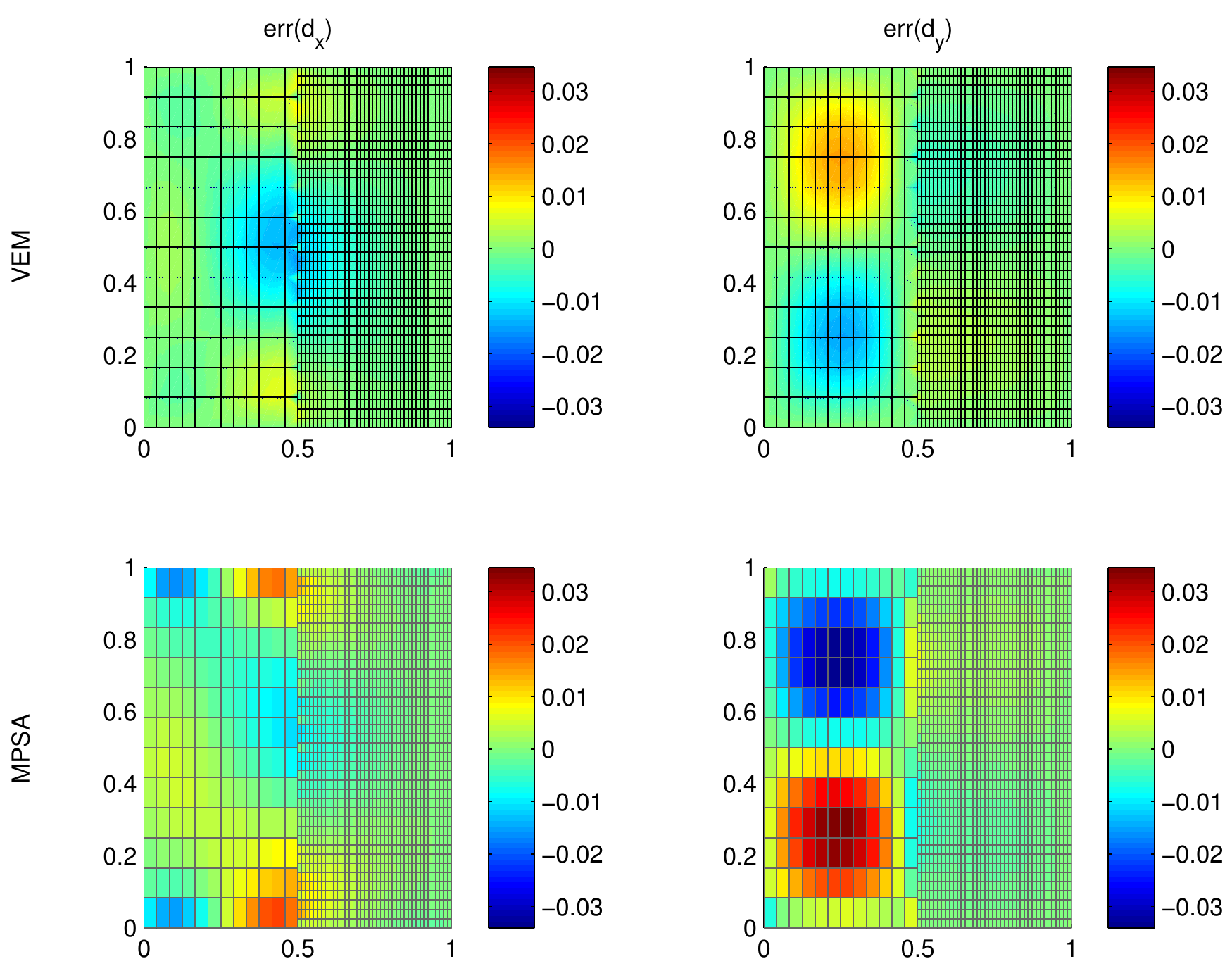}
  \caption{Plot of the displacement in both the vertical and horizontal directions
    (Case \ref{case:twoRegionsIsotropicRefinement})}
  \label{fig:errdboxessa10}
\end{figure}

In Figure \ref{fig:convboxav10}, we look at the error when the refinement on the
right-hand side is increased in both directions. By refinement factor, we mean the
number of sub-intervals that an edge of the initial grid is divided into to obtained
the refined grid. We observe that the error for the VE method increases significantly
for the divergence of the displacement and the stress in the $L^\infty$-norm. In the
$L^2$-norm, the increase is much less severe, which indicates that error is
essentially of local nature. In Figure \ref{fig:stressfacesboxes10}, we plot the
force at the interface between the two regions. For the VE method, the stress is
defined inside the cells so that we obtain two curves at the interfaces, one for the
coarse cells, the other for the fine cells. The force is computed at a cell interface
by integrating the product of the stress in the cell with the normal of the
interface. For the MPSA method, the force is defined on the edges and is therefore
readily available at the interface. We observe that the stress for the VE method is
strongly oscillating in the cells which belong to the refined region. For the
horizontal component of the force, the oscillations take the form of peaks, while the
force computed from the cells belonging to the coarse region is smooth and rather
close to the analytical value which is zero due to the symmetry of the problem. For
the $y$-component of the force, the analytical value is no longer constant. For the
VE method, the value computed from the cells of the fine region still presents
oscillation but, in addition, the value computed from the coarse region presents
strong variations, approximating the smooth analytical values by a staircase
function. Such behavior may be problematic if the solver is coupled with a fracture
model, typically non-linear, based on local value of the stress field. In comparison,
the MPSA method yields much smoother approximations. In Figure
\ref{fig:vemdfacesboxesa10}, we plot the error in displacement at the
interfaces. From this figure, it is clear that the local error concentrates at the
hanging nodes, see also Figure \ref{fig:errdboxessa10}.

\begin{figure}[h]
\includegraphics[width=\textwidth]{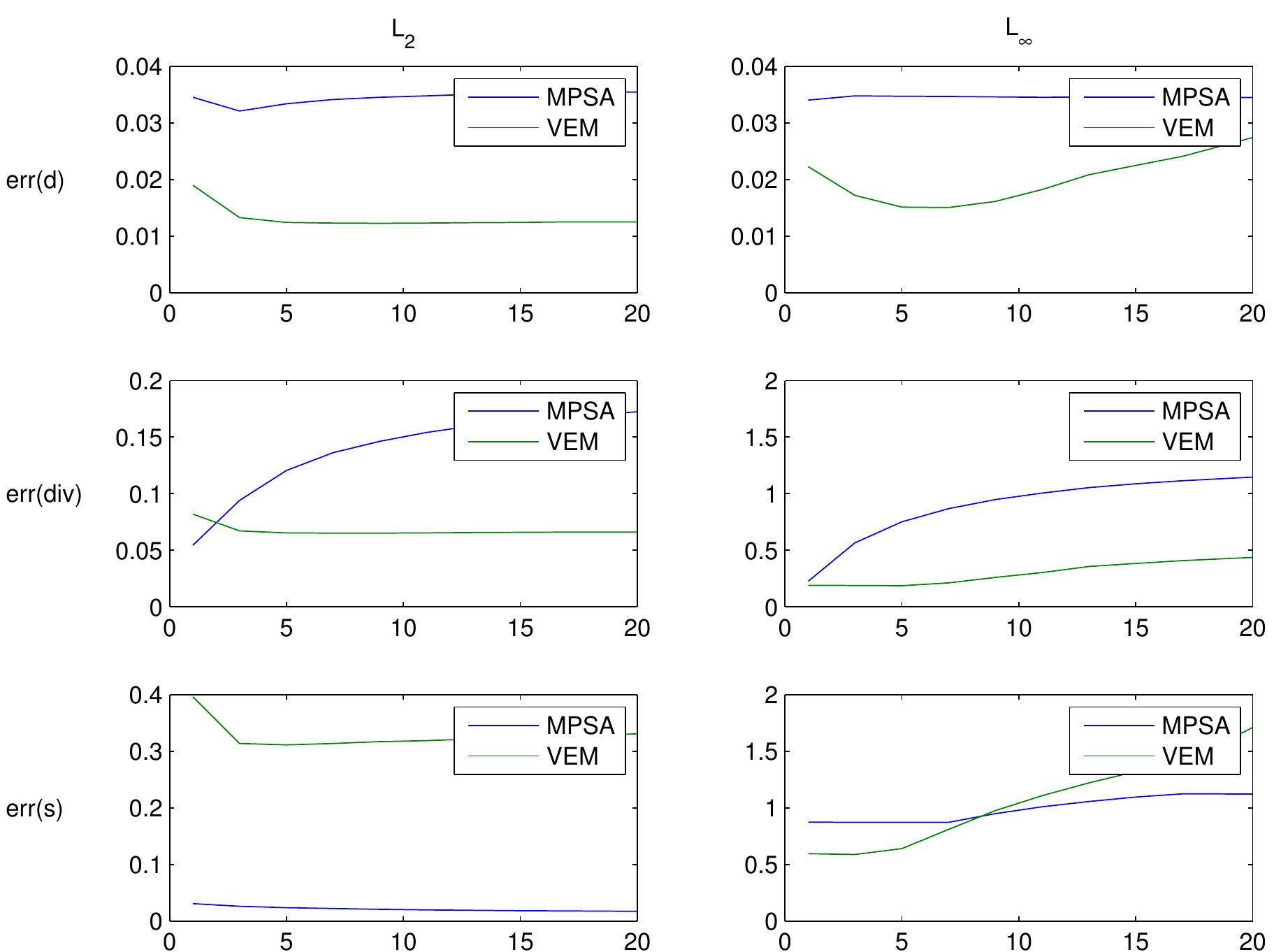}
\caption{The region at the right-hand side is refined only in the $y$ direction (Case
  \ref{case:twoRegionsVerticalDirectionRefinement}). The $L^2$-norm (left) and $L^\infty$-norm
  (right) of the error are plotted for the displacement (upper row), the stress
  (middle row) and the divergence (lower row). The $x$-axis indicates the refinement
  factor in the right-hand side region.}
\label{fig:convboxessv10}
\end{figure}

\begin{figure}[h]
  \begin{center}
    \includegraphics[width=0.6\textwidth]{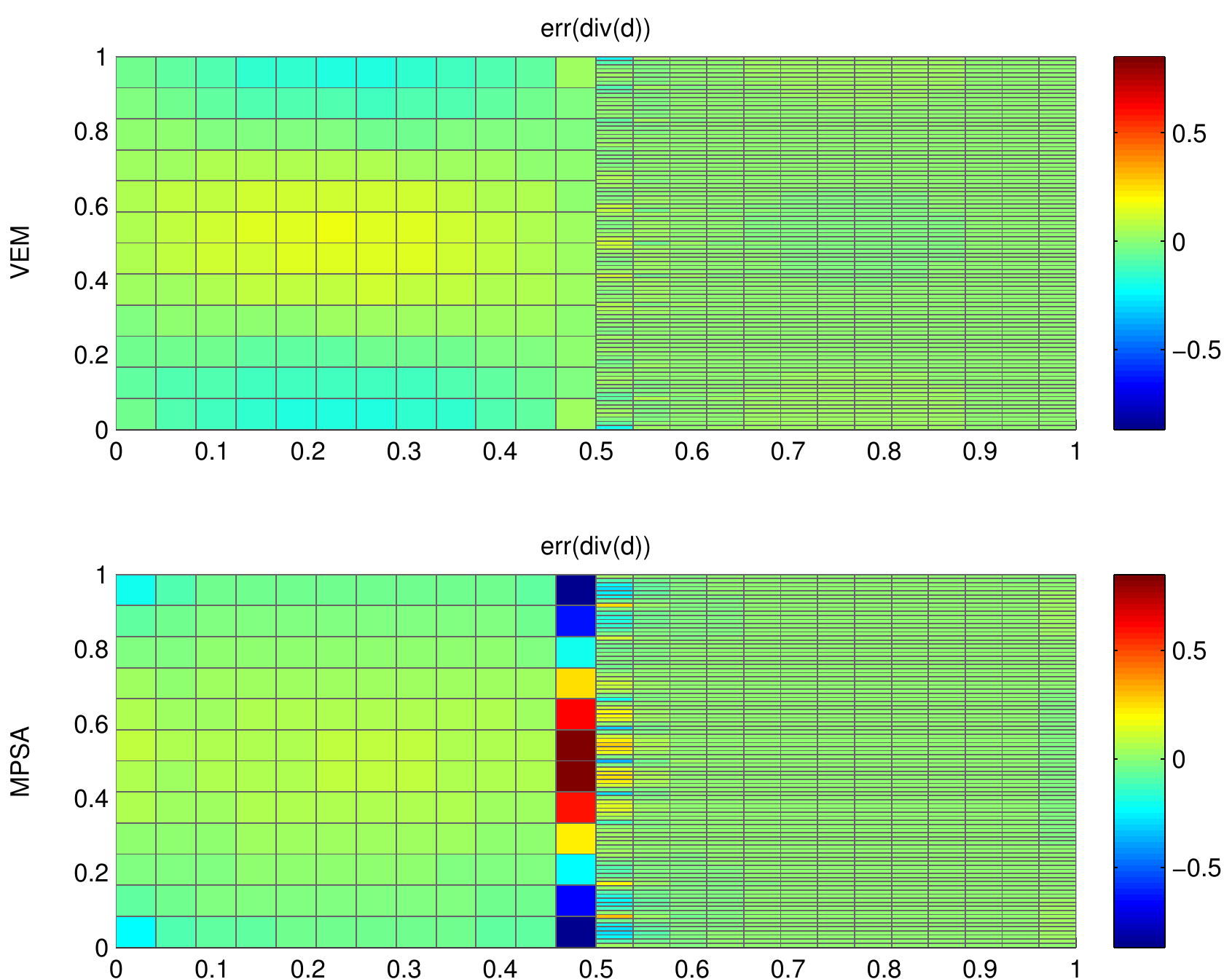}\\
    \includegraphics[width=\textwidth]{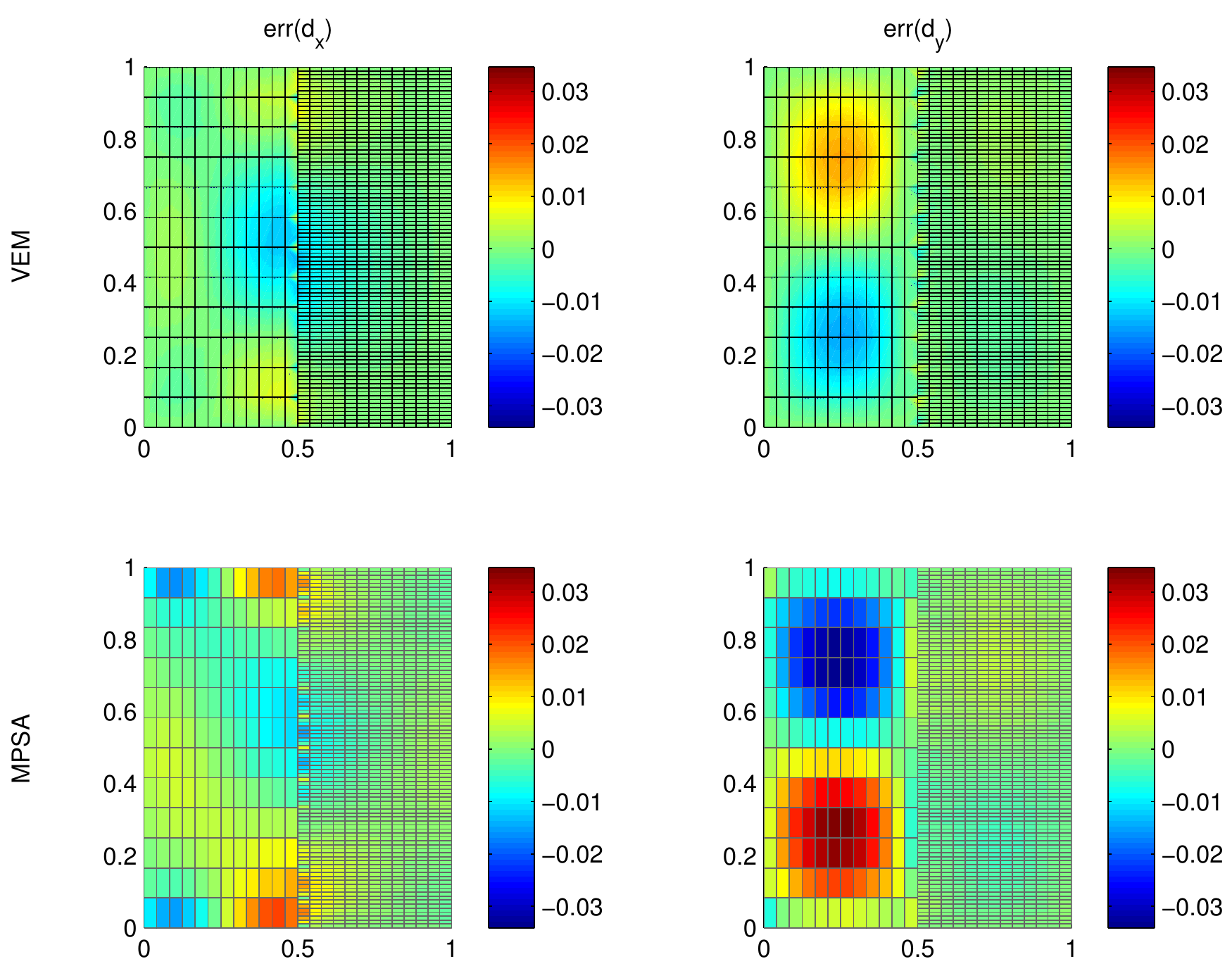}
  \end{center}
  \caption{Plot of the error in the divergence (two upper plots), the horizontal
    displacement (two lower left plots) and the vertical displacement (two lower
    right plots) for Case \ref{case:twoRegionsVerticalDirectionRefinement}. The
    vertical refinement ratio in the region on the right-hand side is equal to 20.}
\label{fig:errdivdboxesv10}
\end{figure}

\begin{figure}[h]
  \centering
  \includegraphics[width=0.9\textwidth]{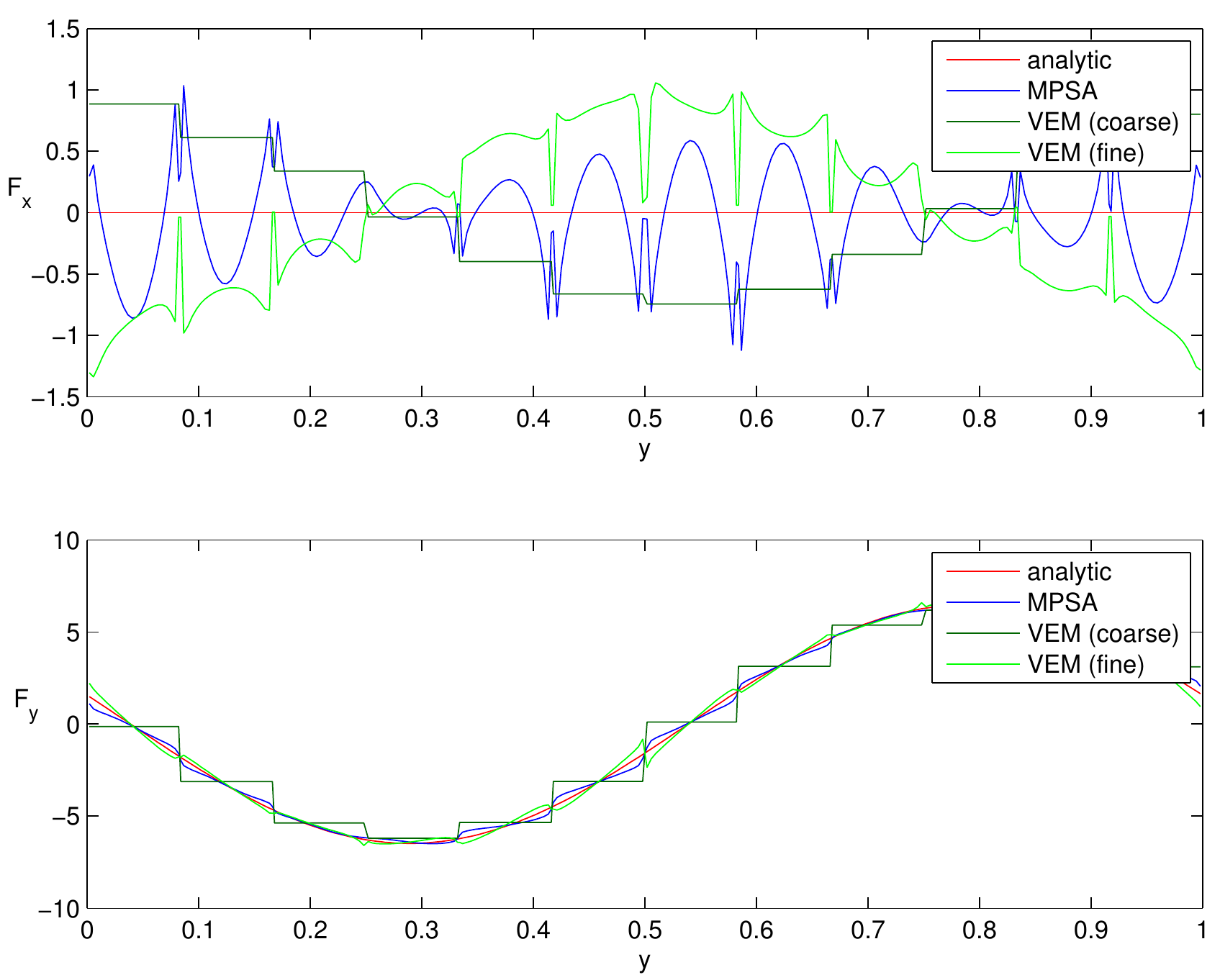}
  \caption{Plot of the forces at the interface, for Case
    \ref{case:twoRegionsVerticalDirectionRefinement}. The values are obtained in the
    same way as in the plot of Figure \ref{fig:stressfacesboxes10}.}
  \label{fig:stressfacesboxesv10}
\end{figure}

\begin{figure}[h]
  \centering
  \includegraphics[width=0.7\textwidth]{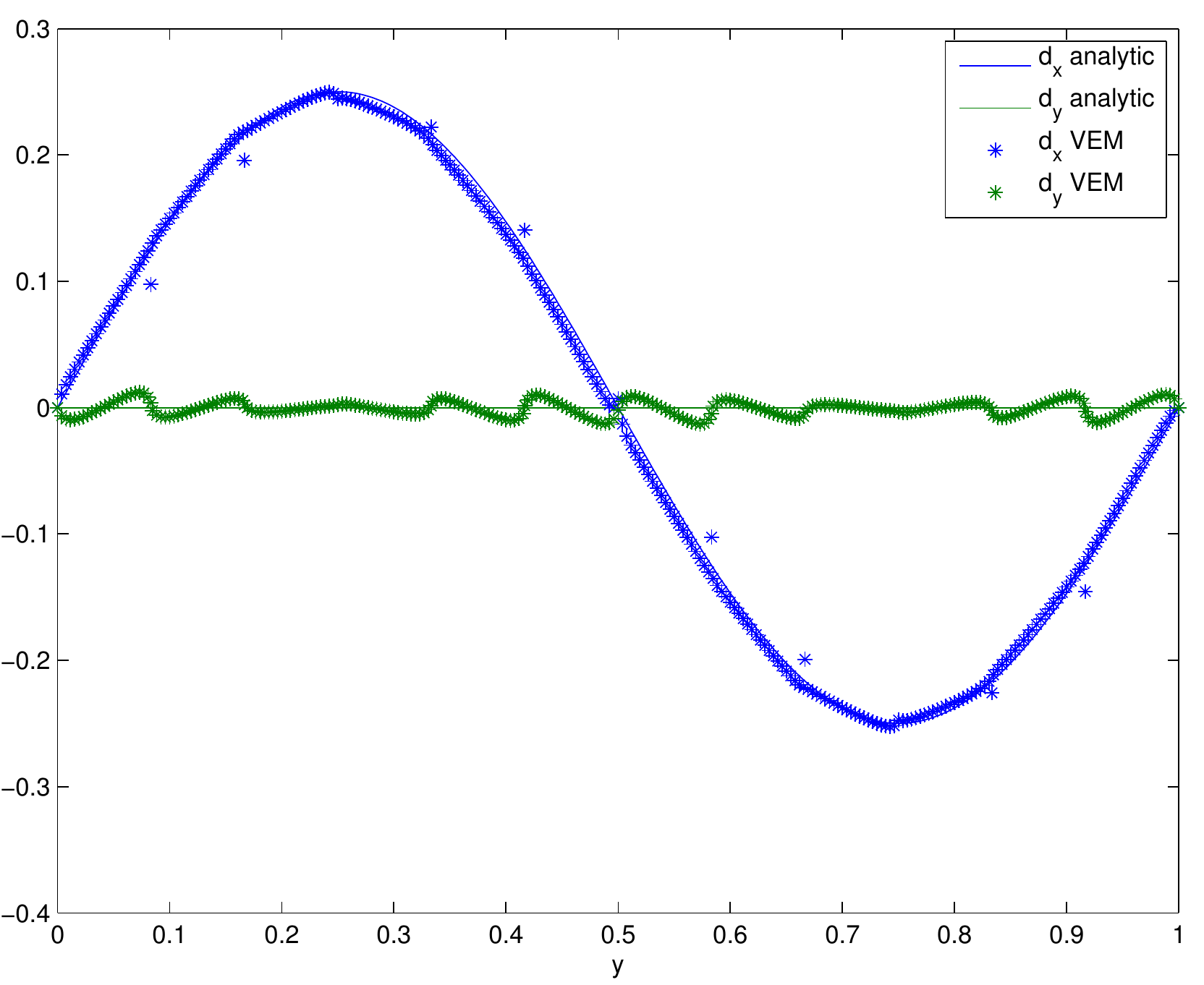}
  \caption{Plot of the error in displacement at the interface, only for the VE
    method (Case \ref{case:twoRegionsVerticalDirectionRefinement})}
  \label{fig:VEMdfacesboxesv10}
\end{figure}

Let us now consider the case where the refinement is done only in the $y$-direction
(Case \ref{case:twoRegionsVerticalDirectionRefinement}). The discretization at the
interface is the same as in the previous case but the cells at the right-hand side
get a relatively larger area and a larger aspect ratio. In Figure
\ref{fig:convboxessv10}, we observe that the error in the $L^\infty$-norm no longer
grows for the VE method. The strong oscillations in the $x$-component of the force
are smaller compared to the previous case, as we can observe by comparing Figures
\ref{fig:stressfacesboxes10} and \ref{fig:stressfacesboxesv10}. We can see in Figure
\ref{fig:VEMdfacesboxesv10} that most of the local error occurs at the non hanging
nodes. For the MPSA method, oscillations that were not present in Case
\ref{case:twoRegionsIsotropicRefinement} now appear in case
\ref{case:twoRegionsVerticalDirectionRefinement}. Moreover, we can see in Figure
\ref{fig:errdivdboxesv10} that the error spreads to the layer of coarse cell lying at
the interface, especially for the error in the divergence term. The calculation of
the divergence term in the MPSA method is based on the continuity points at the
boundary, which is calculated from solving a local problem.

\begin{figure}[h]
  \centering
  \includegraphics[width=0.7\textwidth]{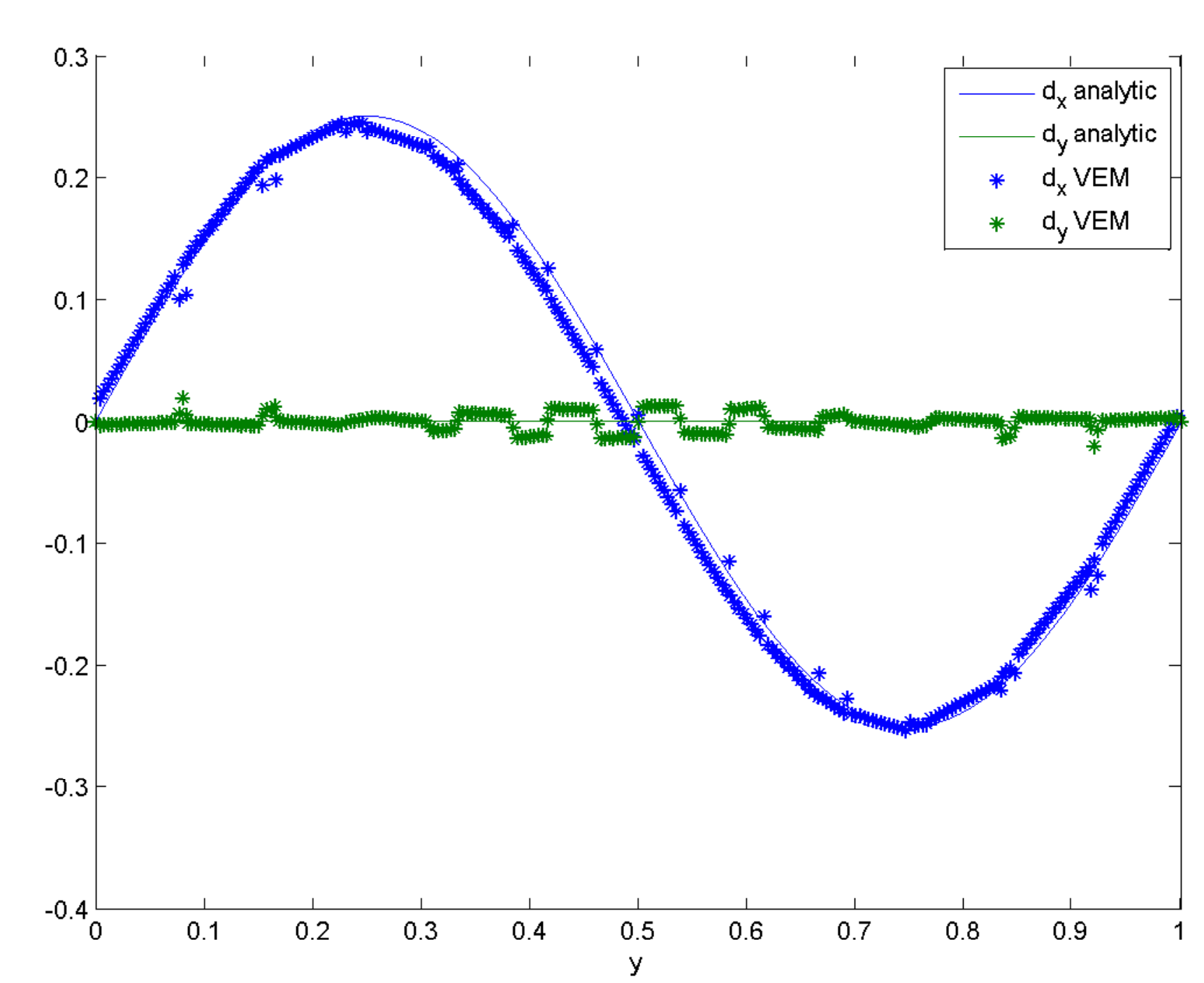}
  \caption{Plot of the error in displacement at the interface, only for the VE
    method in the case where there is no layer but 20 extra nodes on each face at the
    interface (Case \ref{case:refinedInterface}).}
  \label{fig:VEMdfacesboxesv17}
\end{figure}

\begin{figure}[h]
  \centering
  \includegraphics[width=0.3\textwidth]{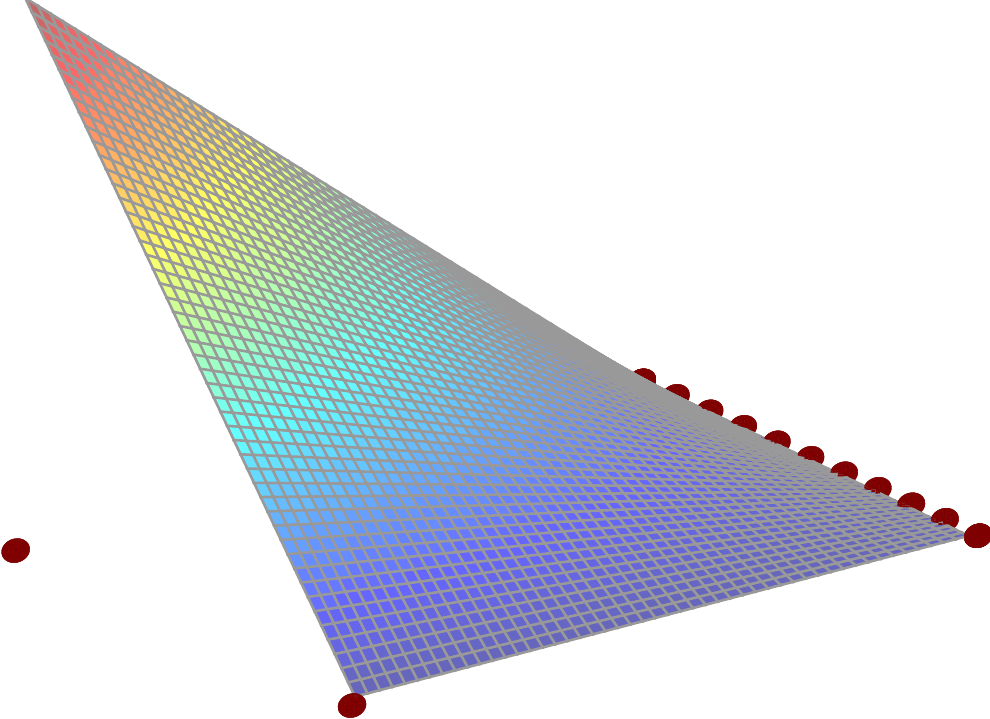}
  \includegraphics[width=0.3\textwidth]{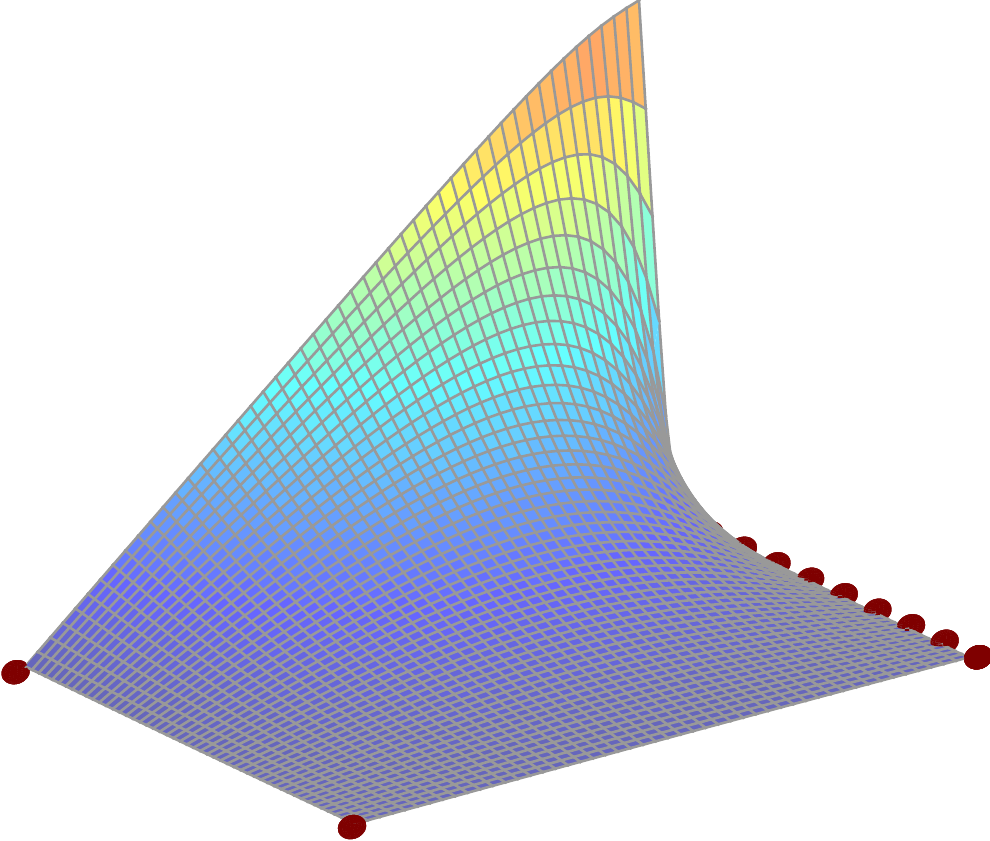}
  \includegraphics[width=0.3\textwidth]{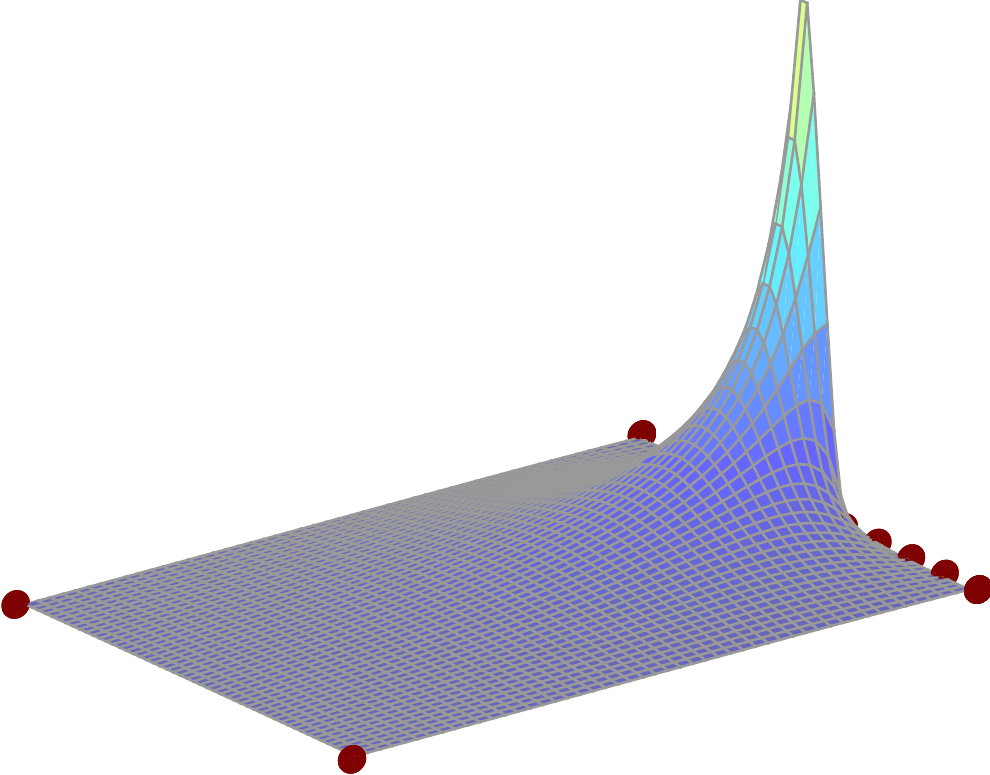}\\[1cm]
  \includegraphics[width=0.3\textwidth]{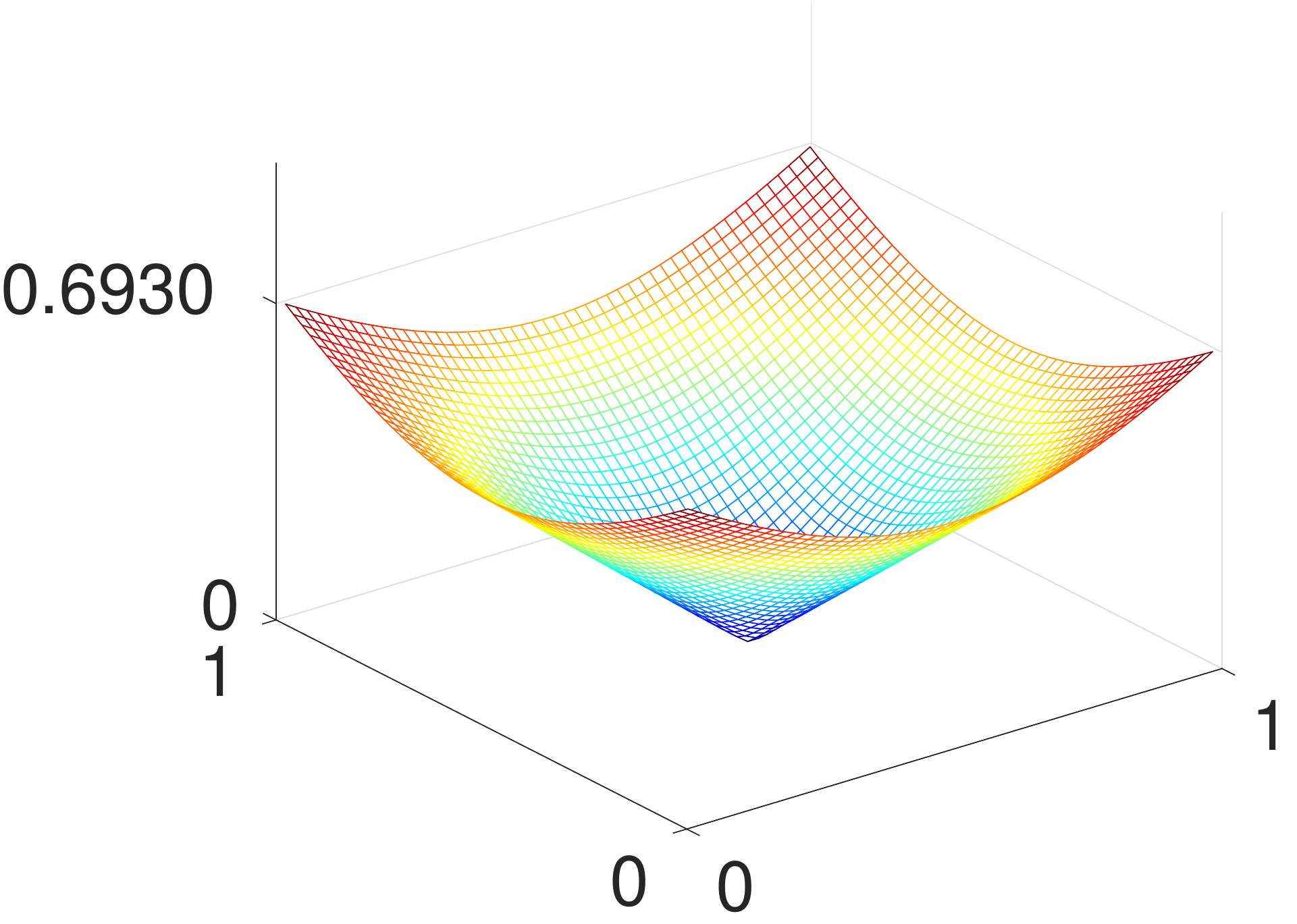}
  \includegraphics[width=0.3\textwidth]{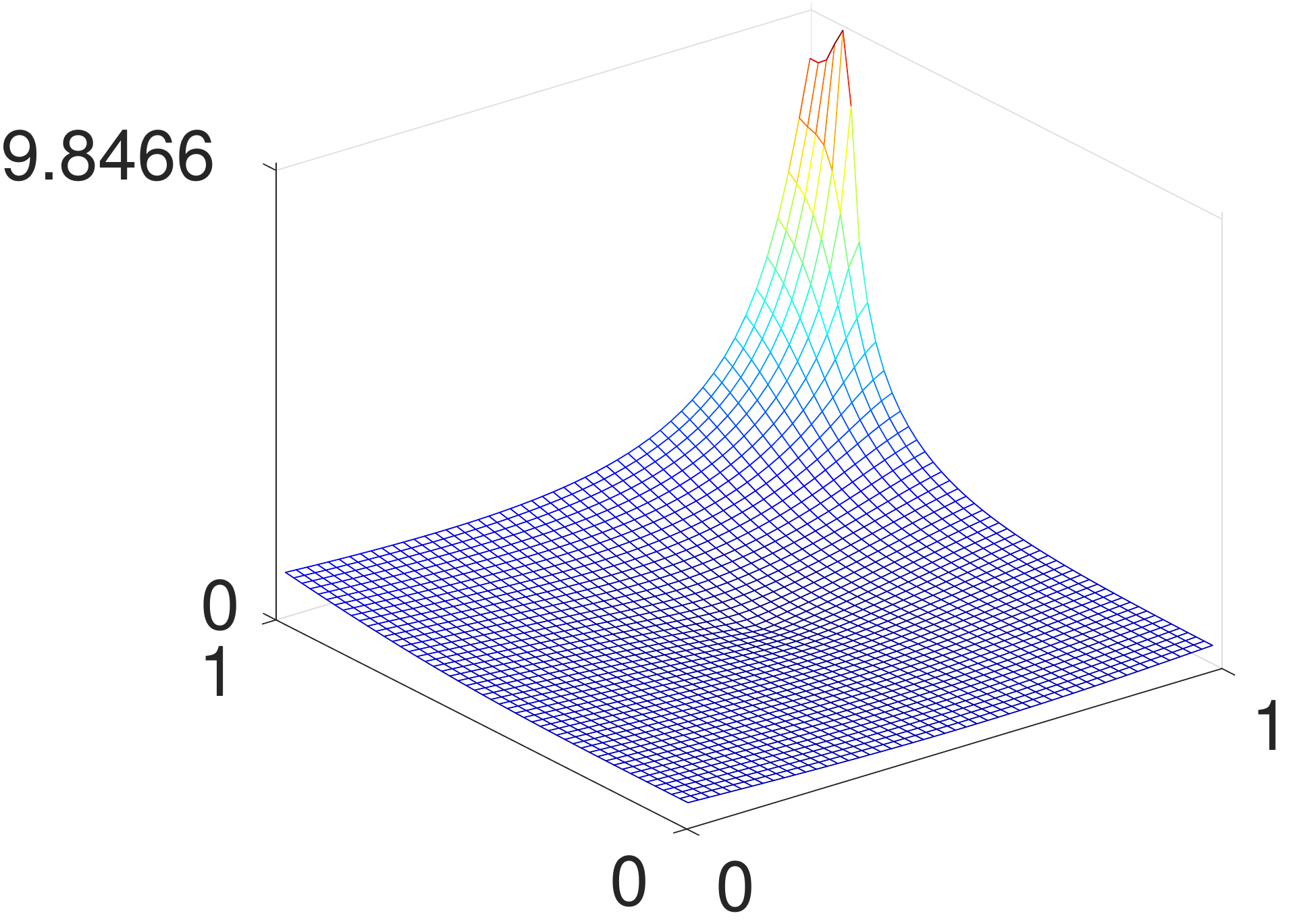}
  \includegraphics[width=0.3\textwidth]{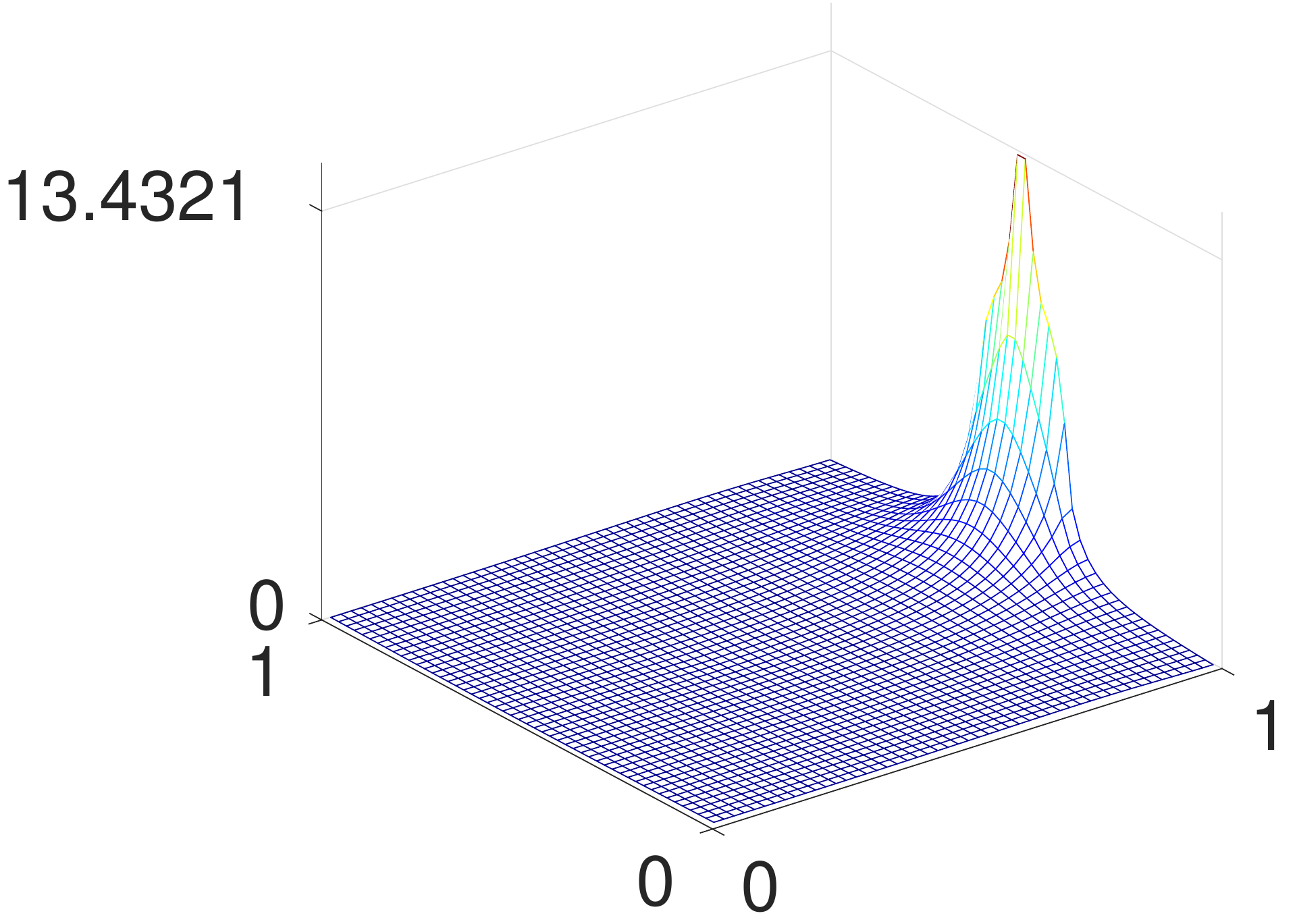}
  \caption{Illustration for Case \ref{case:refinedInterface} using the Laplace
    operator. First row: Plot a three virtual basis elements for a square cell with
    10 edges on one side. The red dots indicates the position of the nodes. We can
    sort the virtual basis in three types: Basis with two large edges (left), basis
    with a large and a small edge (middle), basis with two small edges
    (right). Second row: Plot of the norm of the residual of the projection, that is
    $\abs{\grad(\phi - \Pcal\phi)}^2$, for the basis function represented above.}
  \label{fig:virtualbasis}
\end{figure}

\begin{figure}[h]
  \centering
  \includegraphics[width=0.3\textwidth]{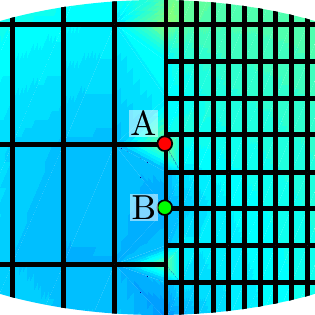}
  \caption{Zoom of the interface region showing the error in horizontal displacement
    for the VE method in case \ref{case:twoRegionsIsotropicRefinement}, see Figure
    \ref{fig:errdboxessa10}. The error is different at nodes with a large edge (type
    A) and nodes which connect to two small edges (type B).}
  \label{fig:zoom}
\end{figure}

Finally, we setup a case where the two regions have the same coarse mesh but we add
extra nodes at the interface (Case
\ref{case:refinedInterface}\refstepcounter{mysubcase}\label{case:refinedInterface}). In
this way, we remove the difference in volumetric refinement between the two regions
and isolate the effect of edge refinements. In Figure \ref{fig:VEMdfacesboxesv17}, we
plot the error displacement for the VE method at the interface and observe that the
nodes which belong to both a long edge and a short edge behave differently than the
nodes that belong to two small edges (in this case, the hanging nodes). This
observation complies with the results observed earlier and presented in Figure
\ref{fig:vemdfacesboxesa10} and Figure \ref{fig:VEMdfacesboxesv10}, which shows that
the VE method spreads the error unevenly between these two types of nodes. Second, it
shows that it is related to edge refinement. In the VE method, the basis elements are
not computed, only the degrees of freedom are used for the assembly and linear
approximations remain exact but, in the case of elements with many nodes, the basis
elements will be highly non-linear. We illustrate this in Figure
\ref{fig:virtualbasis} where we compute some of the virtual basis elements
constructed using \textit{harmonic lifting} as in \cite{gain2014}, see also
\cite{beirao2013basic}. We consider the same type of cells as the ones which lie at
the interface in Case \ref{case:refinedInterface}, reducing the refinement to ten
nodes in order to make the pictures easier to read. For simplicity, this illustration
has been created using the Laplace operator and not the linear elasticity
operator. We can sort the virtual basis in three categories:: Basis with two large
edges (type I), basis with a large and a small edge (type II), basis with two small
edges (type II). The virtual basis elements have very sharp gradients in small
regions and are almost flat elsewhere so that most of their energy is concentrated in
high frequencies. In this case, the projection operator $\Pcal$ over linear function,
see section \ref{subsec:vem}, does not provide a good approximation and most of the
contribution for this basis element will be handled by the regularization term,
$s_K$, which is only a poor substitute for $a_K$. We have computed the residual part
for the three basis,
\begin{equation*}
  \frac{a_K(\phi - \Pcal\phi, \phi - \Pcal\phi)}{a_K(\phi, \phi)} =
  \begin{cases}
    0.49&\text{ if $\phi$ is of type I,}\\
    0.90&\text{ if $\phi$ is of type II,}\\
    0.99&\text{ if $\phi$ is of type III.}
  \end{cases}
\end{equation*}
Since in this case the length of the large and small edges are $L=1$ and $l=0.1$,
respectively, these computations confirm the following orders of magnitude,
\begin{equation*}
  \frac{a_K(\Pcal\phi, \Pcal\phi)}{a_K(\phi, \phi)} \approx
  \begin{cases}
    1&\text{ if $\phi$ is of type I,}\\
    l/L&\text{ if $\phi$ is of type II,}\\
    (l/L)^2&\text{ if $\phi$ is of type III,}
  \end{cases}
\end{equation*}
which can be obtained by roughly estimating the area of the support of the gradient
of the basis function. In Figure \ref{fig:vemdfacesboxesa10}, we observe that, at the
interface region, the displacement values obtained at nodes that are connected to a
large edge (which we denote type A) have different errors that the other nodes (which
we denote type B), see a zoom on this region in Figure \ref{fig:zoom}.


\clearpage

\subsection{Case \ref{case:layer}: Layer between two domains}
\refstepcounter{mycase}\label{case:layer}

\begin{figure}[h]
  \centering
  \includegraphics[width=0.9\textwidth]{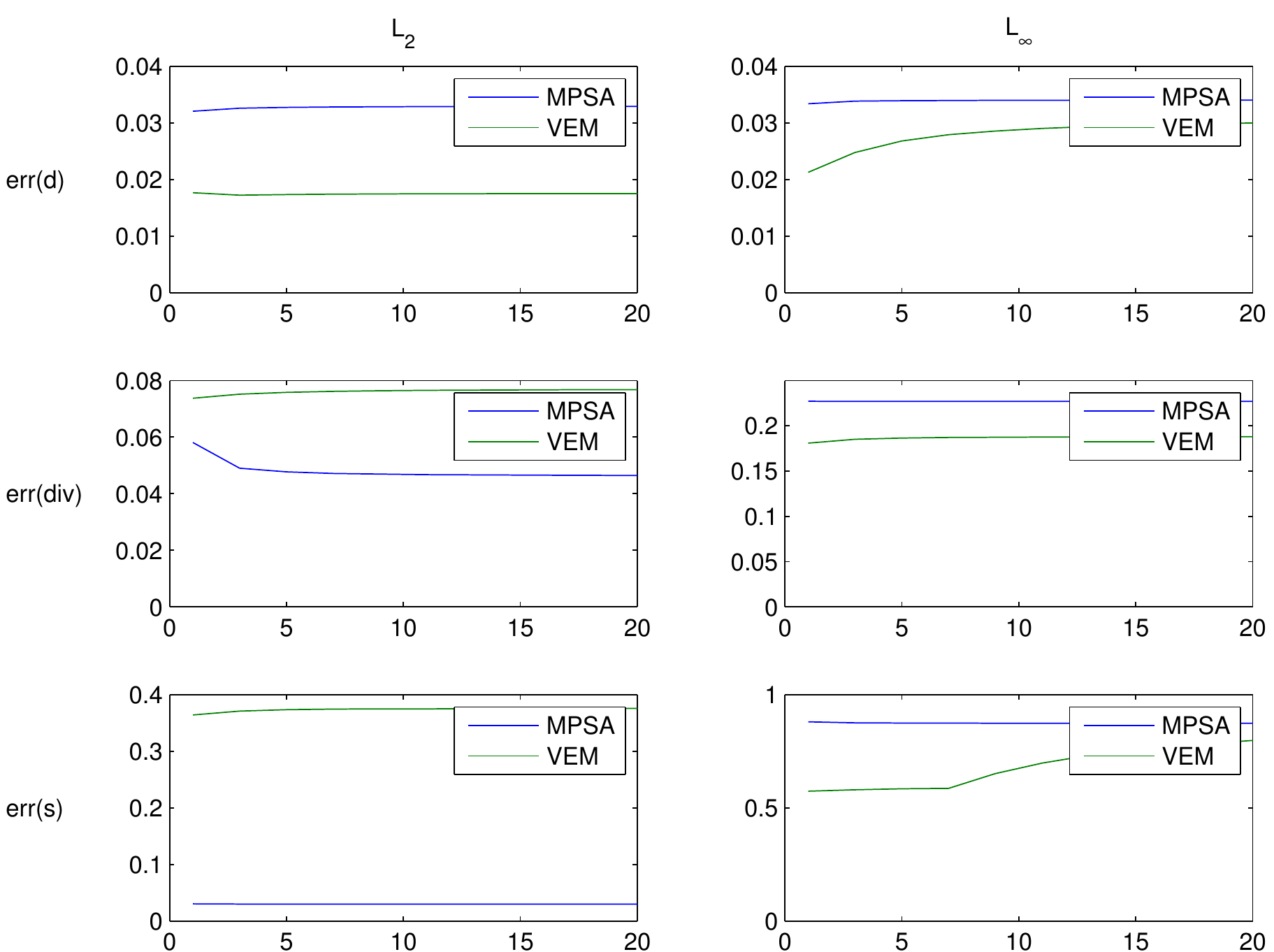}
  \caption{The width of the layer is decreased and no vertical refinement in the
    layer is used (Case \ref{case:layerNoRefinement}). The $L^2$-norm (left) and
    $L^\infty$-norm (right) of the error are plotted for the displacement (upper
    row), the stress (middle row) and the divergence (lower row). The $x$-axis
    indicates the reduction factor for the width of the layer.}
  \label{fig:convboxes2eps0}
\end{figure}

\begin{figure}[h]
  \centering
  \includegraphics[width=0.6\textwidth]{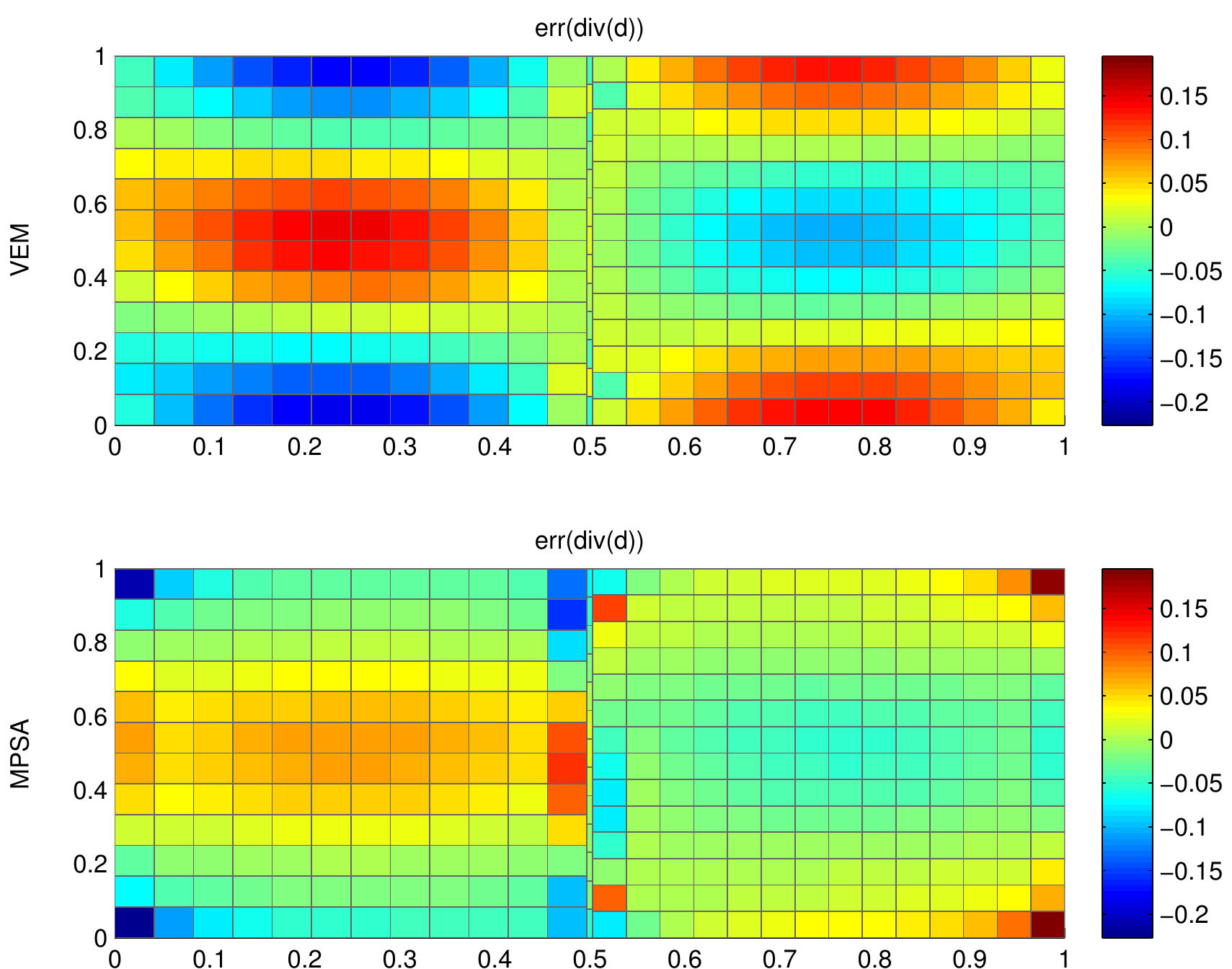}\\
  \includegraphics[width=0.8\textwidth]{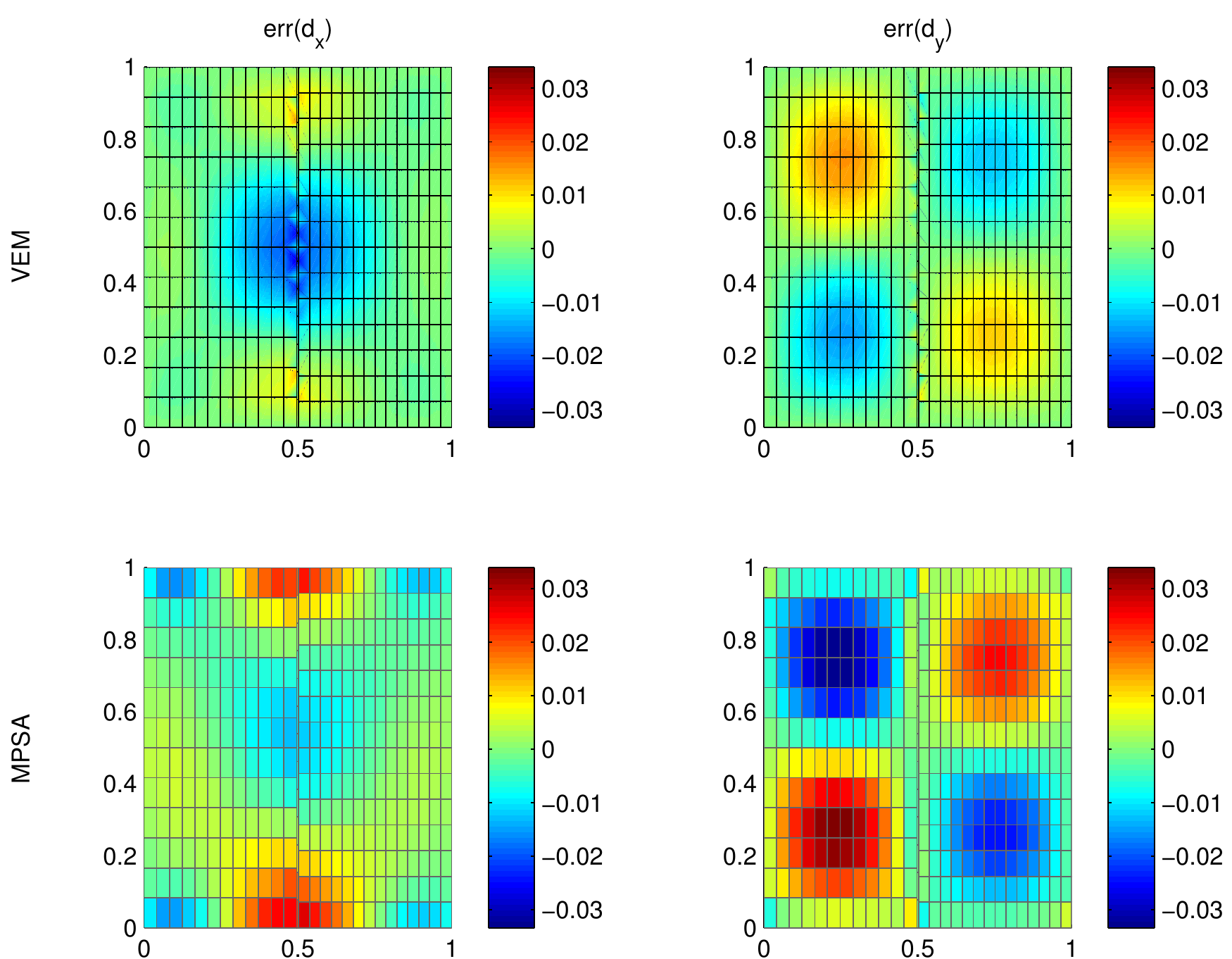}
  \caption{Error in the divergence (two upper figures) and the displacement in both
    $x$ and $y$ directions (Case \ref{case:layerNoRefinement}). The width of the layer is 20
    times smaller than the adjacent cells.}
  \label{fig:errdboxes2eps0}
\end{figure}

\begin{figure}[h]
  \centering
  \includegraphics[width=0.9\textwidth]{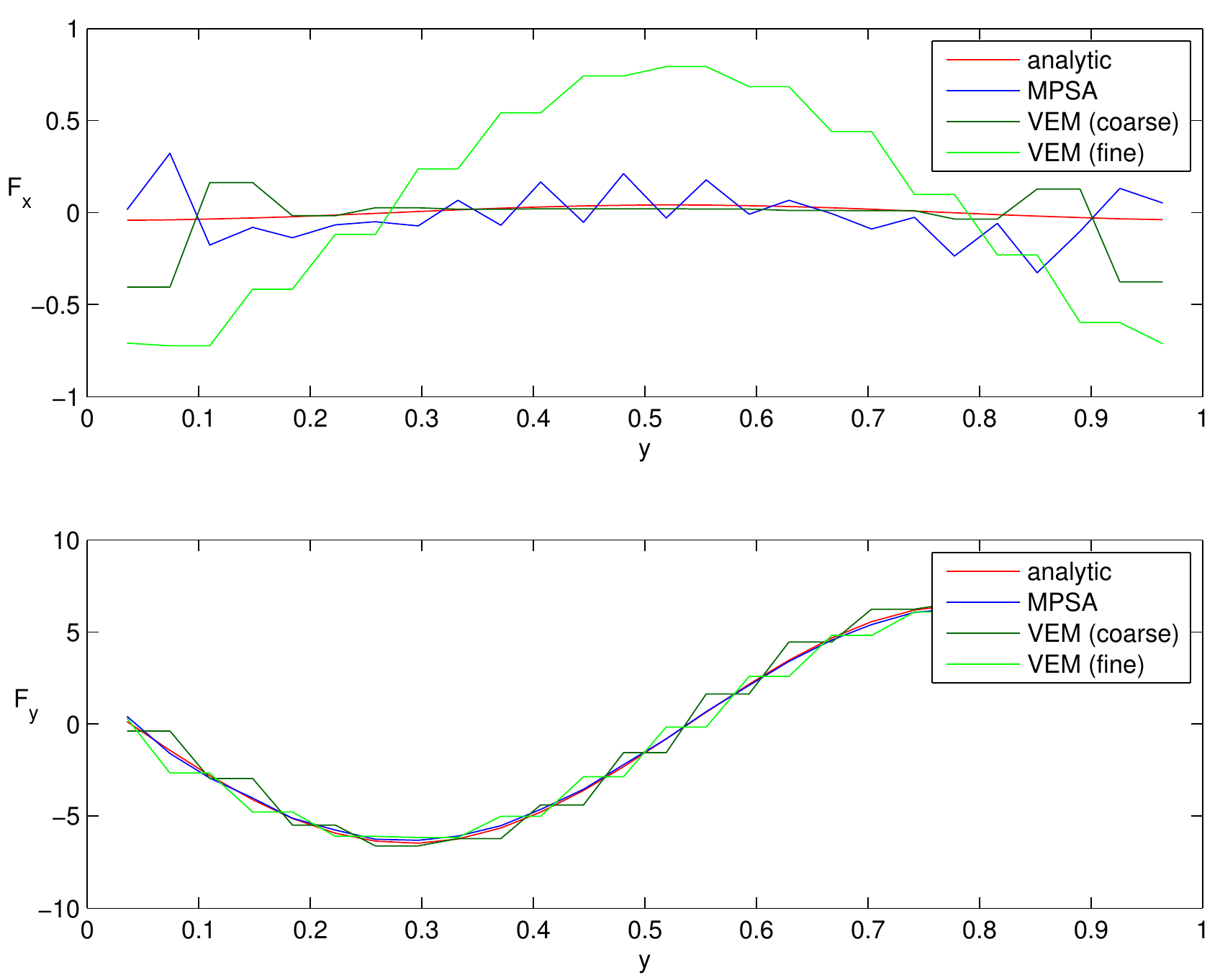}
  \caption{Plot of the forces at the interface, for Case
    \ref{case:layerNoRefinement}. The values are obtained in the same way as in the plot of
    Figure \ref{fig:stressfacesboxes10}.}
  \label{fig:stressfacesboxes2eps0}
\end{figure}

In subsurface flow, an important part of the flow is concentrated in the fractures of
the rock. We set up test cases that reproduce the geometry of a fracture by
introducing a thin layer in an otherwise uniform Cartesian grid. If this layer were a
fracture or a damaged zone, then it would have very different mechanical properties
than the rest of the matrix, but, in our test case, the layer is assigned the same
properties as the rest of the domain so that the analytical solution given by
\eqref{eq:exact_disp_field}. In this way, we isolate the errors of the two numerical
methods which are induced by a typical geometrical discrete representation of a
fracture, without including the mechanical effects of the fracture itself. First, we
consider a test case (Case
\ref{case:layerNoRefinement}\refstepcounter{mysubcase}\label{case:layerNoRefinement}),
where the layer is discretized with the same level of refinement in the $y$ direction
as the rest of the matrix. In Figure \ref{fig:convboxes2eps0}, we let the width of
the layer get thinner and thinner. We observe that the error does not grow,
indicating the robustness of both methods with respect to the thickness of the
layer. Figure \ref{fig:errdboxes2eps0} gives the error in displacement and divergence
of displacement. We observe that the error in displacement is more localized for the
VE method and more spread for the MPSA method, which is consistent with the results
of Figure \ref{fig:convboxes2eps0} when comparing the $L^\infty$ norm and
$L^2$-norms. In Figure \ref{fig:stressfacesboxes2eps0}, we present a plot of the
forces at the interface between the layer and the region with coarse cells. We choose
the left interface of the layer but the results on the other interface have the same
characteristics. We observe that the MPSA method gives slightly but not significantly
better approximation of the force.

\begin{figure}[h]
  \centering
  \begin{tabular}[h]{cc}
    \includegraphics[width=0.45\textwidth]{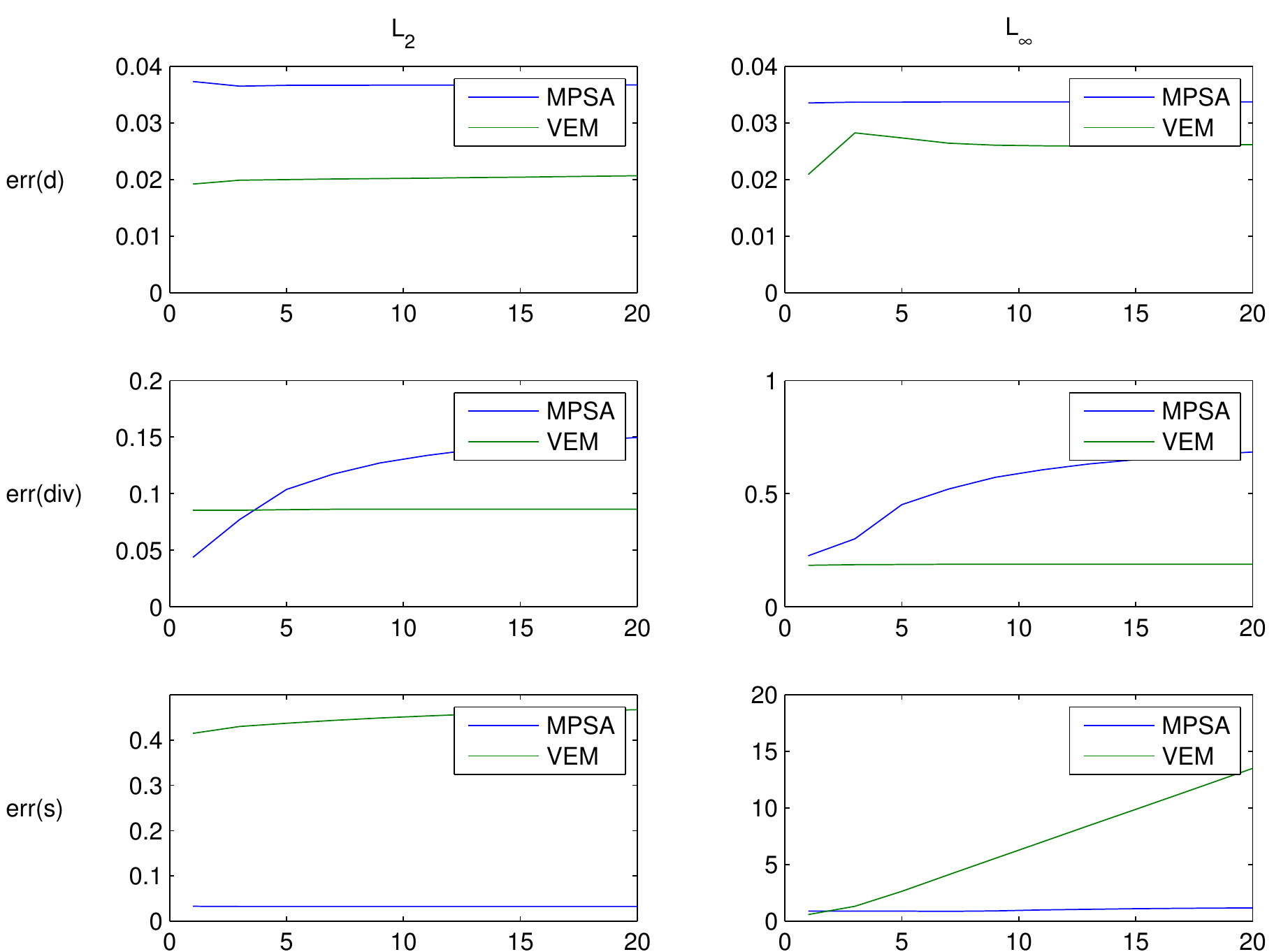}&\includegraphics[width=0.45\textwidth]{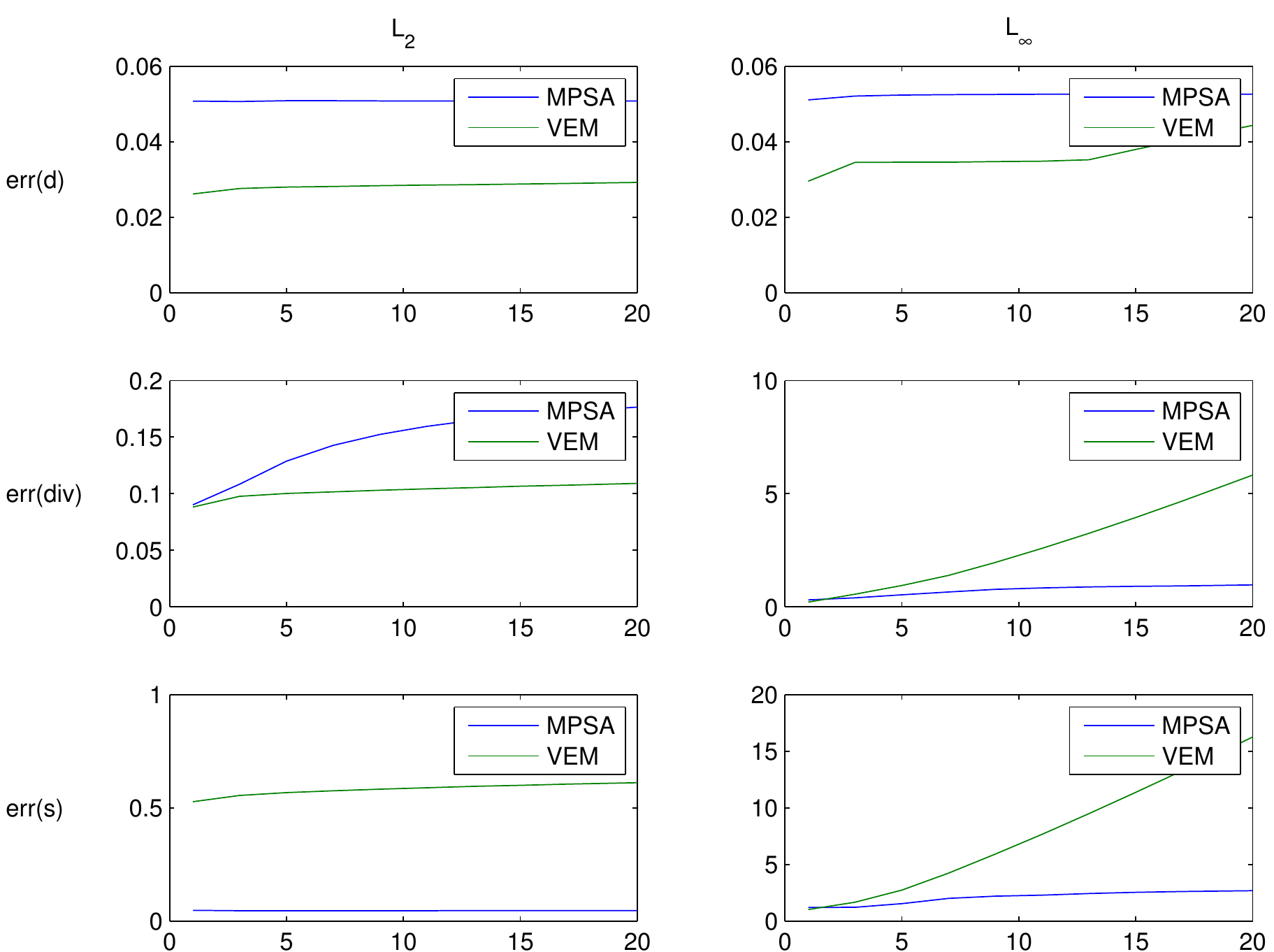}\\
    Case \ref{case:layerRefinementCartesian} & Case \ref{case:layerRefinementTwisted}
  \end{tabular}
  \caption{The width of the layer is decreased but the aspect ratio of the cells in
    the layer is preserved.  The $L^2$-norm (left columns) and $L^\infty$-norm (right
    columns) of the error are plotted for the displacement (upper row), the stress
    (middle row) and the divergence (lower row). On the left, we have the Cartesian
    case (Case \ref{case:layerRefinementCartesian}) and the twisted case on the right (Case
    \ref{case:layerRefinementTwisted}). In all plots, the $x$-axis indicates the value of the
    reduction factor of the layer width.}
  \label{fig:convboxesv2eps7}
\end{figure}

\begin{figure}[h]
  \centering
  \includegraphics[width=0.9\textwidth]{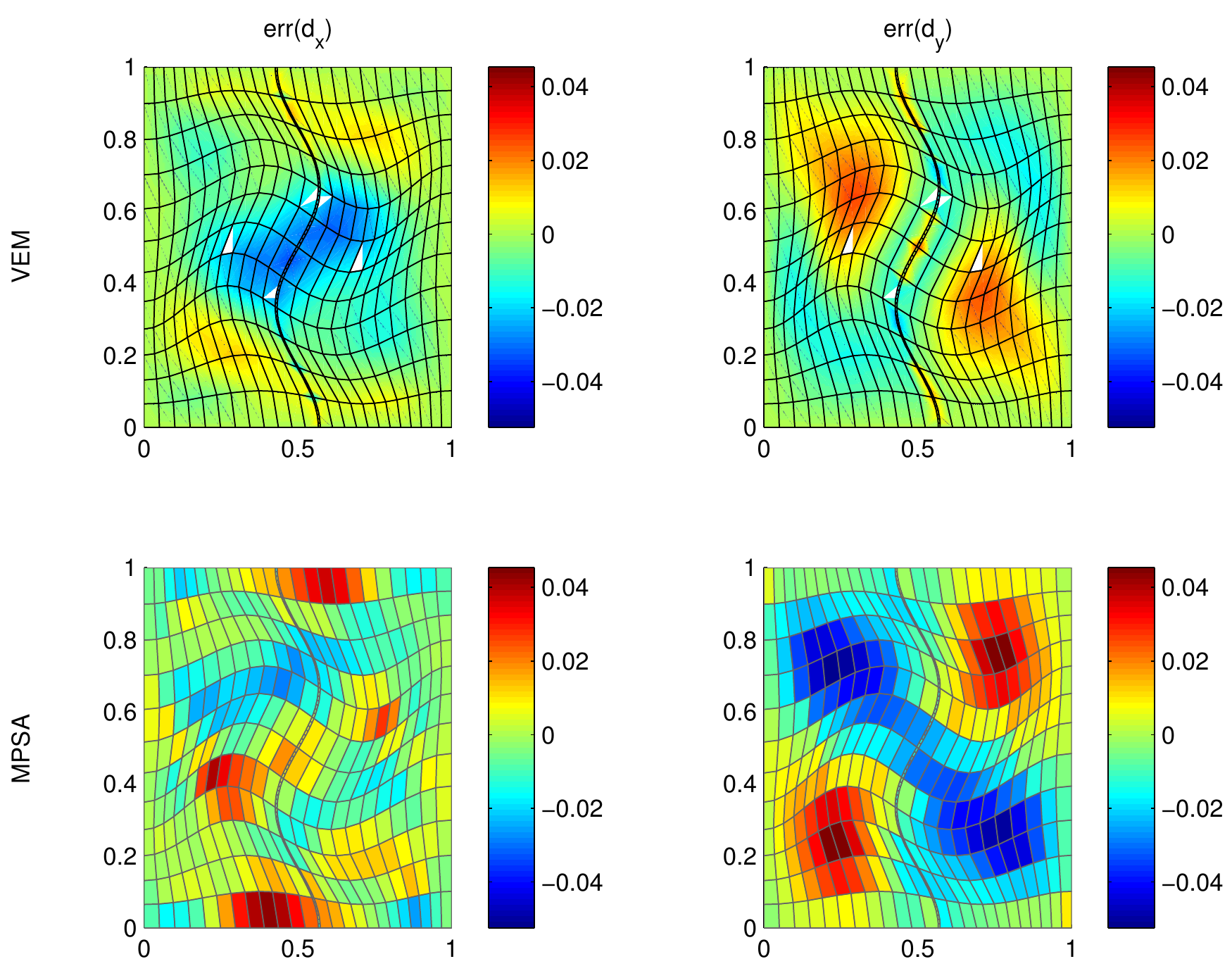}
  \caption{Plot of the error of the displacement in $x$ and $y$ directions (Case
    \ref{case:layerRefinementTwisted}).}
  \label{fig:errdboxesv2eps7}
\end{figure}

\begin{figure}[h]
  \centering
  \includegraphics[width=0.9\textwidth]{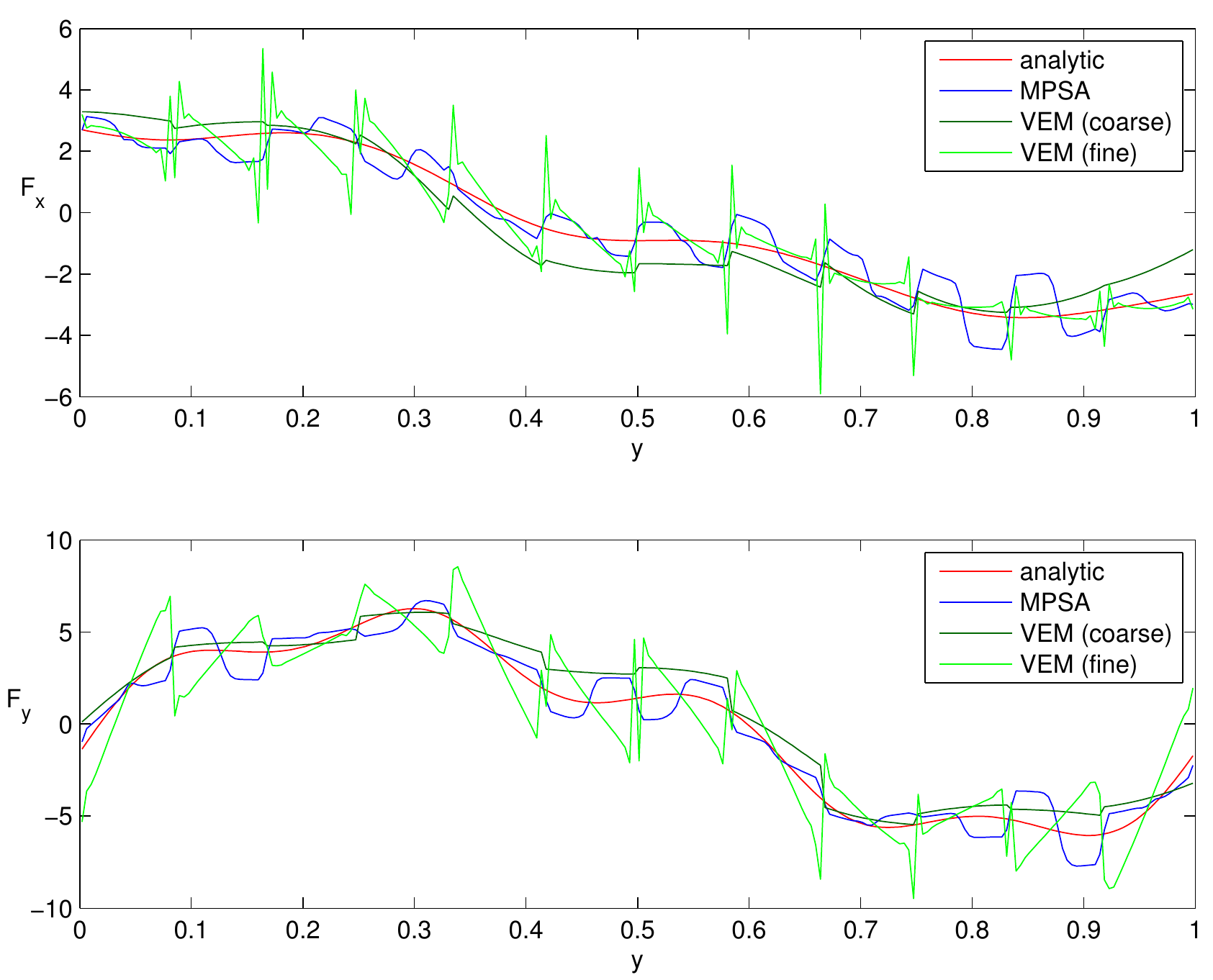}
  \caption{Plot of the forces at the interface, for Case \ref{case:layerRefinementTwisted}. The
    values are obtained in the same way as in the plot of Figure
    \ref{fig:stressfacesboxes10}.}
  \label{fig:stressfacesboxesv2eps7}
\end{figure}

\begin{figure}[h]
  \centering
  \includegraphics[width=0.45\textwidth]{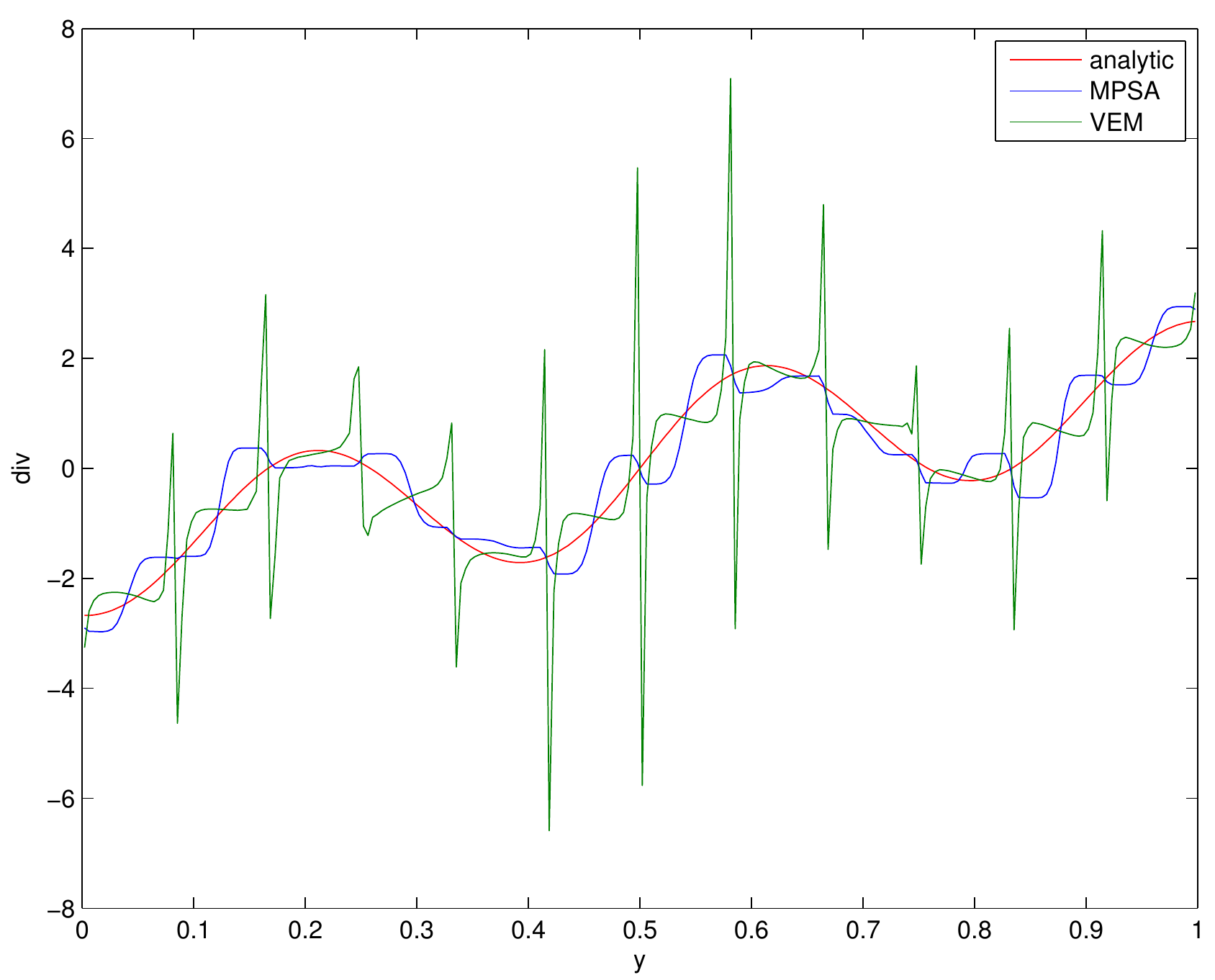}\\
  \includegraphics[width=0.6\textwidth]{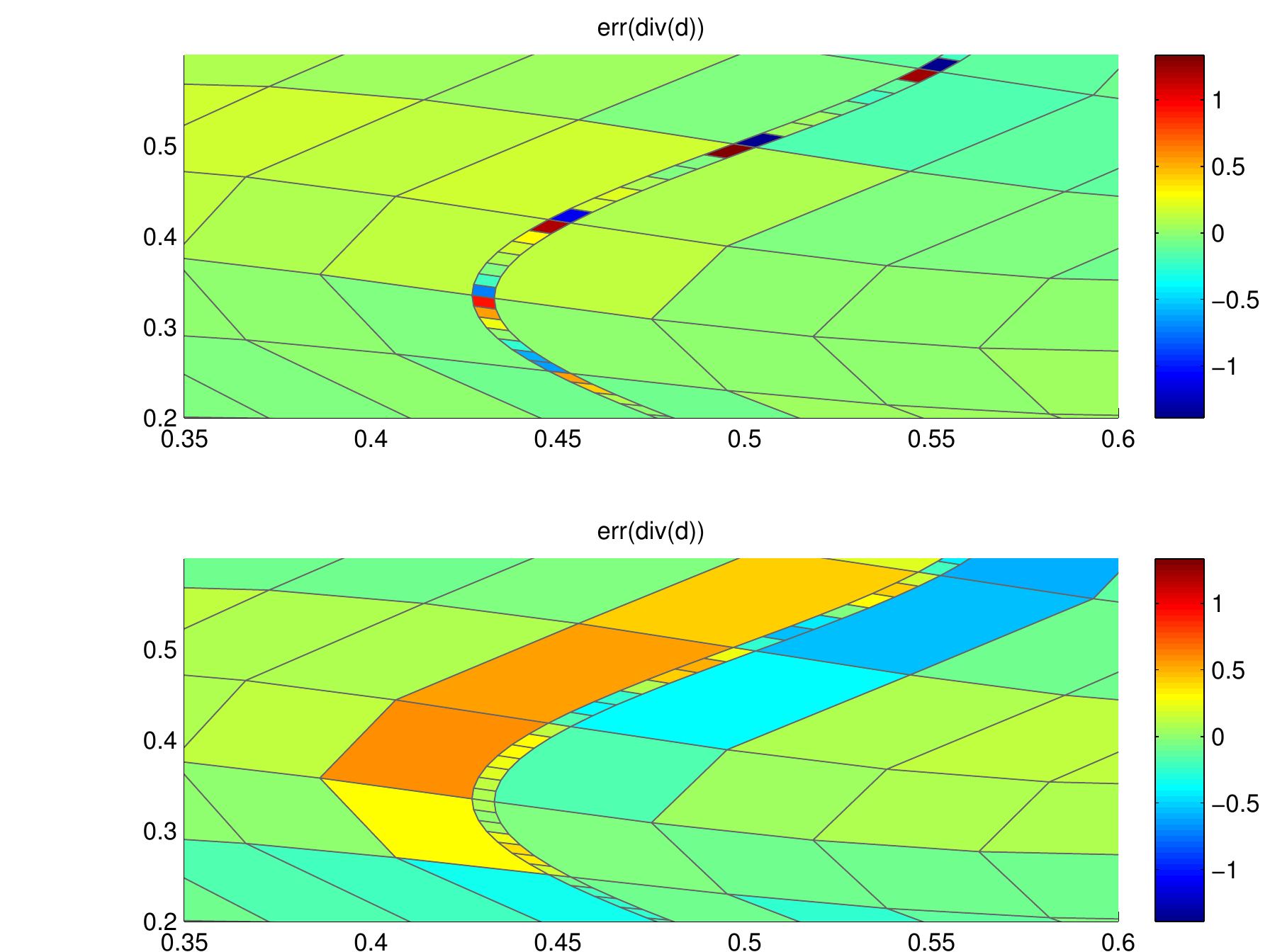}
  \caption{Plot of the divergence at the interface of the layer for Case
    \ref{case:layerRefinementTwisted} (top). We observe large oscillations for the VE
    method (middle) inside the layer but they remain confined in the
    layer. Comparatively, the oscillations are weaker for MPSA (bottom) but the error
    spreads more to the neighboring cells.}
  \label{fig:errdivzoomboxesv2eps7}
\end{figure}

In applications, the discretization level of the fracture layer may typically not
match the one of the matrix. We investigate this situation by setting up a test case
where the refinement in the layer is increased (Case
\ref{case:layerRefinementCartesian}\refstepcounter{mysubcase}\label{case:layerRefinementCartesian}). We
also consider the same case but we twist the grid (Case
\ref{case:layerRefinementTwisted}\refstepcounter{mysubcase}\label{case:layerRefinementTwisted})
to break eventual symmetry effects.  In Figure \ref{fig:convboxesv2eps7}, we can
observe the $L^\infty$-norm of the error for the stress grows significantly for the
VE methods. The error in divergence remains zero in the Cartesian case (Case
\ref{case:layerRefinementCartesian}) but, by looking at the twisted case (Case
\ref{case:layerRefinementTwisted}), we conclude that this is only due to a symmetry
effect. The error in displacement for this test case is plotted in Figure
\ref{fig:errdboxesv2eps7}. In Figure \ref{fig:stressfacesboxesv2eps7} where a plot of
the force is given at the interface, we observe that same oscillations for the VE
method as previously in the case of two adjacent regions with different
discretization levels (namely Case \ref{case:twoRegionsIsotropicRefinement}). Also in
this case, the MPSA method gives a smoother approximation closer to the analytical
solution. In Figure \ref{fig:errdivzoomboxesv2eps7}, we present a plot of the
divergence along the interface together with a zoom on the region around of the error
in divergence. Note that for the MPSA method, we use the finite volume definition of
the divergence, that is, the value of the divergence in the cell is equal the sum of
the normal component for each face. Again, we observe how the error in the VE method
remains highly localized and concentrates inside the layer while the error for the
MPSA methods spreads more to the coarse cell.

\clearpage

\subsection{Case \ref{case:locking}: Stability near incompressibility}
\label{sec:stabnearincomp}
\refstepcounter{mycase}\label{case:locking}

In the following experiments, we set up cases to test the method with respect to the
near incompressibility limit. The boundary conditions are zero displacement on the
left and right side and free displacement on the top and bottom side. The external
force is a constant volumetric vertical force, like for example a gravitational
force. We present test results for three grid types which highlights the main
features. We consider grids made of hexahedrons (Case
\ref{case:limincomphexa}\refstepcounter{mysubcase}\label{case:limincomphexa}),
triangles (Case
\ref{case:limincomptriangle}\refstepcounter{mysubcase}\label{case:limincomptriangle})
and quadrilaterals (Case
\ref{case:limincompquad}\refstepcounter{mysubcase}\label{case:limincompquad}). To
generate each grid, we start with a uniform tessellation. Then, the grid is twisted
in order to remove any side-effects that may arise from symmetry. In the case of
hexahedrons (Case \ref{case:limincomphexa}), we observe numerical locking for VEM and
MPSA while VEM-relax, VEM-relax-extra and MPSA-relax-extra give a good approximation
of the solution. The numerical locking for VEM can be observed directly in the
displacement field where we can see that the medium appears to be much stiffer than
it actually is, hence the name of \textit{locking}. The numerical locking for MPSA is
more visible in the divergence field as artificial oscillations. Note that we could
have chosen the parameter $\nu$ closer to 0.5 where the effects of numerical locking
lead to a completely different solution, but we prefered to show examples where we
can both recognize the solution and see the beginning of the perturbations caused by
locking. In a grid made of hexahedrons, the ratio between the number of nodes and
number of cells is typically large. As discussed in Section \ref{subsubsec:locking},
such configuration is favorable for the VE method and disadvantageous for the MPSA
method. Each node does not belong to more than three edges and the stability
condition given in \cite{da2014mimetic} is fulfilled. Therefore, the solution given
by VEM-relax is free from locking effect and VEM-relax-extra does not bring any
improvement. Let us now consider the grid made of triangles (Case
\ref{case:limincomptriangle}). In this case, we observe that numerical locking
affects the VEM and VEM-relax methods but MPSA, MPSA-relax-extra and VEM-relax-extra
remain unaffected. In the case of a triangular grid the ratio between the number of
nodes and number of cells is typically small and equal to the inverse of the same
ratio for a hexahedral grid. Then, we explained in section \ref{subsubsec:locking}
why this situation favors MPSA and penalizes VEM. This analysis is consolidated by
the numerical results. We observe that the relaxation of the VE method is not enough
to get rid of the numerical locking and extra degrees of freedom are required. We
also note that the solution obtained from MPSA-relax-extra contains slightly more
oscillations and is not as good as the one obtained by standard MPSA. This result
highlights the relaxation effect of the method, which means that, if numerical
locking is not an issue, the standard MPSA is expected to give less error than
MPSA-relax-extra. In the setting of VEM, it corresponds to the fact that, when there
is no locking, VEM has in general a smaller error constant than VEM-relax, as the
latter considers a less accurate approximation of the divergence, see
\eqref{eq:akhookapprox} compared to \eqref{eq:akhook2}. Finally, we consider a grid
made of quadrilaterals (Case \ref{case:limincompquad}). In this case the ratio
between the number of nodes and number of cells is close to one, so that neither the
VEM or MPSA method is a priori favored. We observe that the VEM and VEM-relax methods
suffer from locking. In comparison, the MPSA handles remarkably well this case and
does not present any sign of locking. As predicted by the theory, VEM-relax-extra is
free from locking.

\begingroup
\def\imagewidth{0.3\textwidth}
\def\insertimage #1{\parbox[c][][c]{\imagewidth}{\centering\includegraphics[width=\imagewidth]{#1}}}
\def\insertcolorbar #1{\parbox[c][][c]{\imagewidth}{\centering\hspace*{1mm}\rotatebox{90}{\includegraphics[height=\imagewidth]{#1}}}}
\def\insertlegend #1{\parbox[c][][c]{\imagewidth}{\centering #1}}
\def\insertcase #1{\parbox[c][][c]{0.13\textwidth}{\centering #1}}
\begin{figure}[h]
  \centering
  \begin{tabular}{@{}c@{}c@{}c@{}c@{}c@{}c@{}}

    &
    \insertlegend{$d_x$}&
    \insertlegend{$d_y$}&
    \insertlegend{$\dive\vec{u}$}\\

    \insertcase{VEM}&
    \insertimage{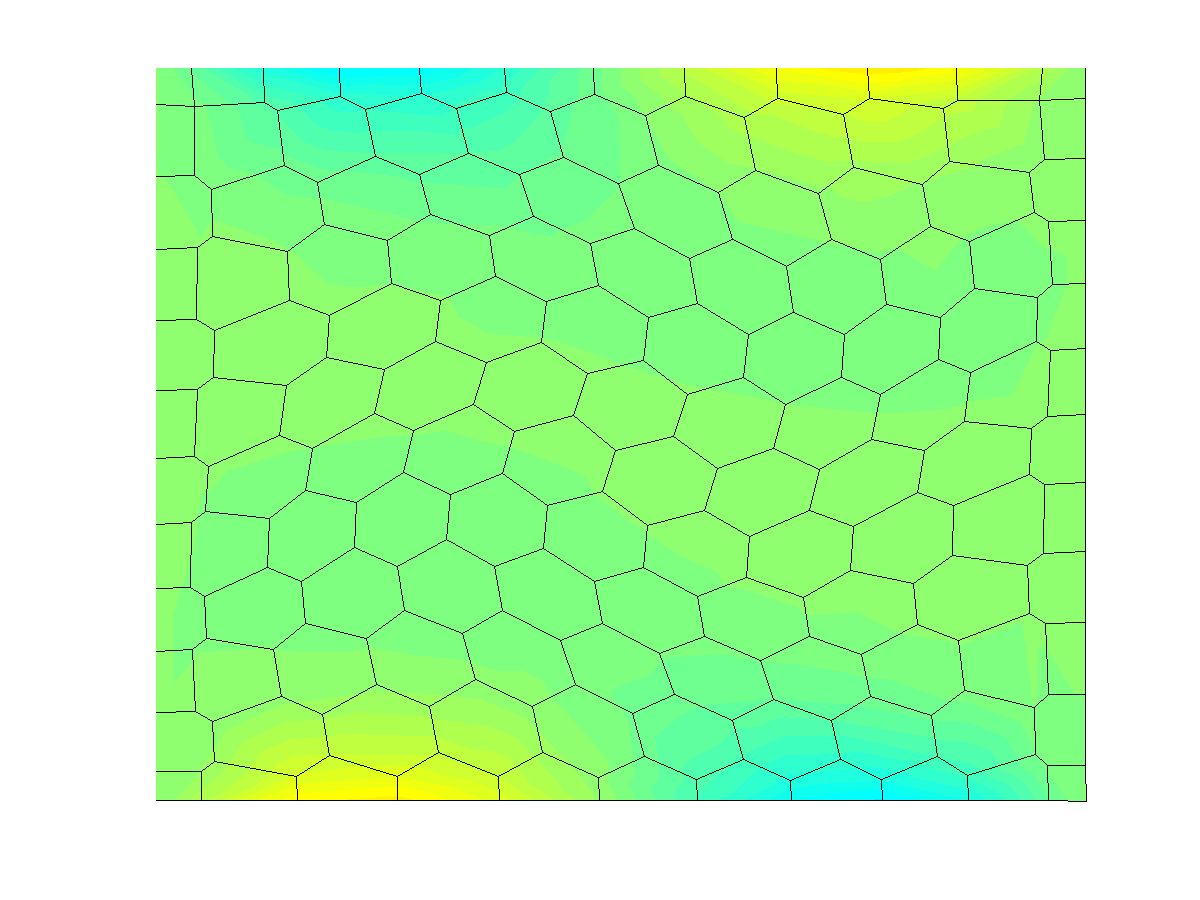}&
    \insertimage{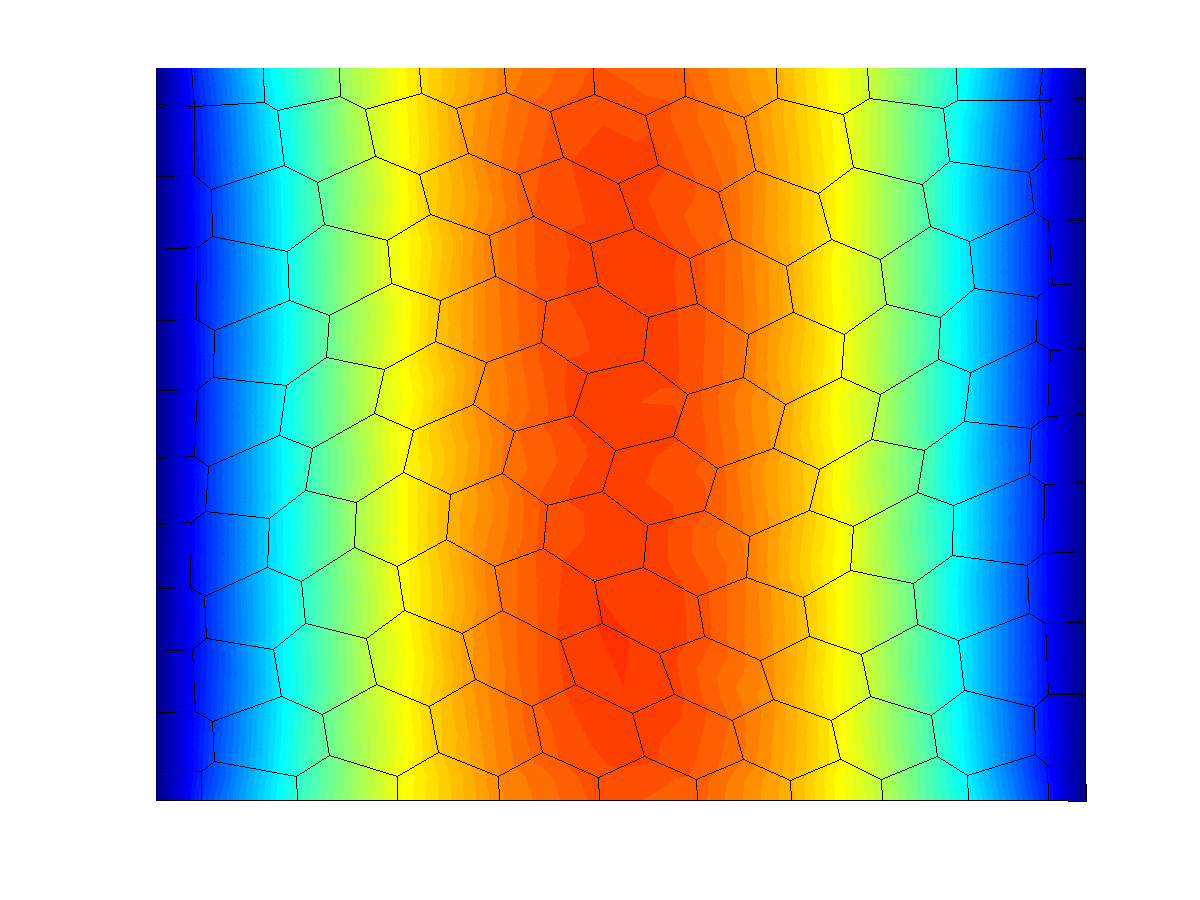}&
    \insertimage{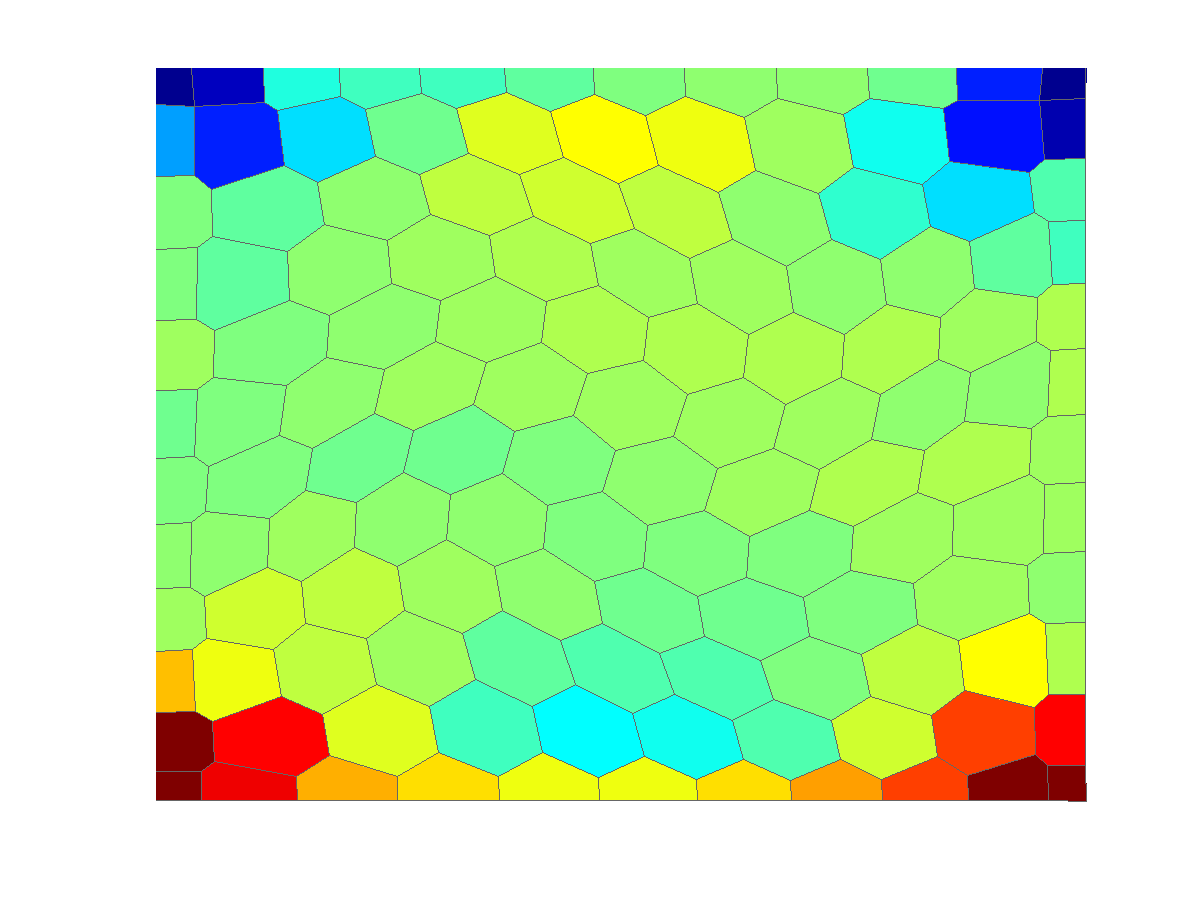}\\

    \insertcase{MPSA}&
    \insertimage{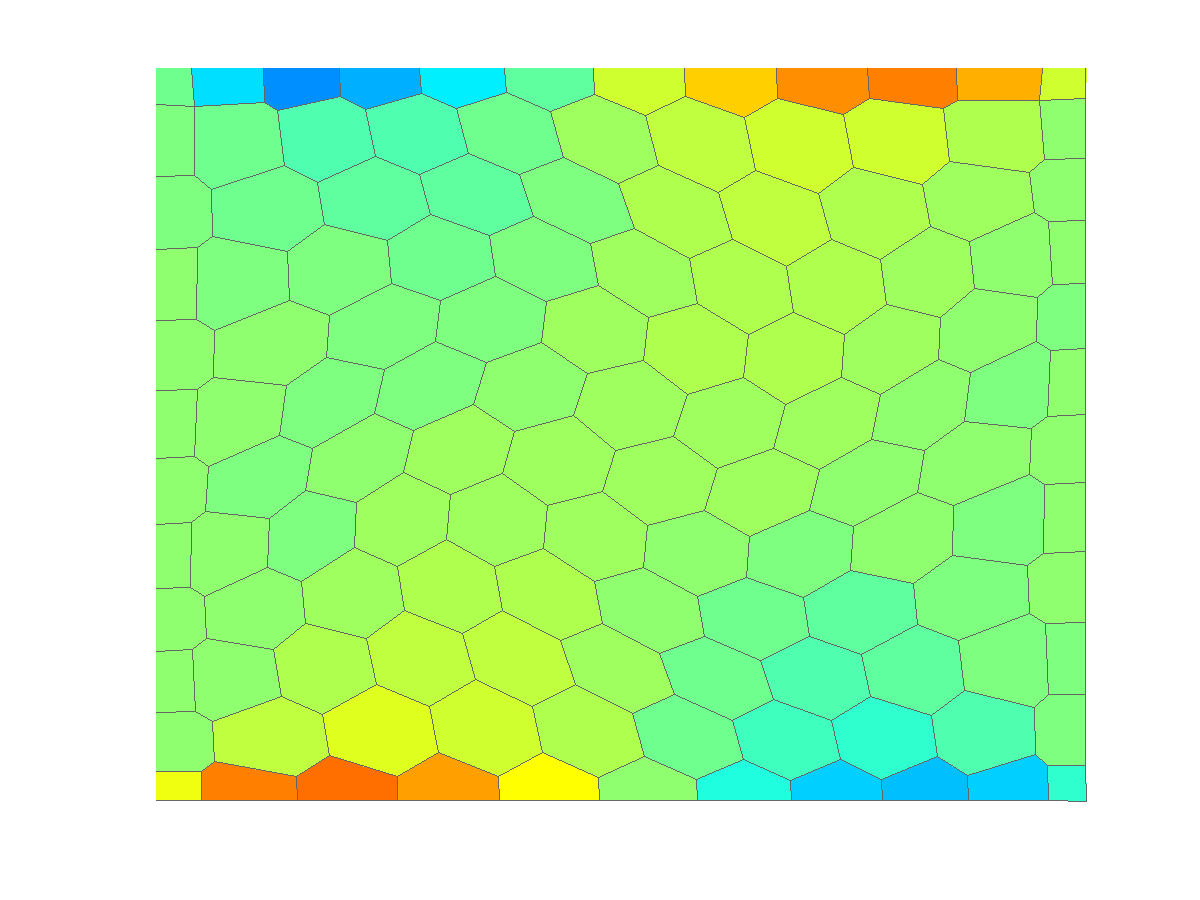}&
    \insertimage{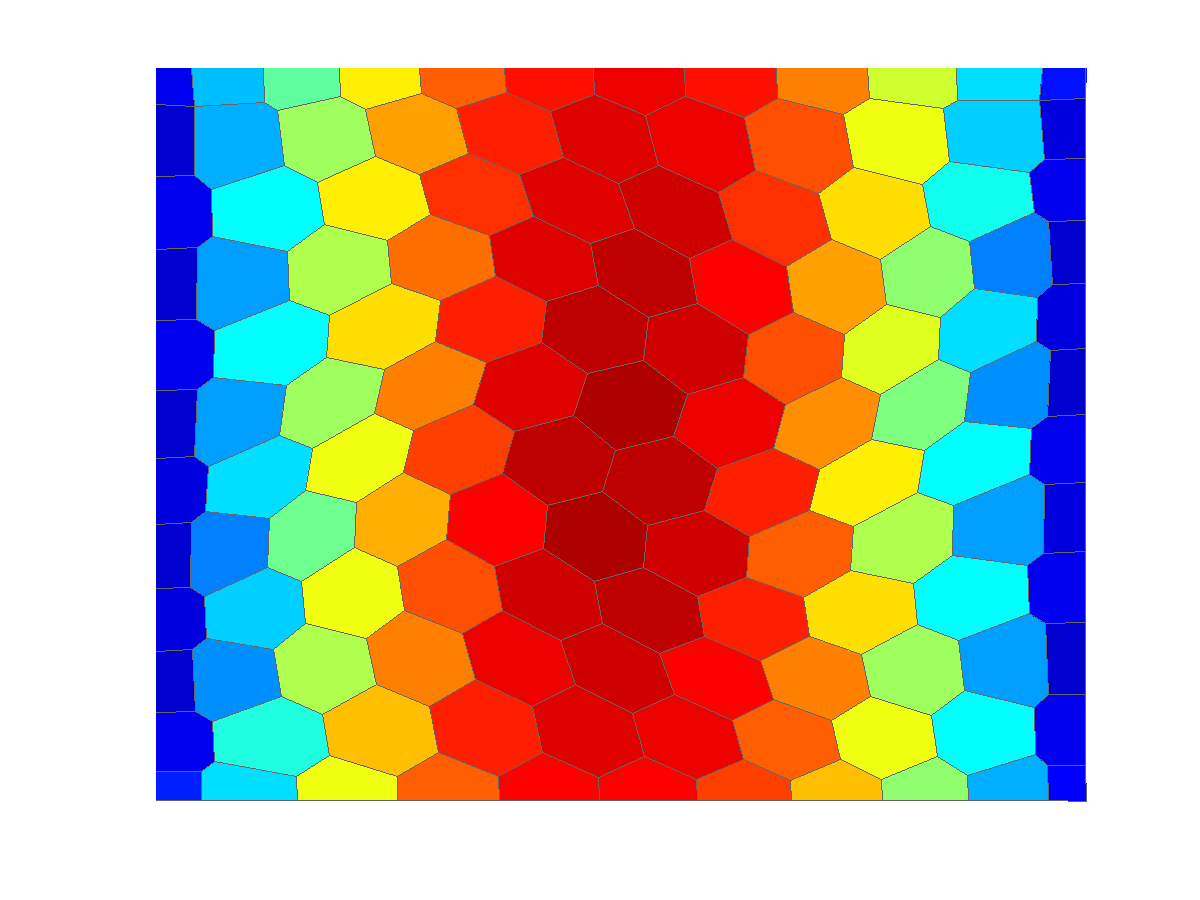}&
    \insertimage{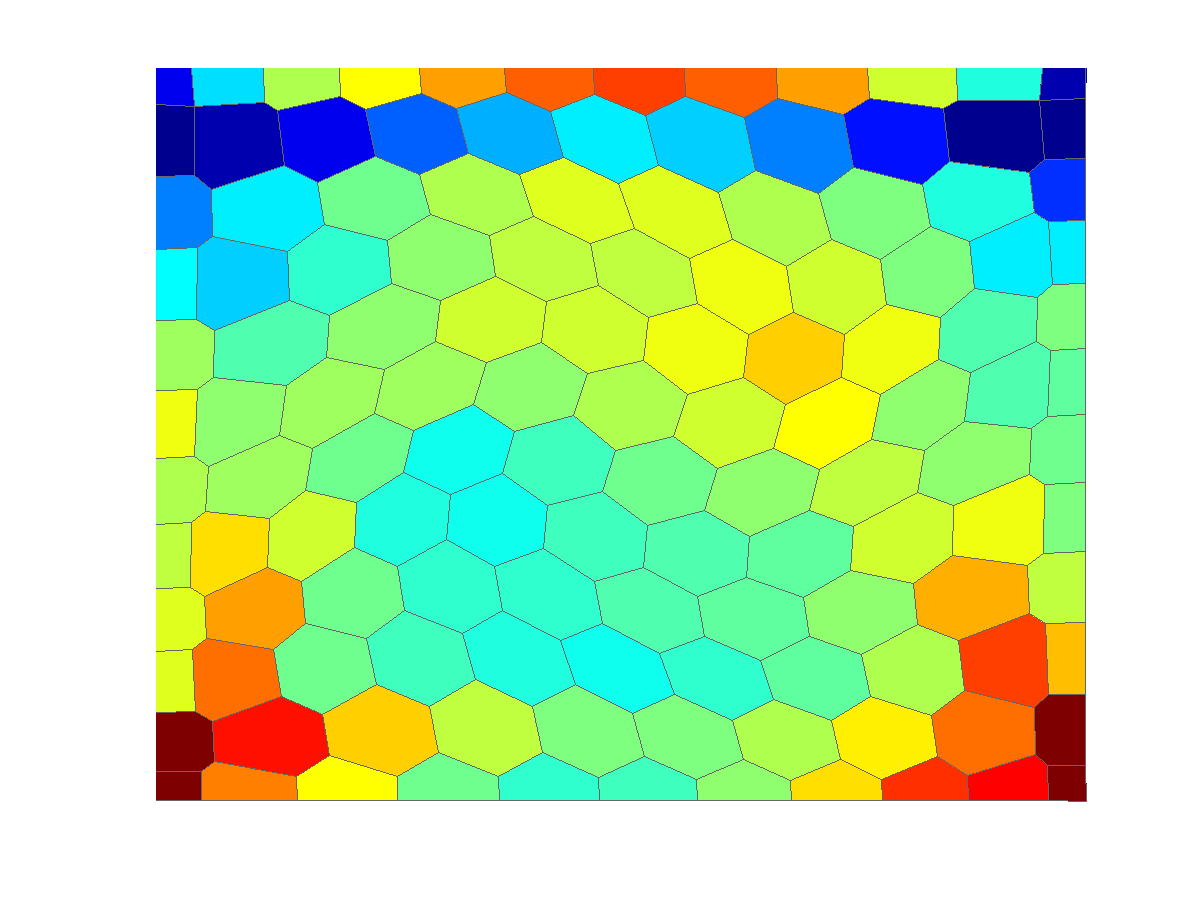}\\

    \insertcase{VEM relax}&
    \insertimage{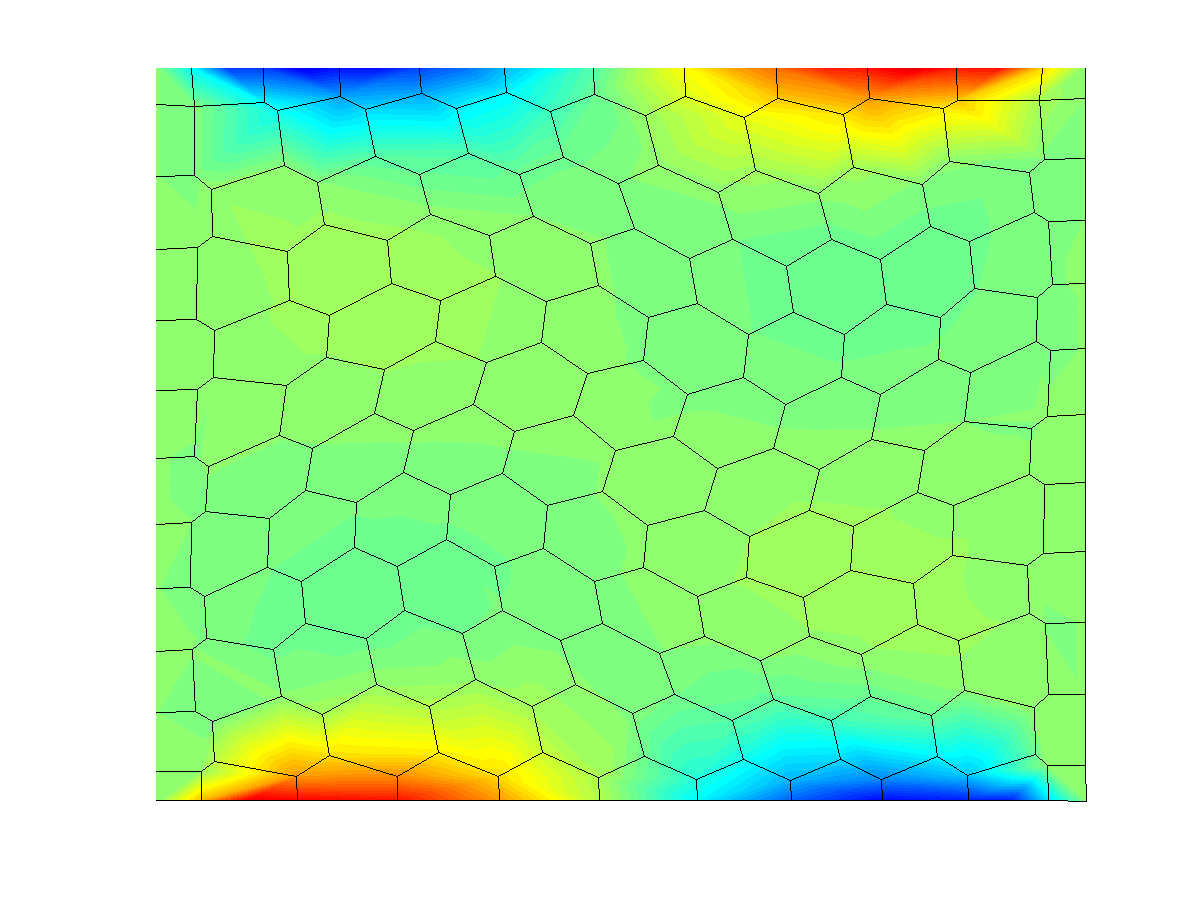}&
    \insertimage{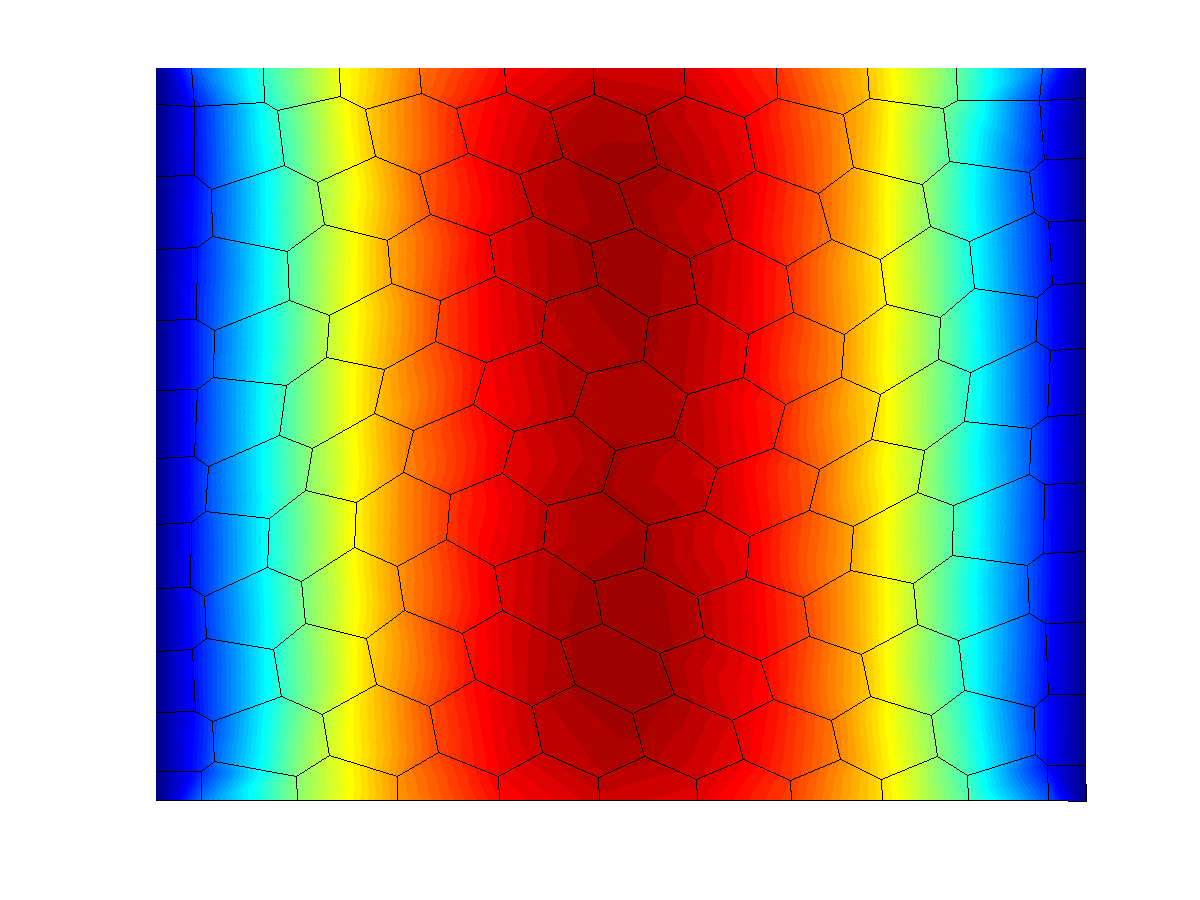}&
    \insertimage{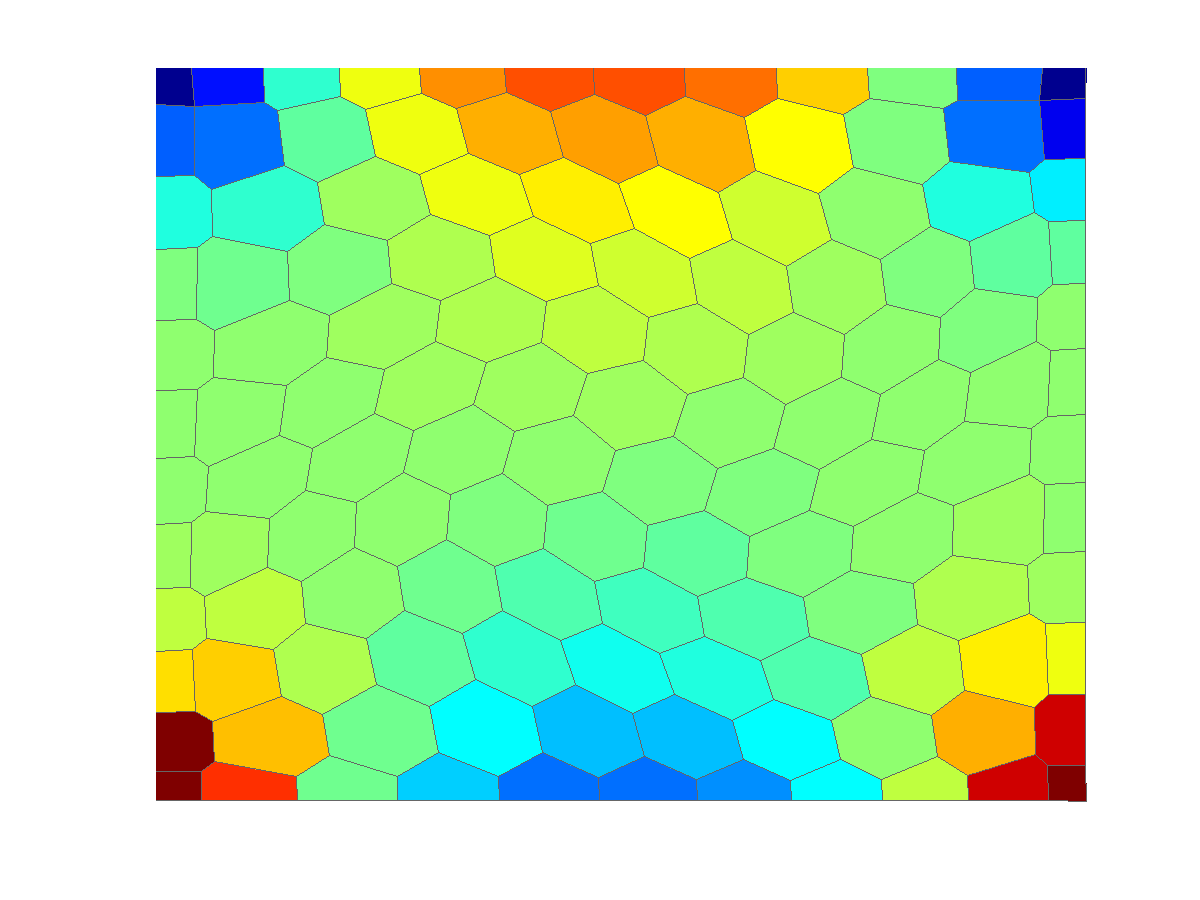}\\

    \insertcase{VEM relax extra}&
    \insertimage{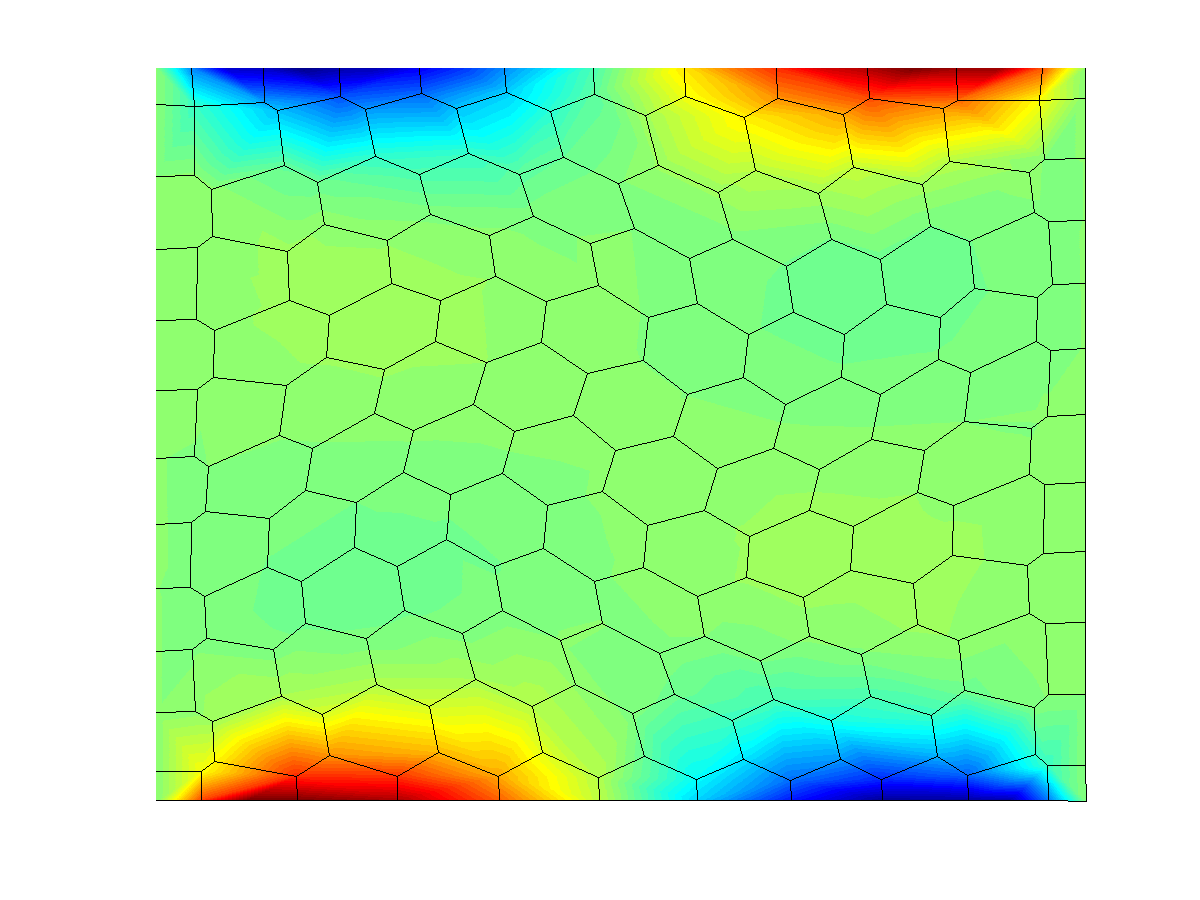}&
    \insertimage{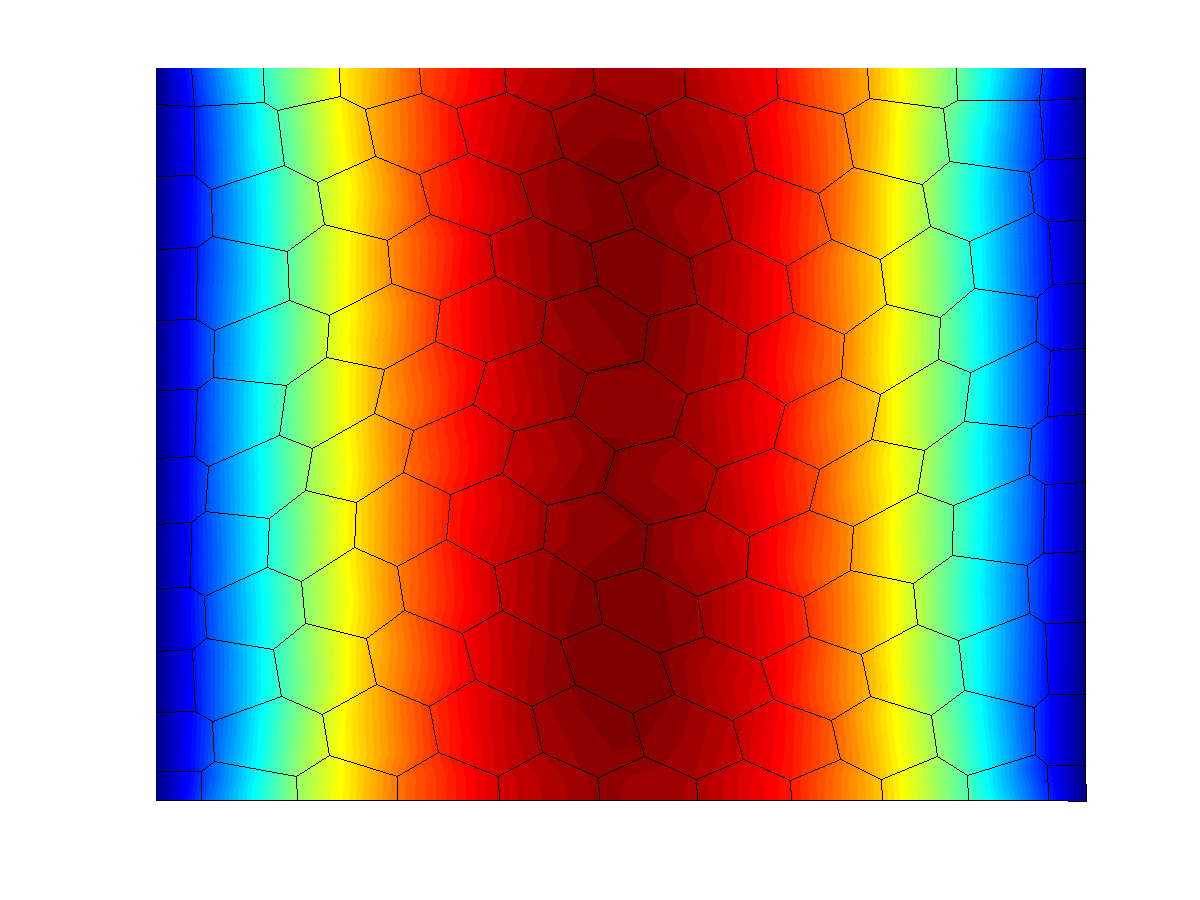}&
    \insertimage{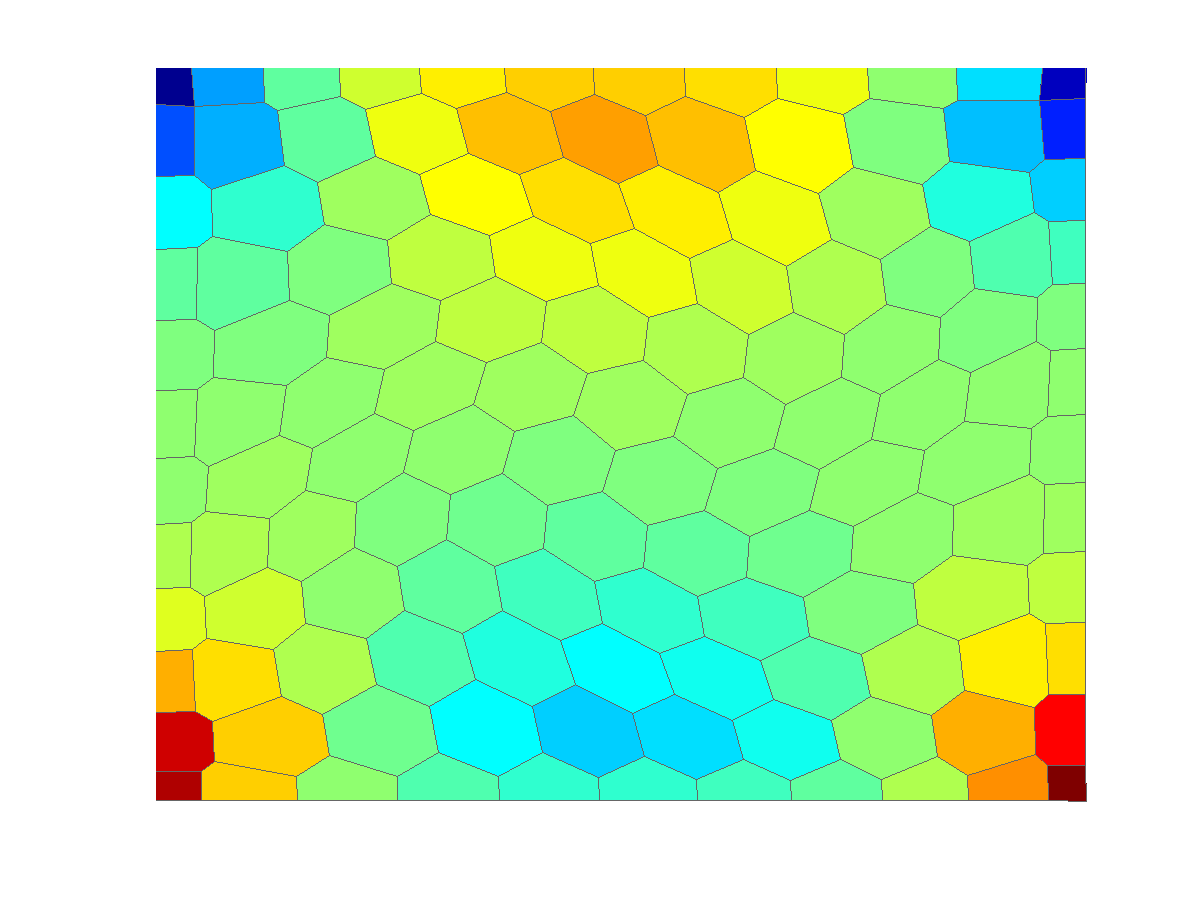}\\

    \insertcase{MPSA relax extra}&
    \insertimage{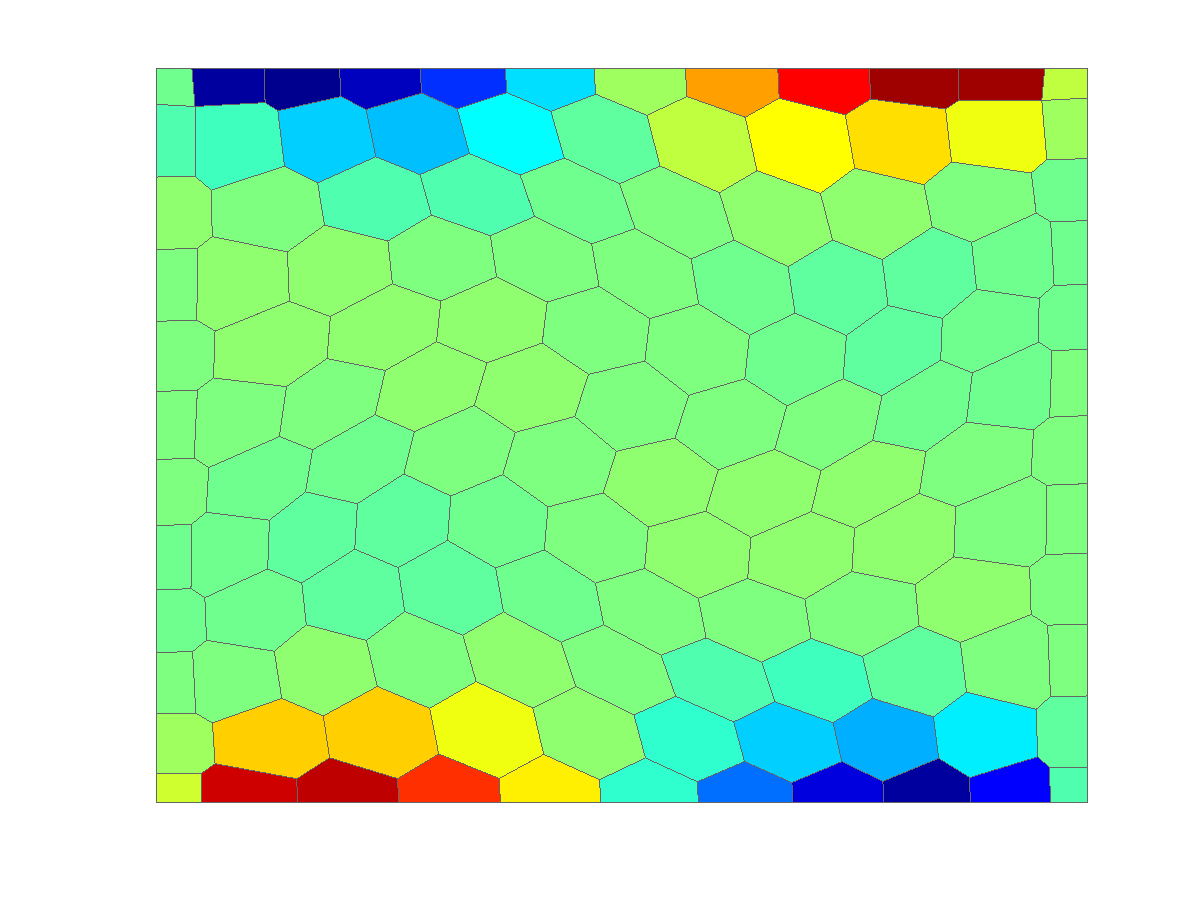}&
    \insertimage{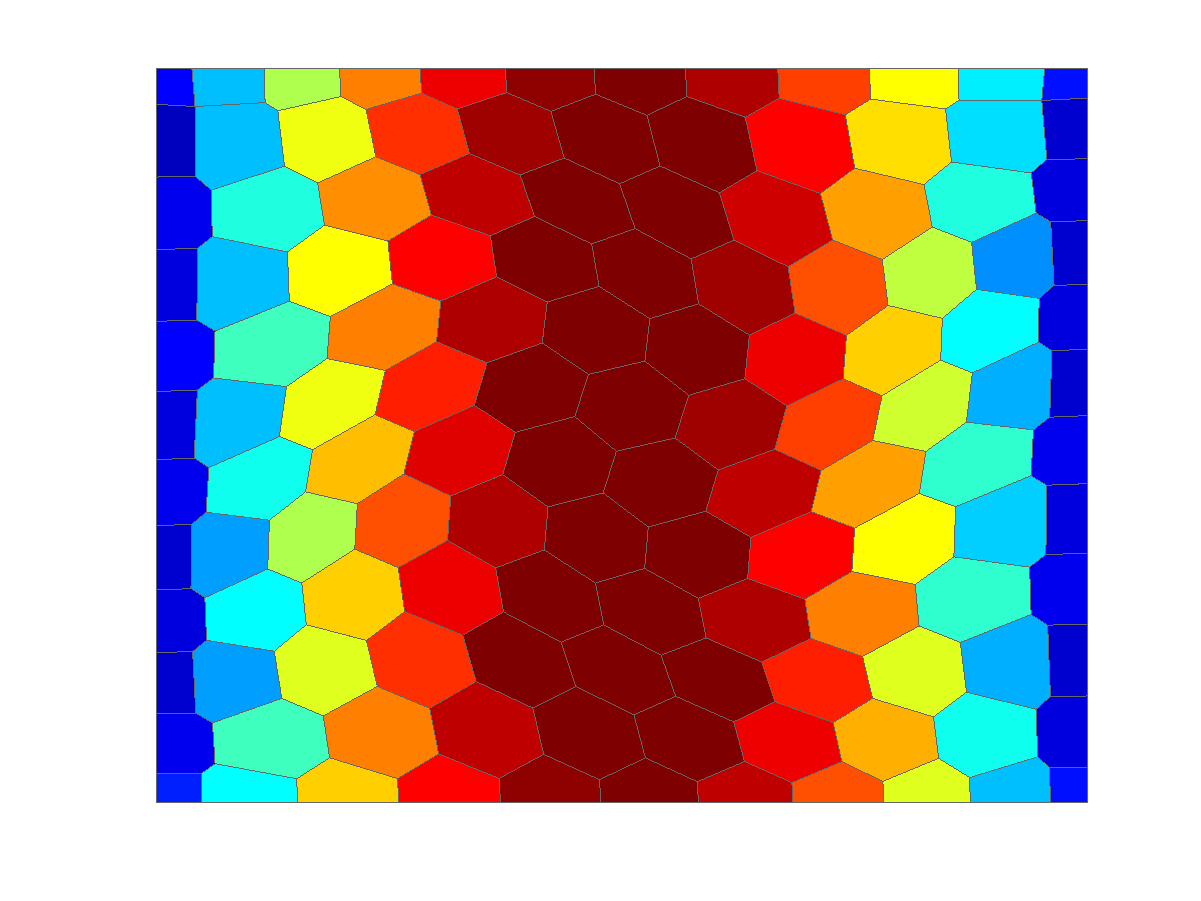}&
    \insertimage{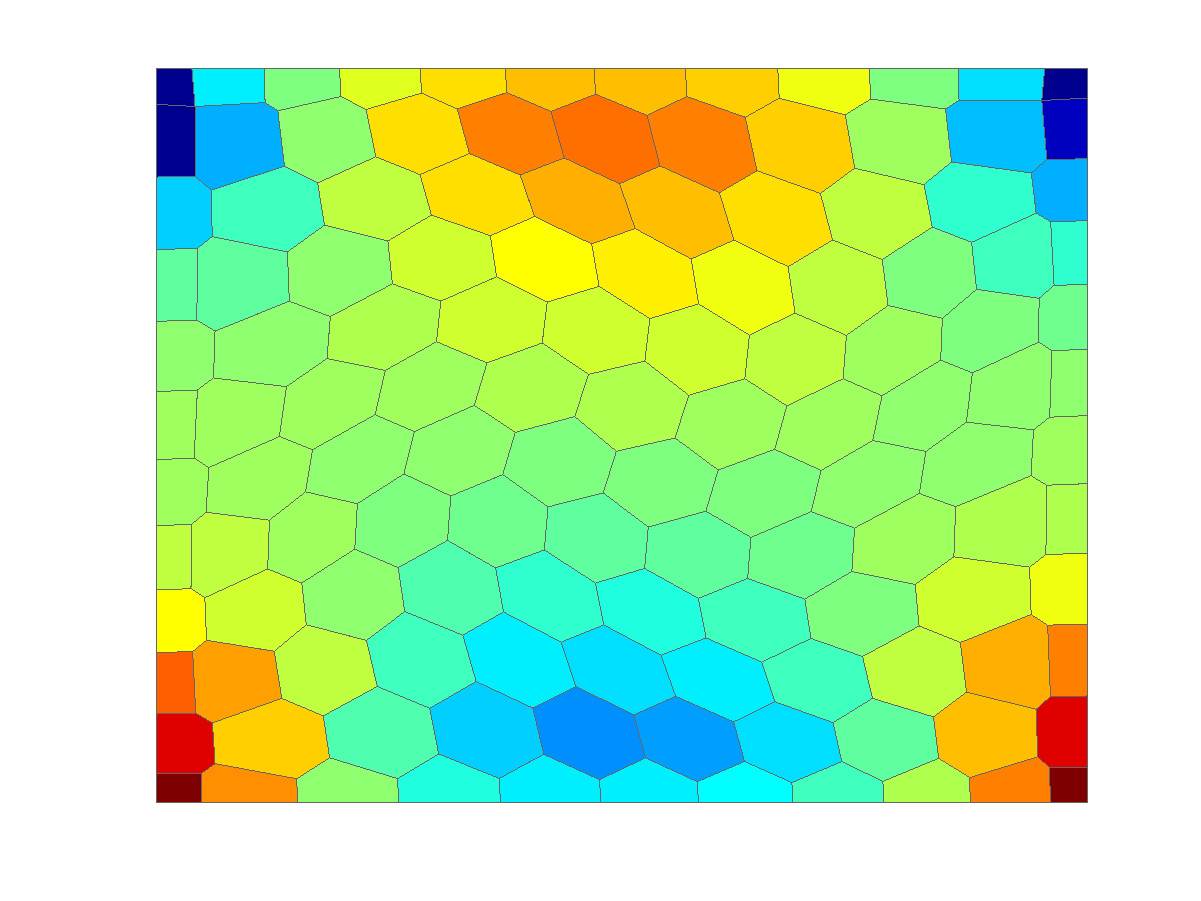}\\

    &
    \insertcolorbar{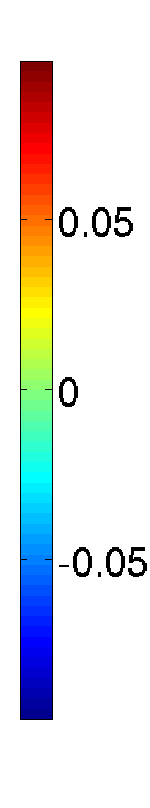}&
    \insertcolorbar{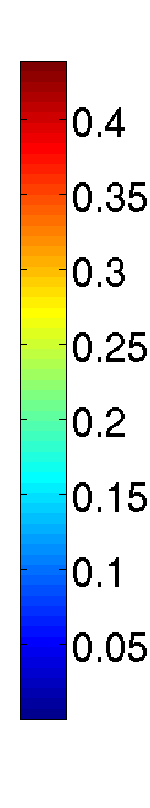}&
    \insertcolorbar{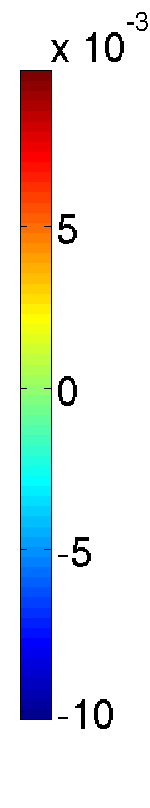}

  \end{tabular}
  \caption{Hexahedral grid and $\nu=0.495$ (Case \ref{case:limincomphexa}). MPSA and VEM suffers from numerical
    locking. VEM-relax and VEM-relax-extra are free from numerical locking.}
\end{figure}

\begin{figure}[h]
  \centering
  \begin{tabular}{c@{}c@{}c@{}c@{}c@{}c}
    &
    \insertlegend{$d_x$}&
    \insertlegend{$d_y$}&
    \insertlegend{$\dive\vec{u}$}\\

    \insertcase{VEM}&
    \insertimage{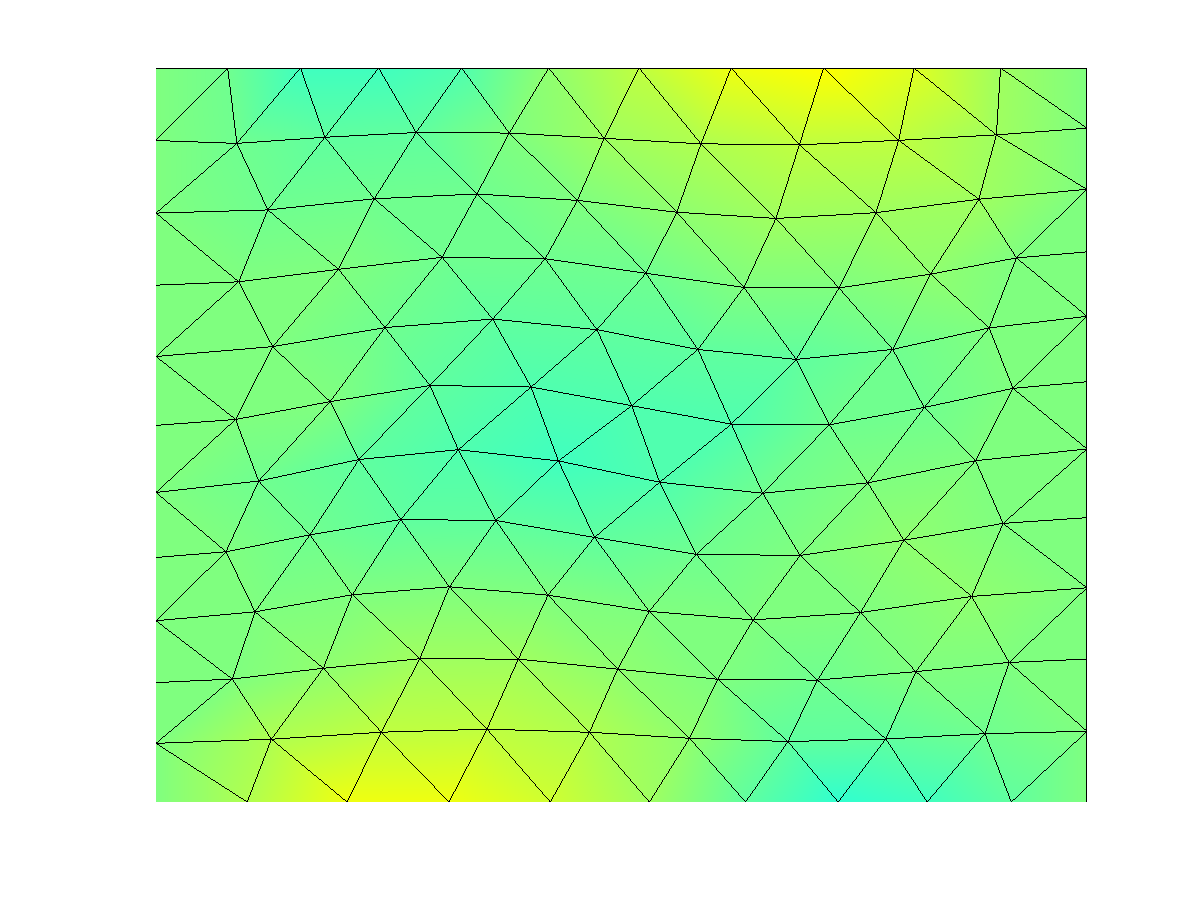}&
    \insertimage{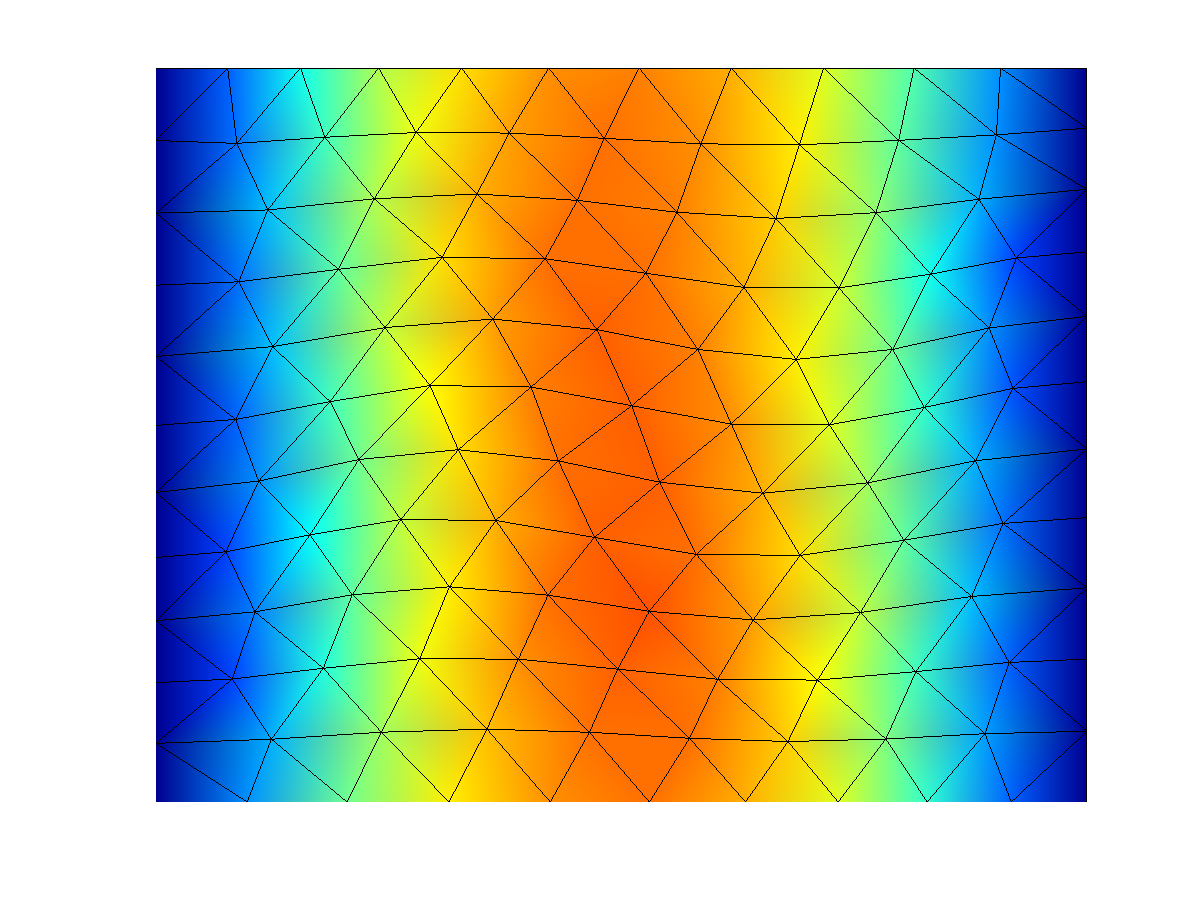}&
    \insertimage{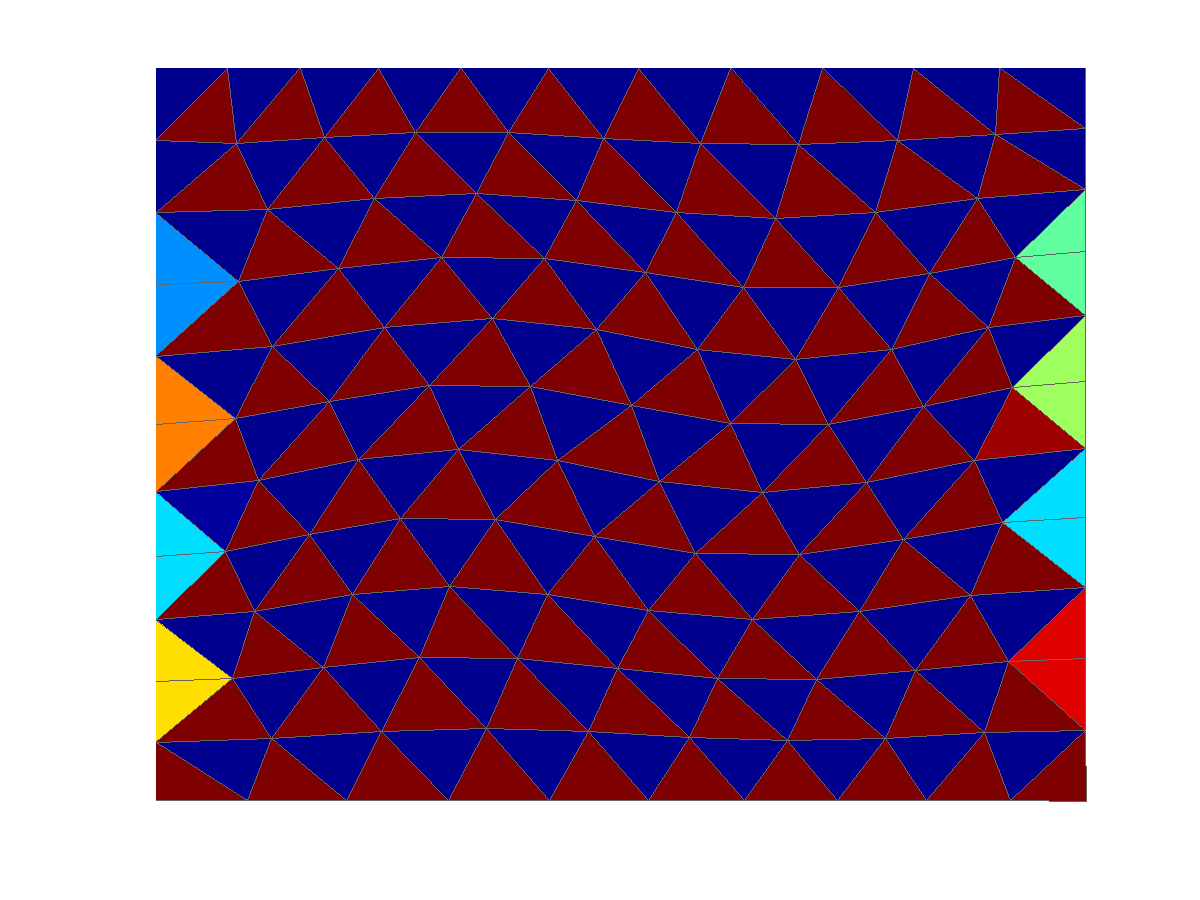}\\

    \insertcase{MPSA}&
    \insertimage{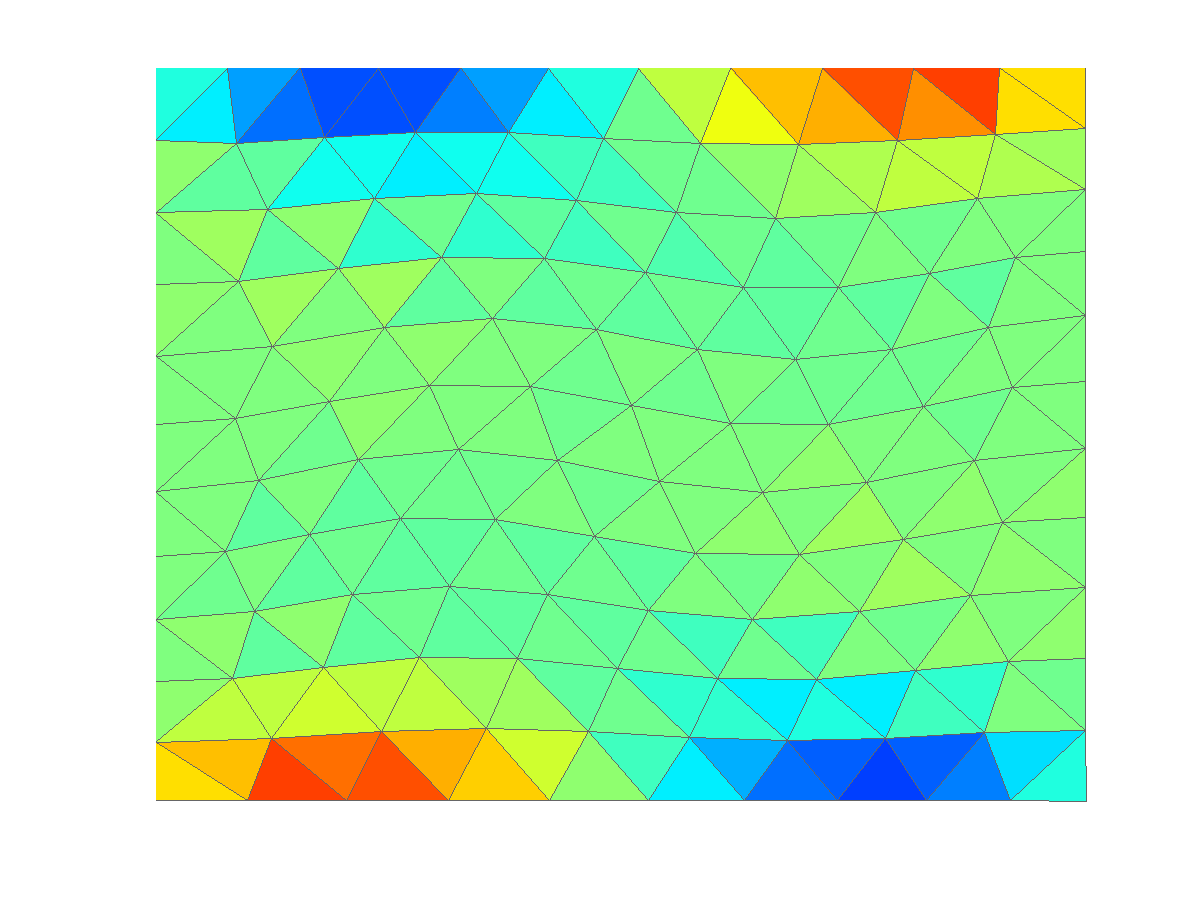}&
    \insertimage{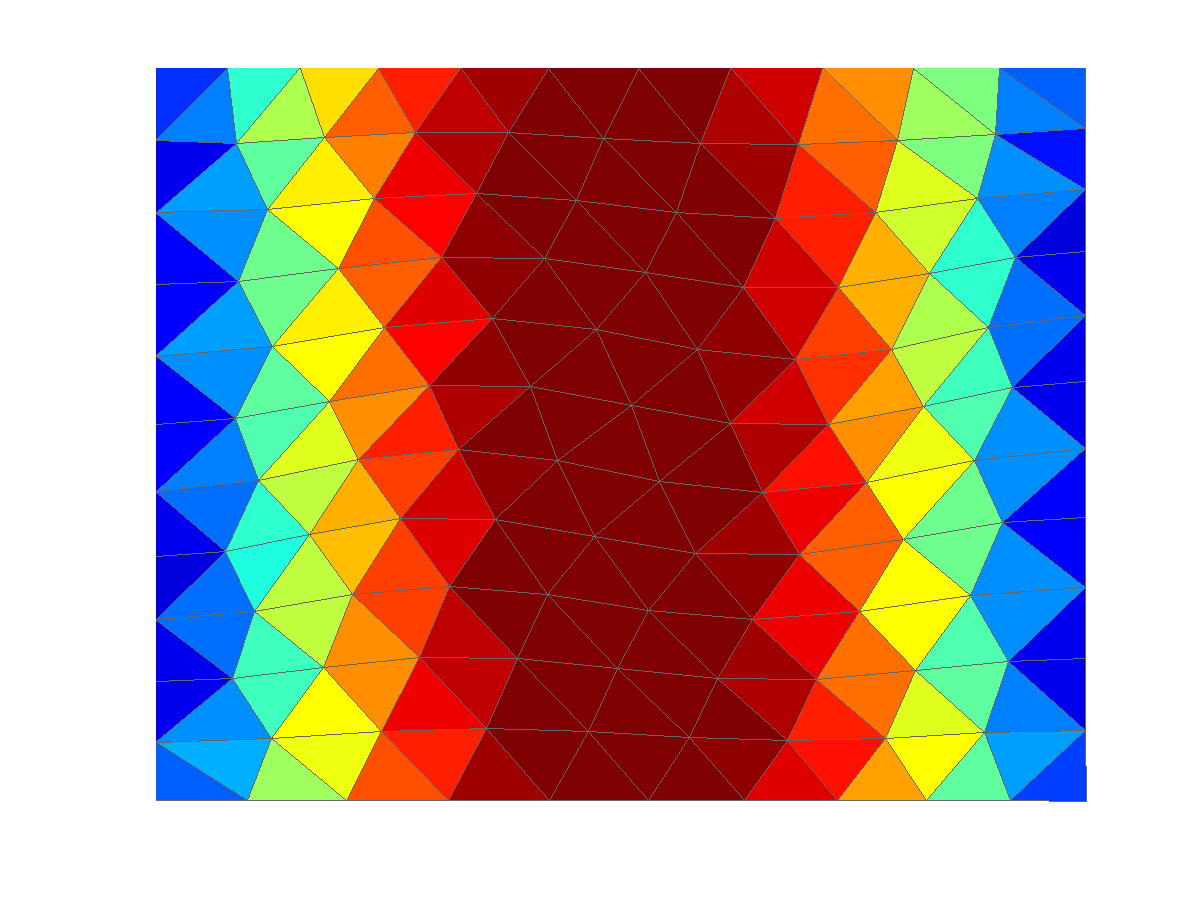}&
    \insertimage{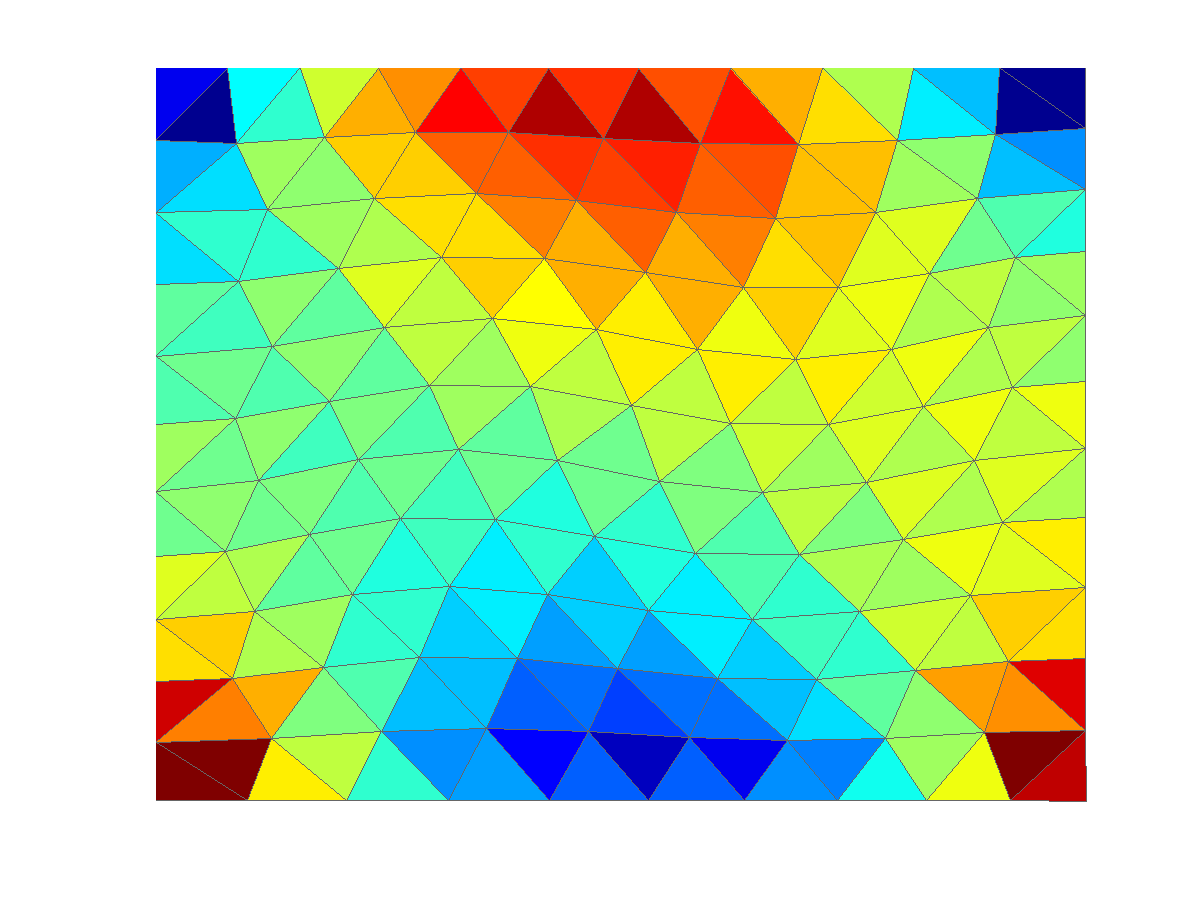}\\

    \insertcase{VEM relax}&
    \insertimage{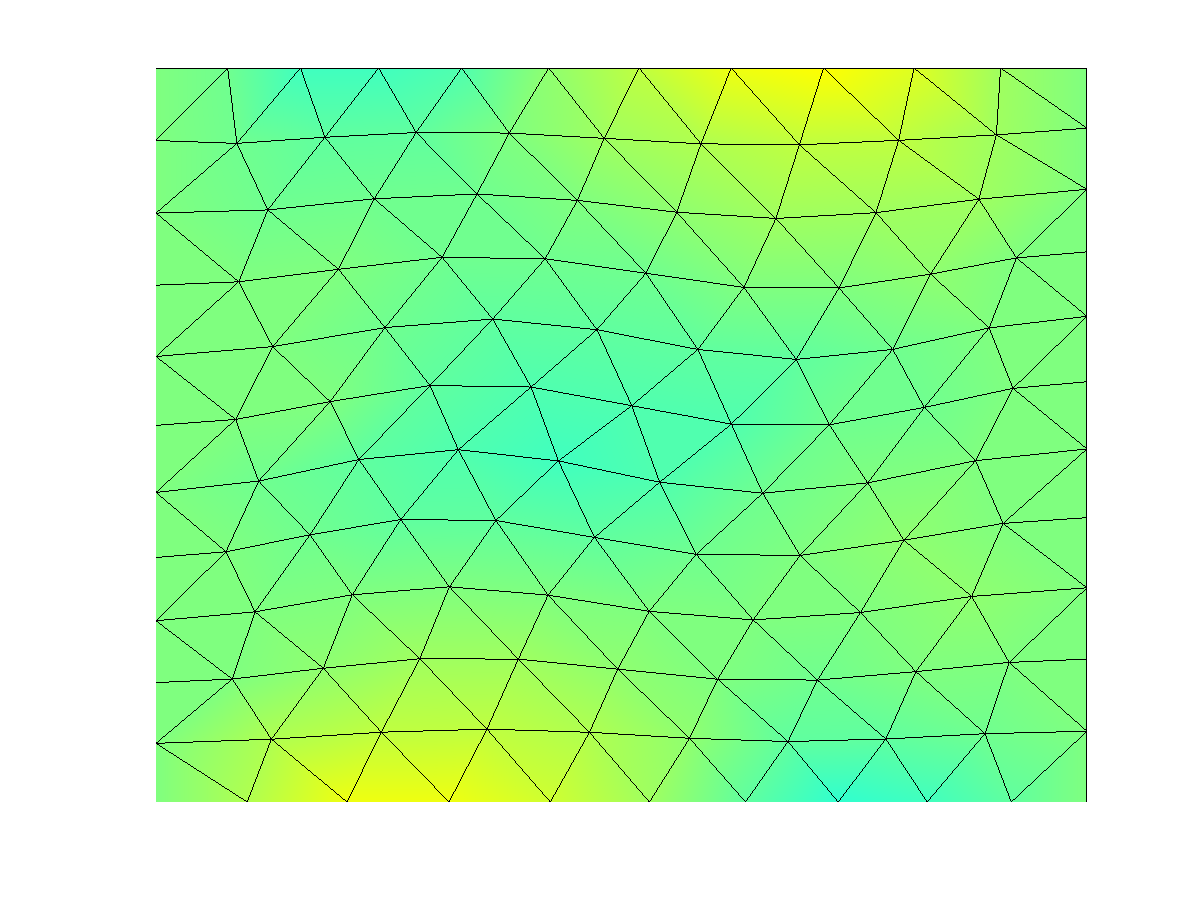}&
    \insertimage{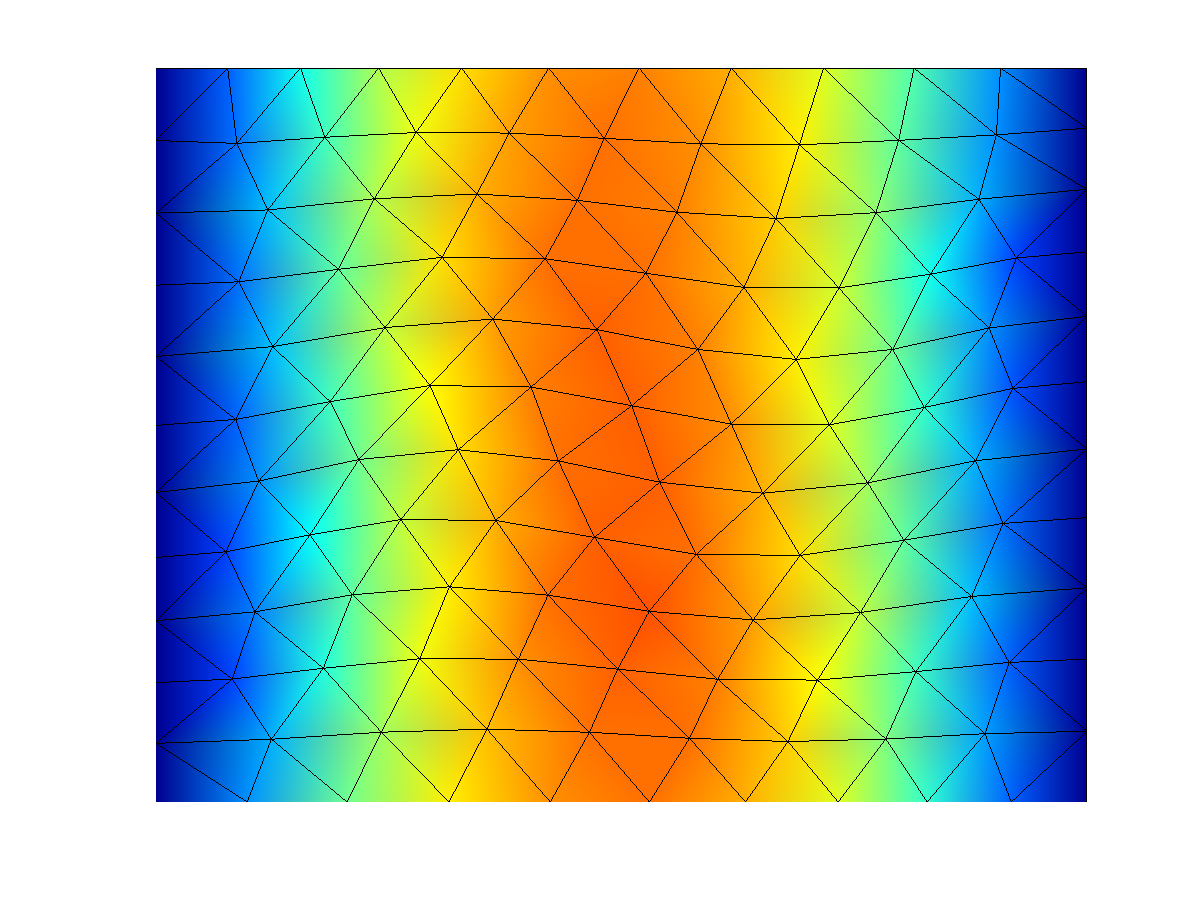}&
    \insertimage{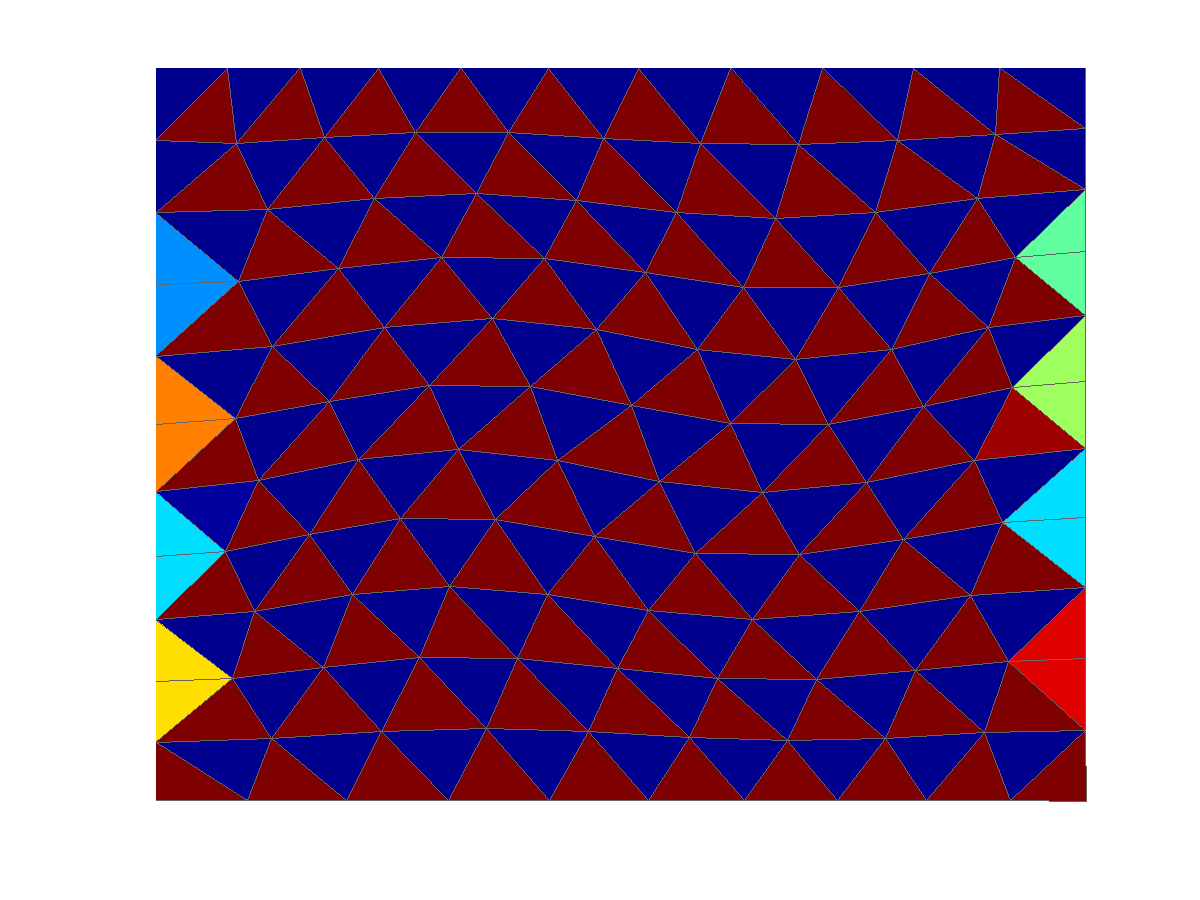}\\

    \insertcase{VEM relax extra}&
    \insertimage{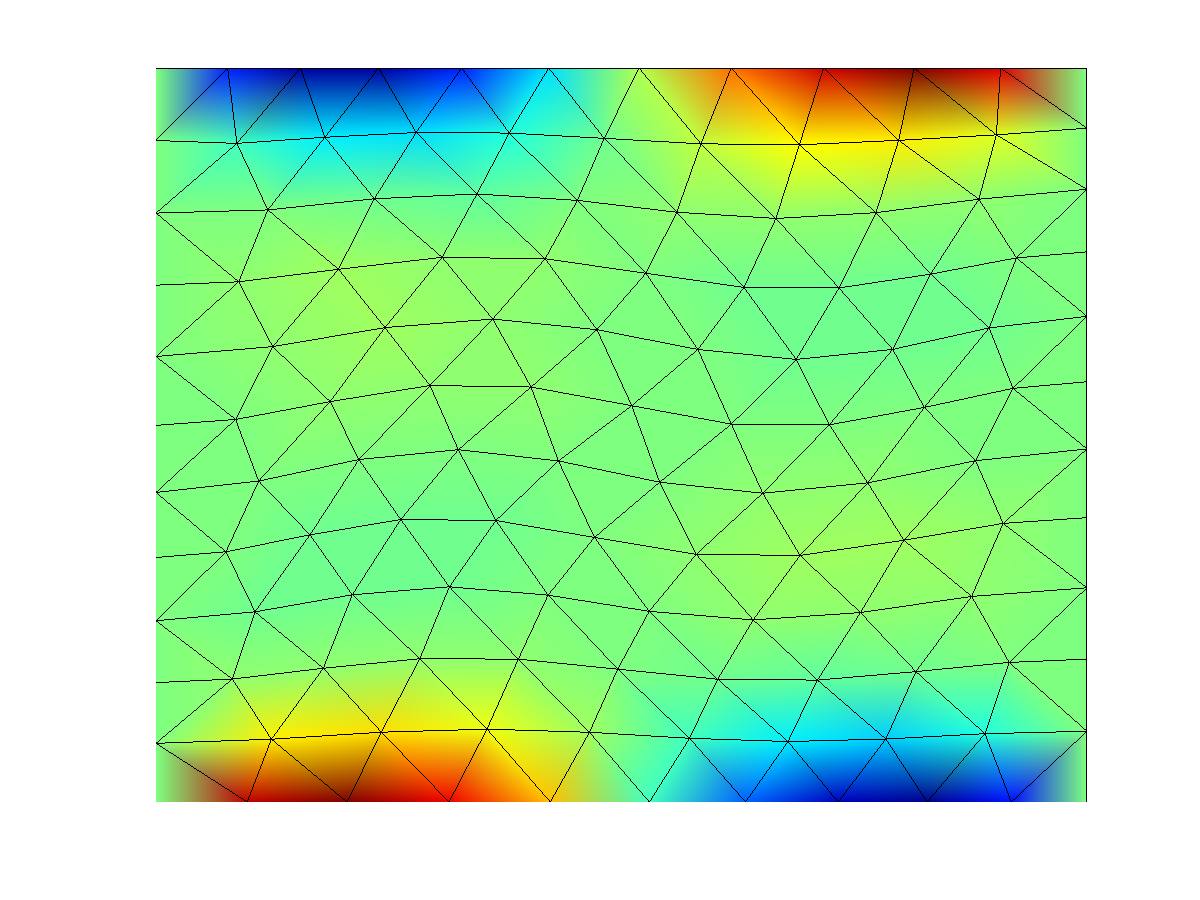}&
    \insertimage{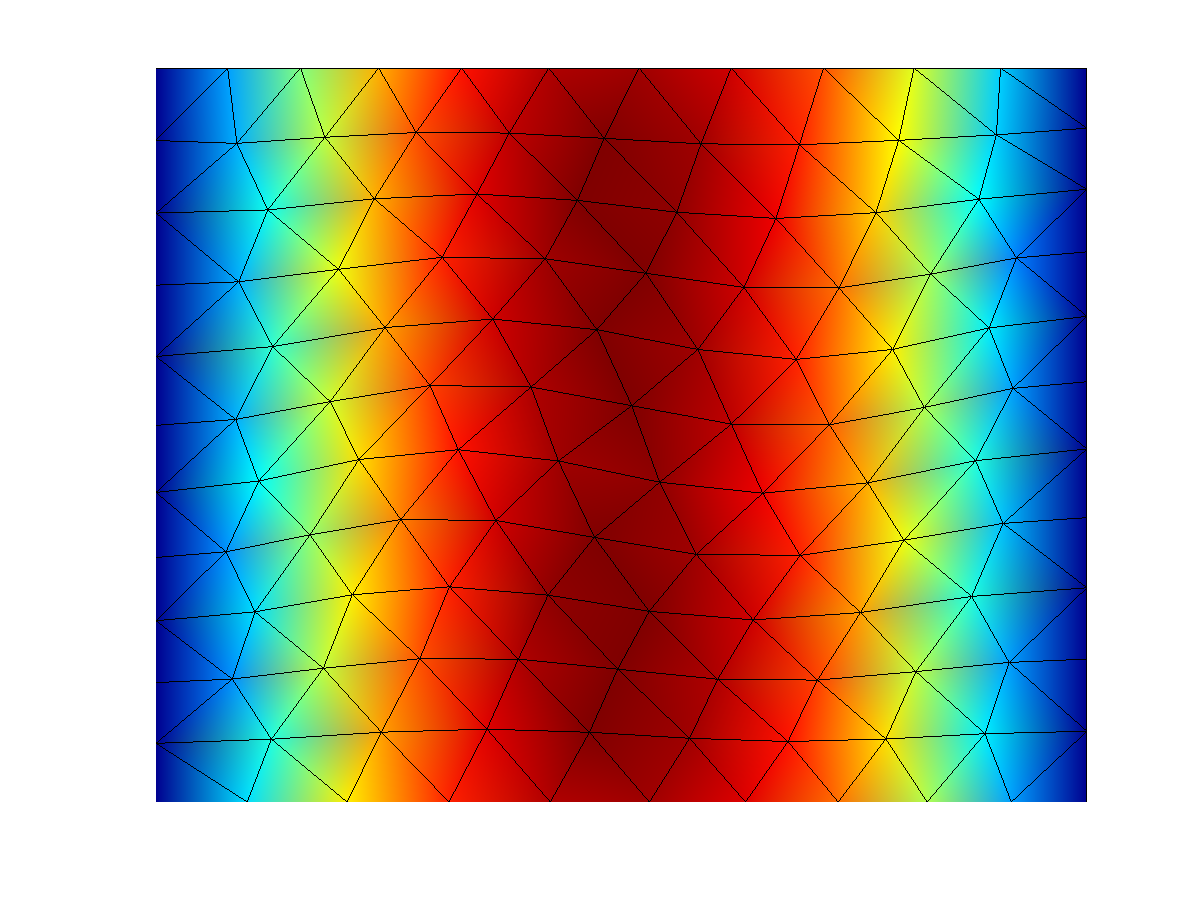}&
    \insertimage{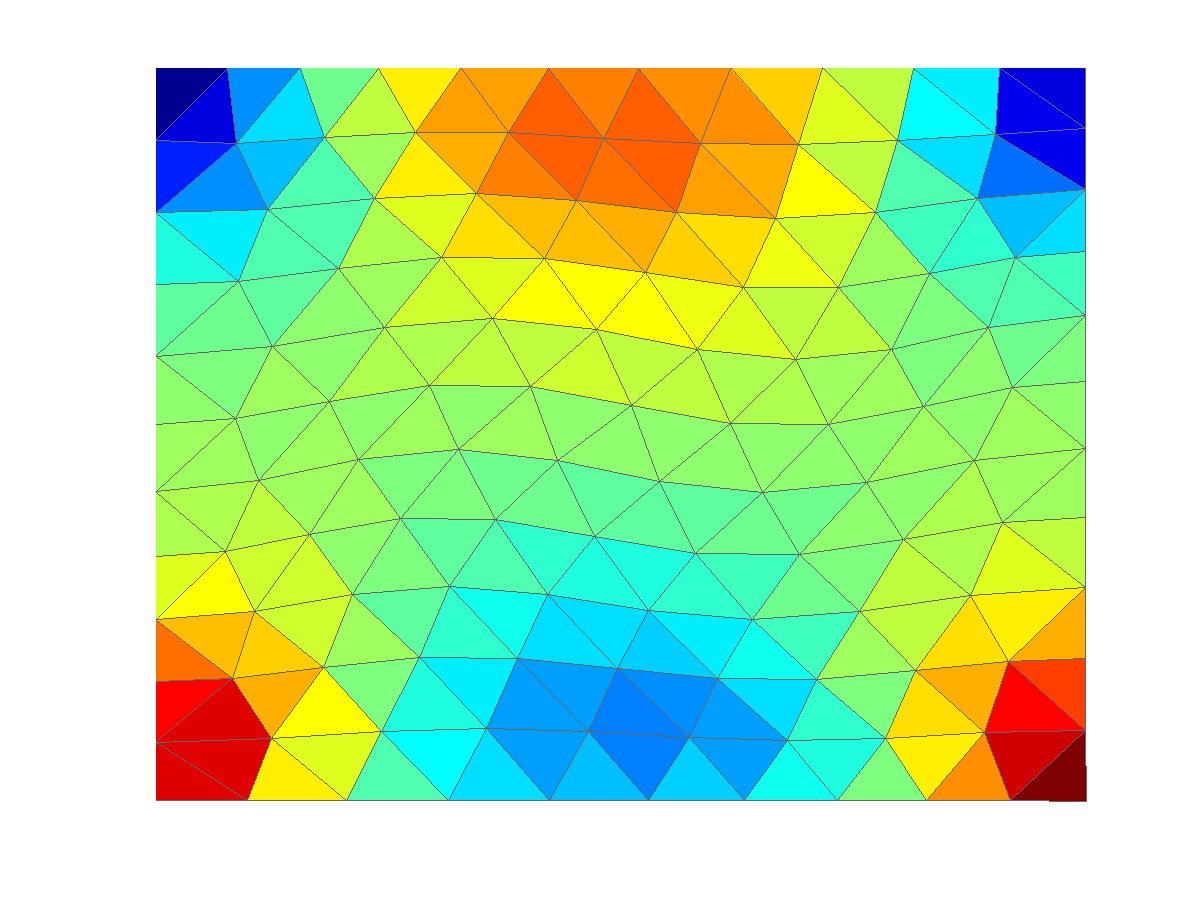}\\

    \insertcase{MPSA relax extra}&
    \insertimage{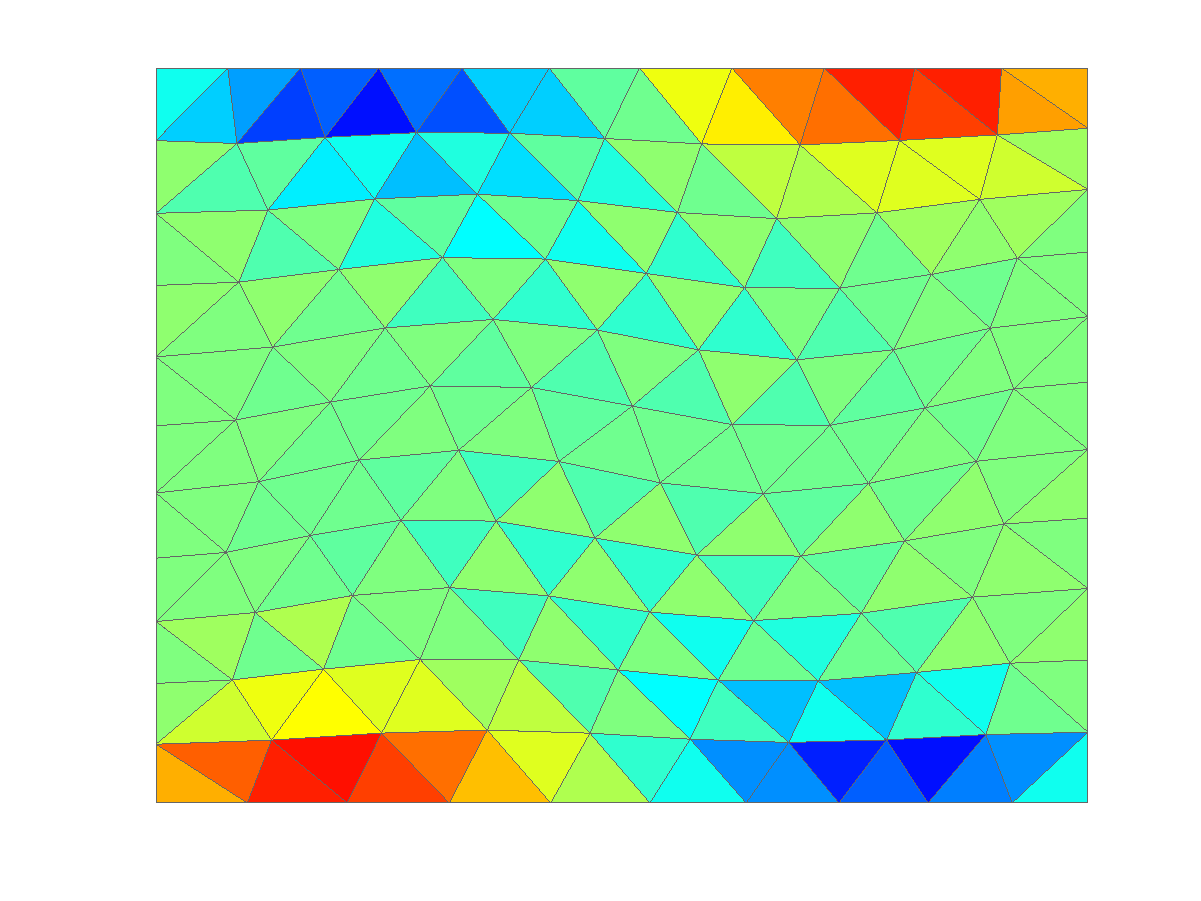}&
    \insertimage{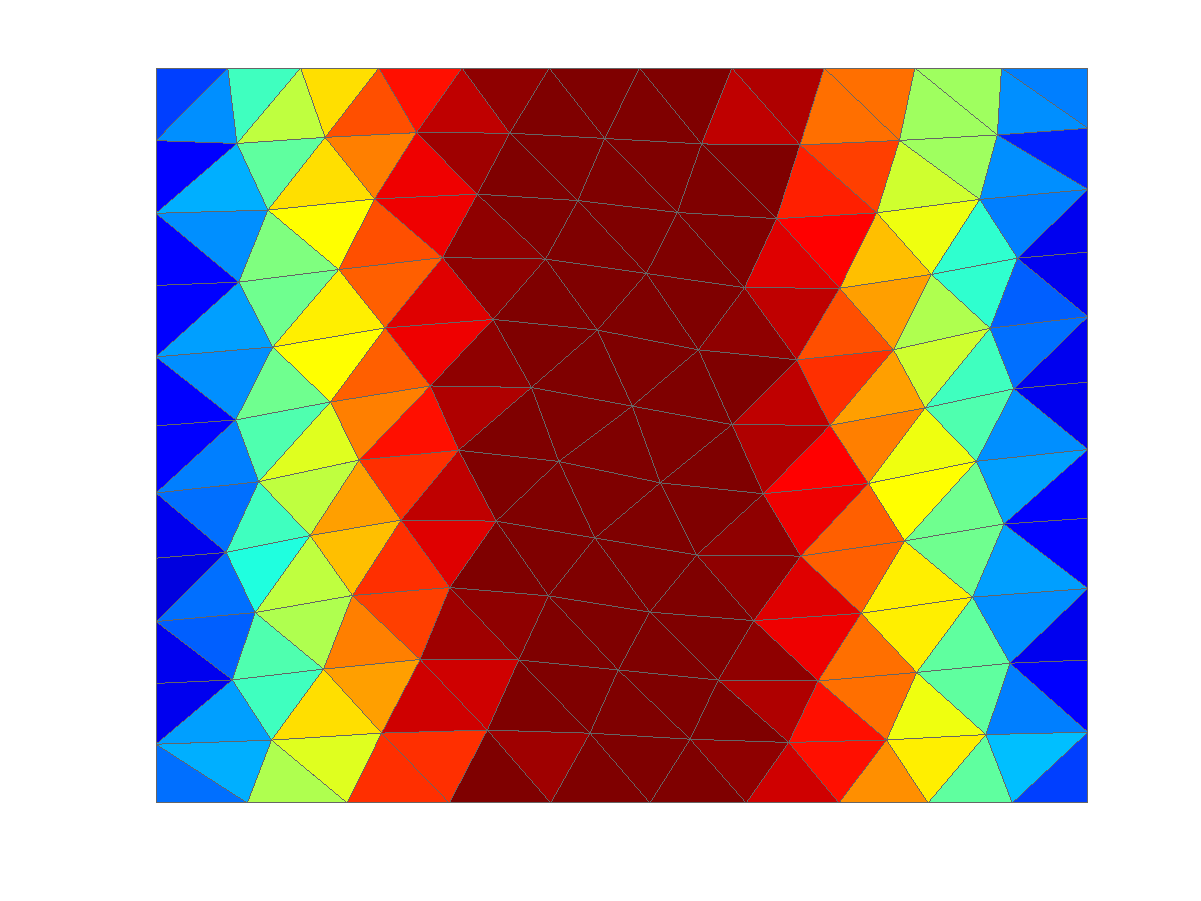}&
    \insertimage{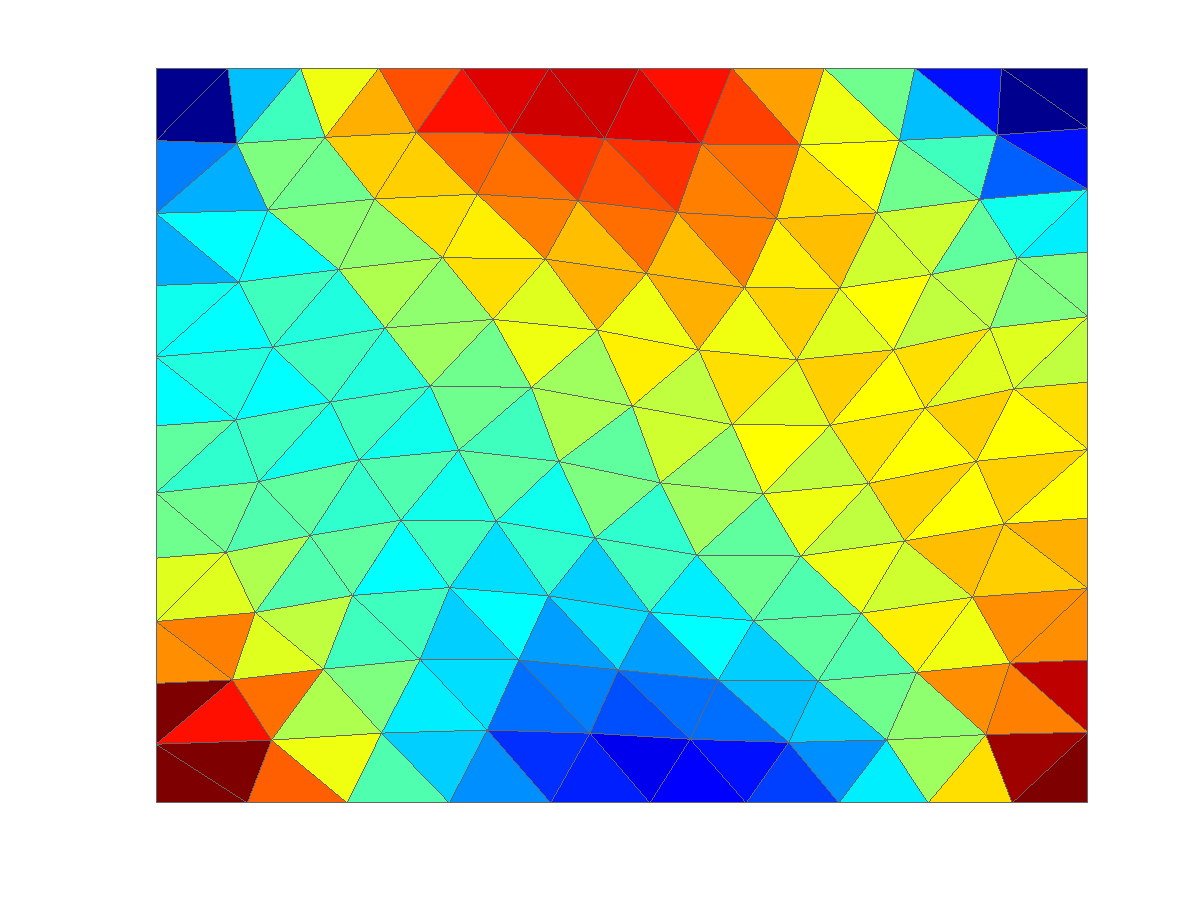}\\

    &
    \insertcolorbar{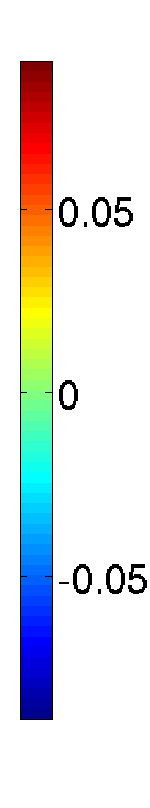}&
    \insertcolorbar{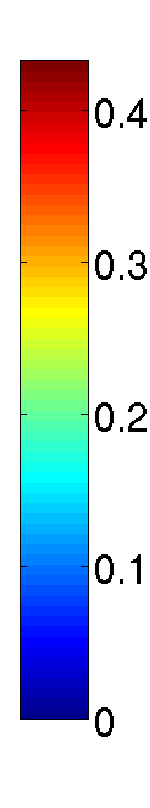}&
    \insertcolorbar{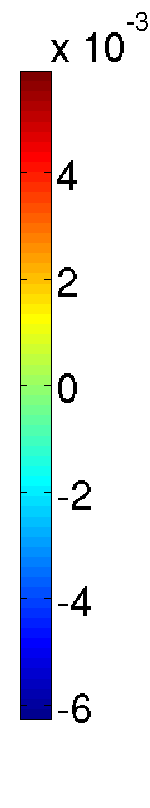}

  \end{tabular}
  \caption{Triangular grid and $\nu=0.495$ (Case \ref{case:limincomptriangle}). VEM
    and VEM-relax suffers from numerical locking. MPSA and VEM-relax-extra are free
    from numerical locking.}
\end{figure}

\begin{figure}[h]
  \centering
  \begin{tabular}{c@{}c@{}c@{}c@{}c@{}c}

    &
    \insertlegend{$d_x$}&
    \insertlegend{$d_y$}&
    \insertlegend{$\dive\vec{u}$}\\

    \insertcase{VEM}&
    \insertimage{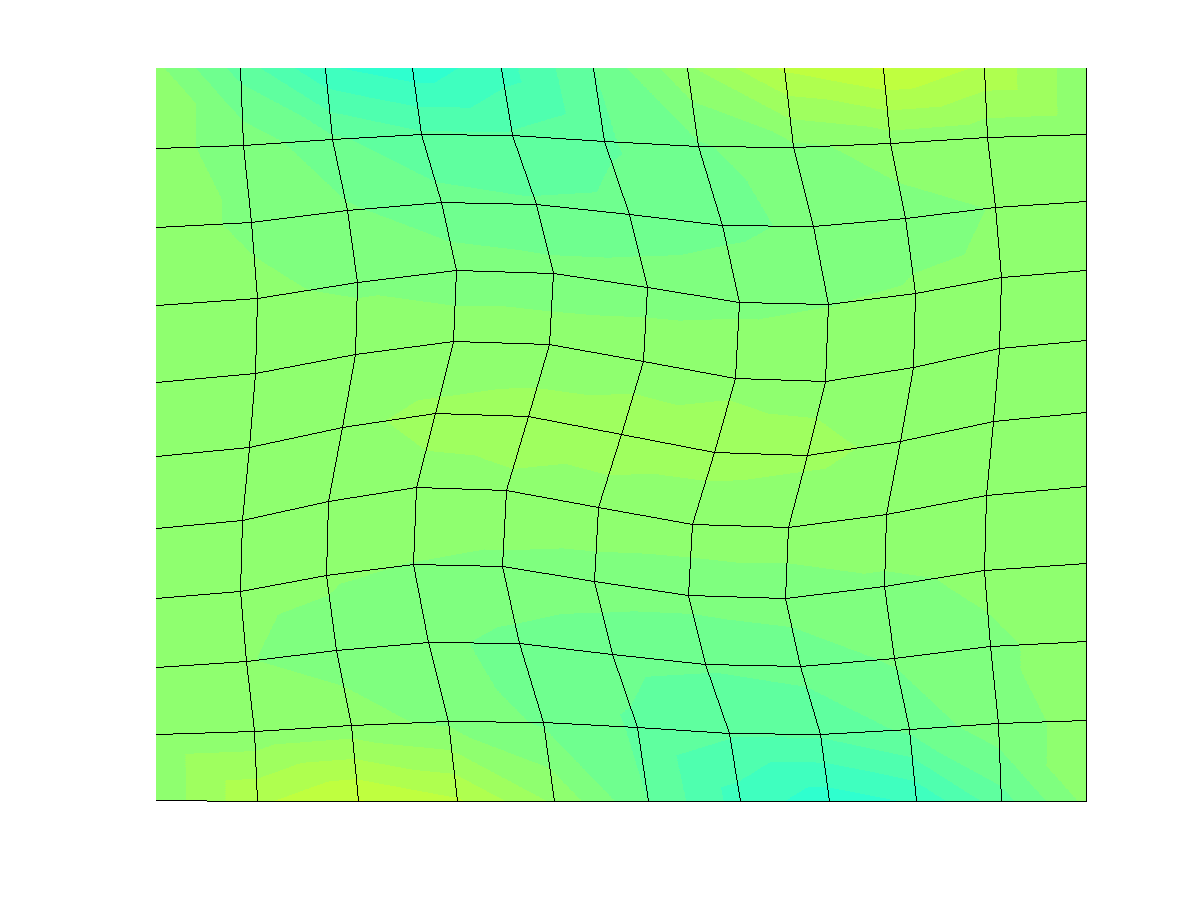}&
    \insertimage{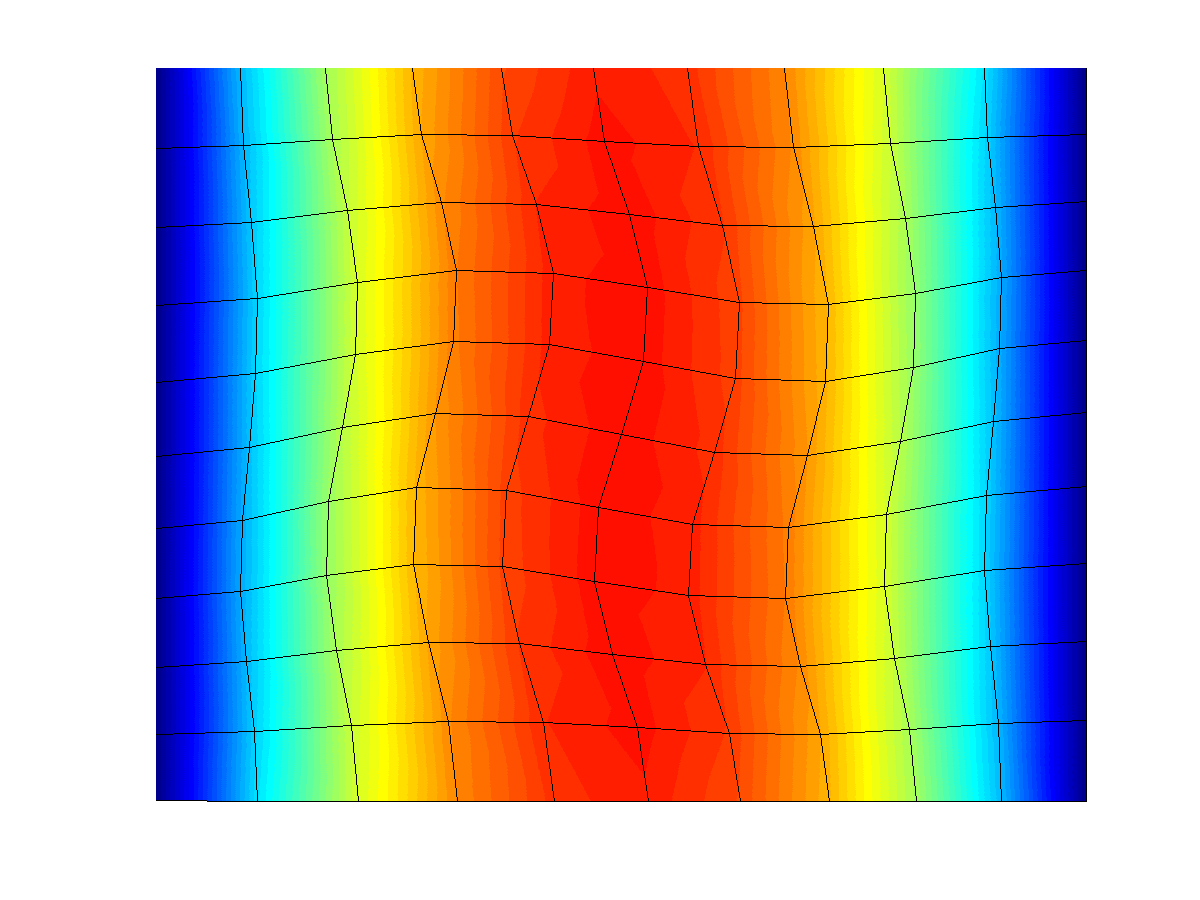}&
    \insertimage{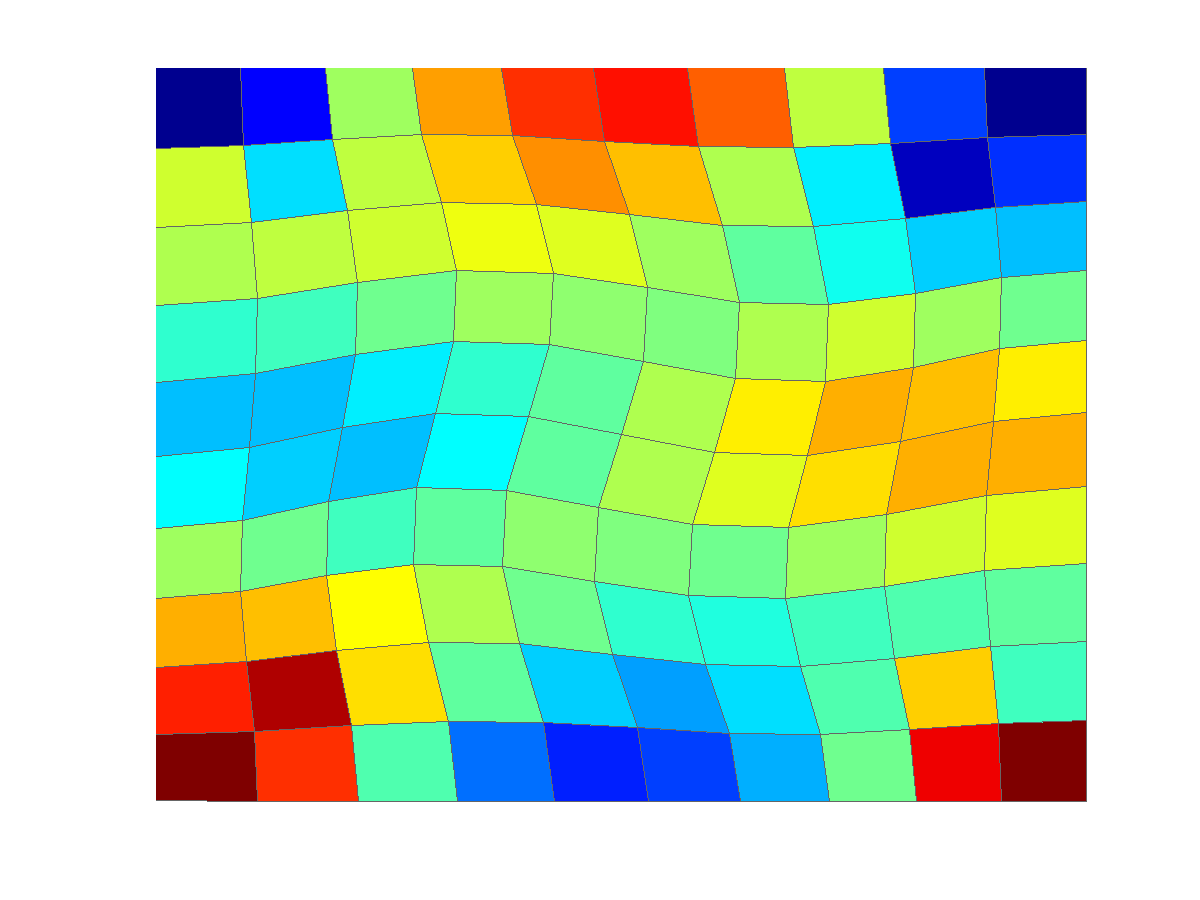}\\

    \insertcase{MPSA}&
    \insertimage{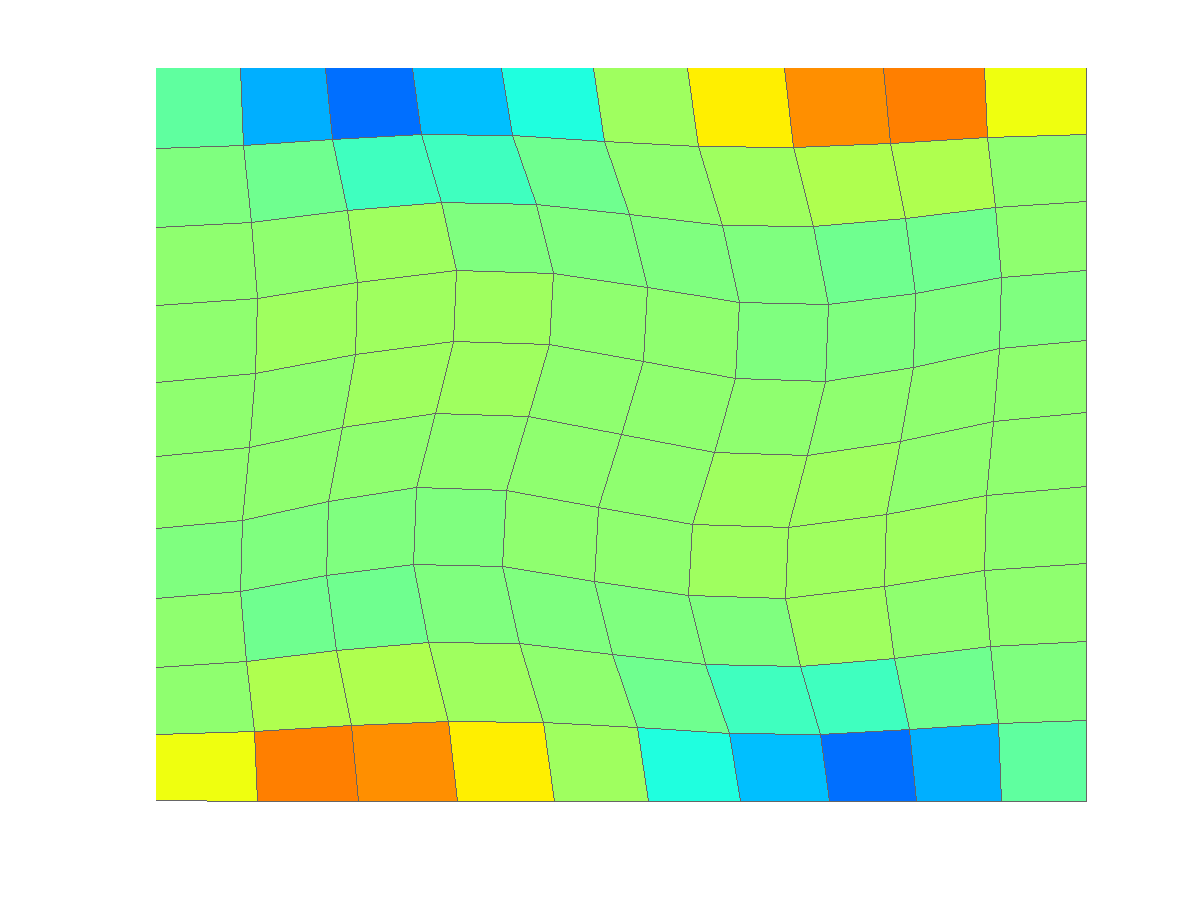}&
    \insertimage{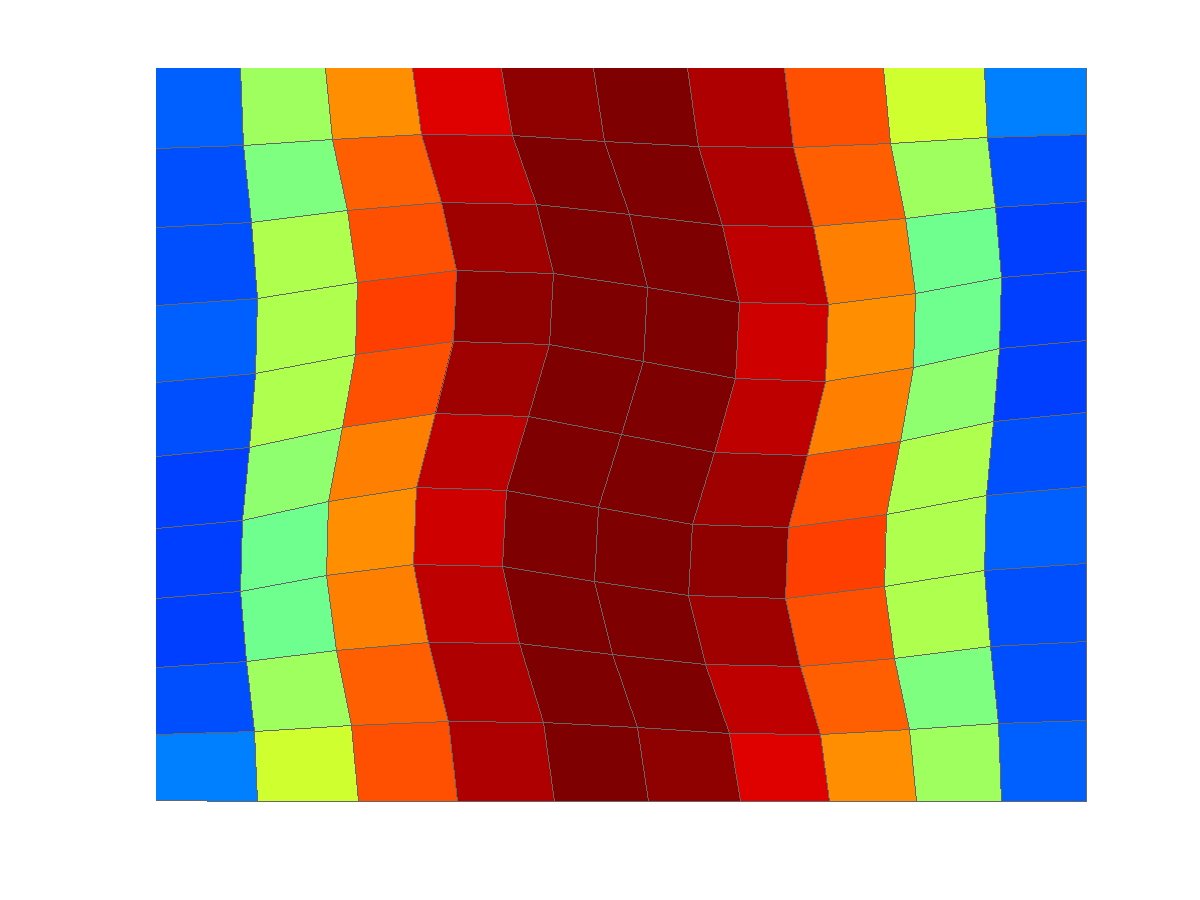}&
    \insertimage{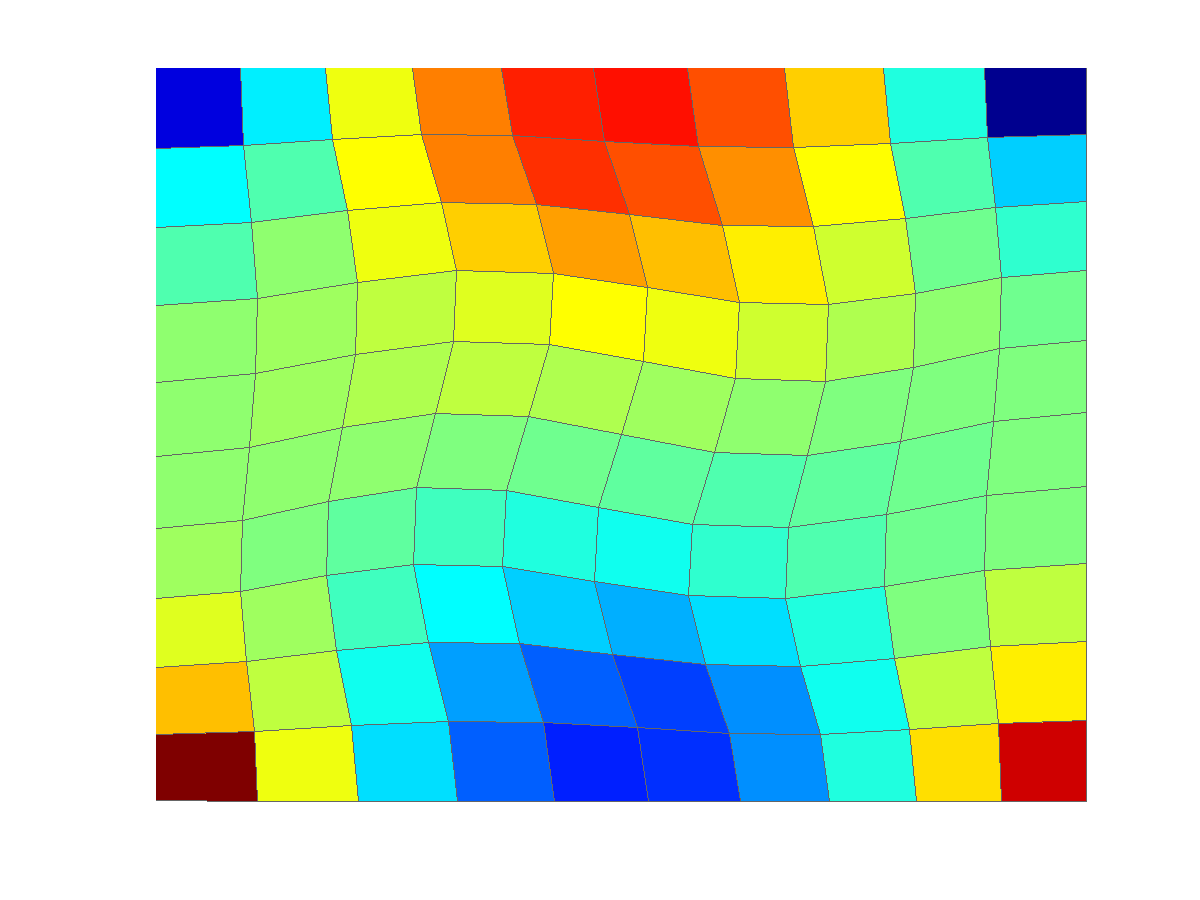}\\

    \insertcase{VEM relax}&
    \insertimage{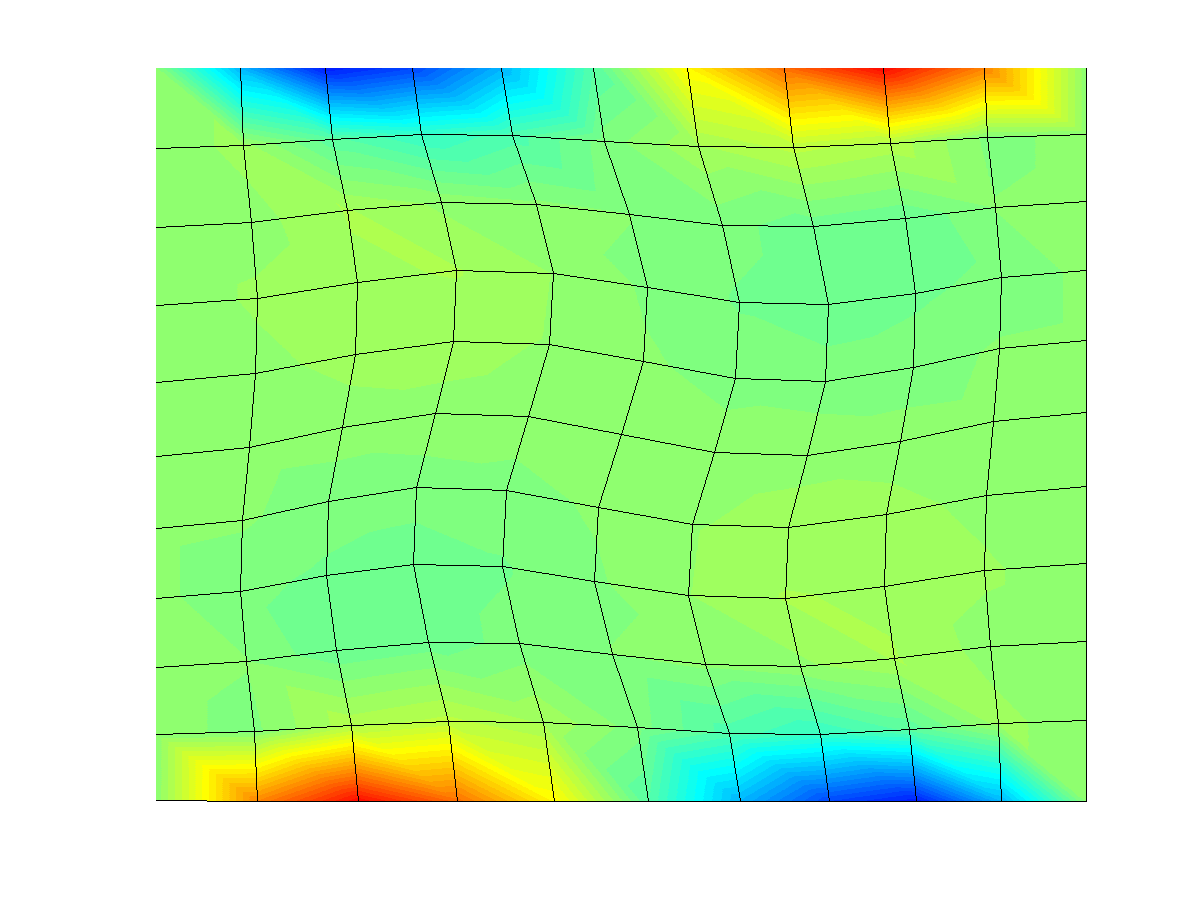}&
    \insertimage{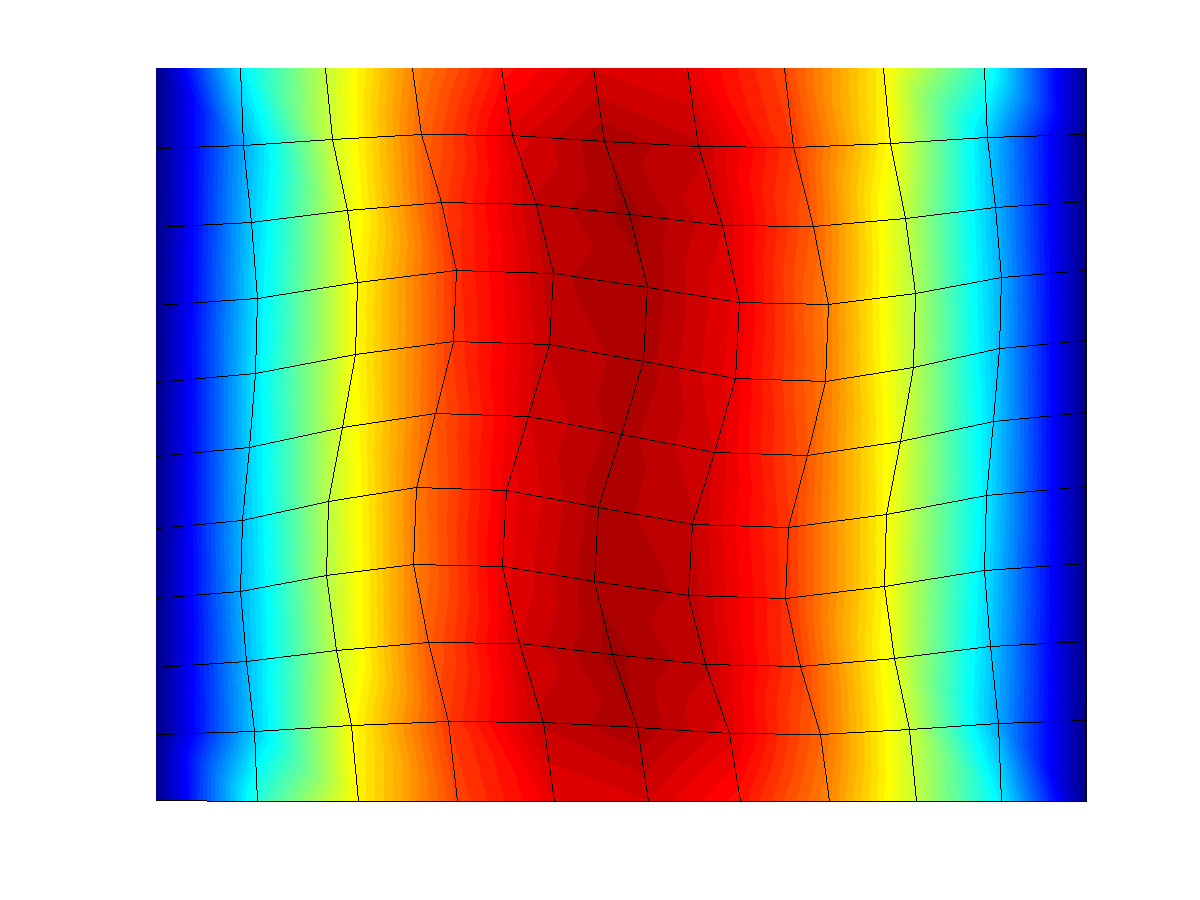}&
    \insertimage{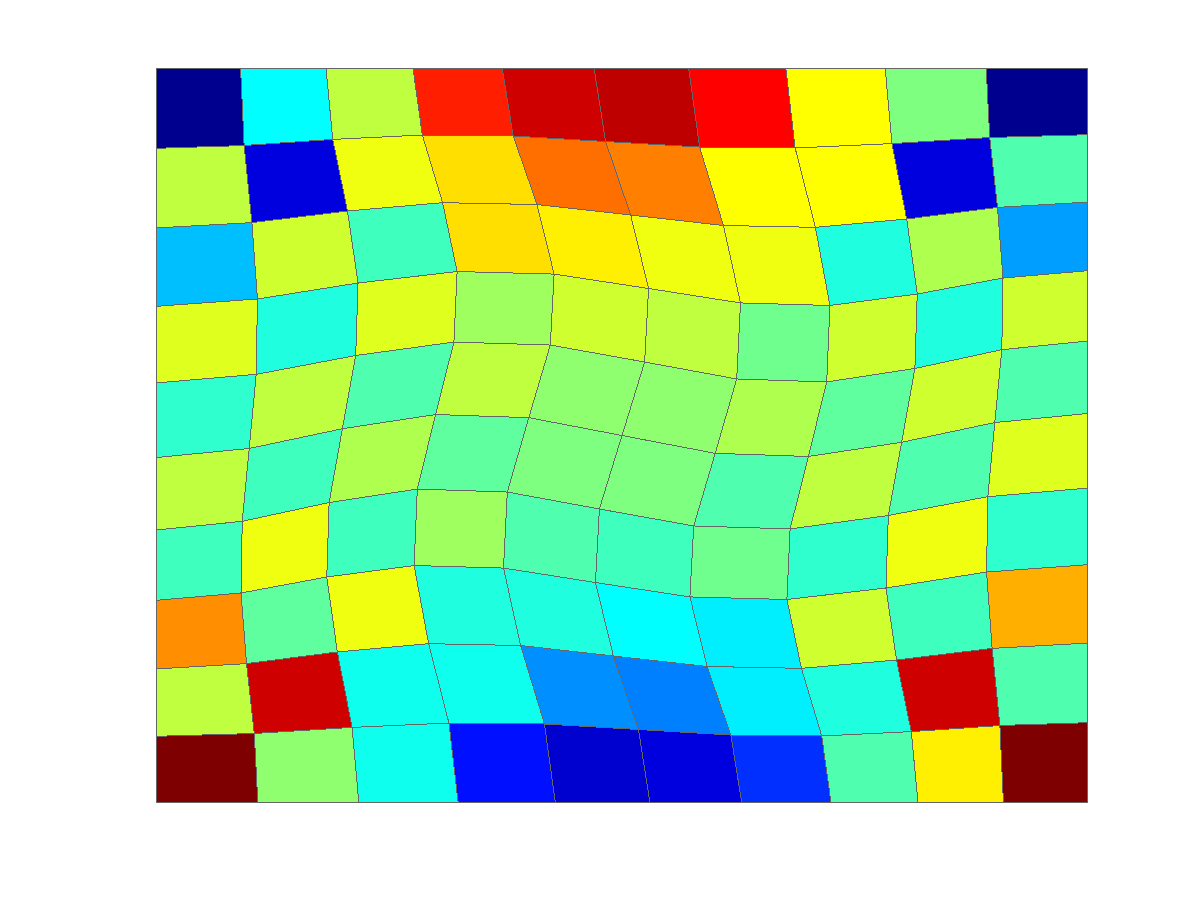}\\

    \insertcase{VEM relax extra}&
    \insertimage{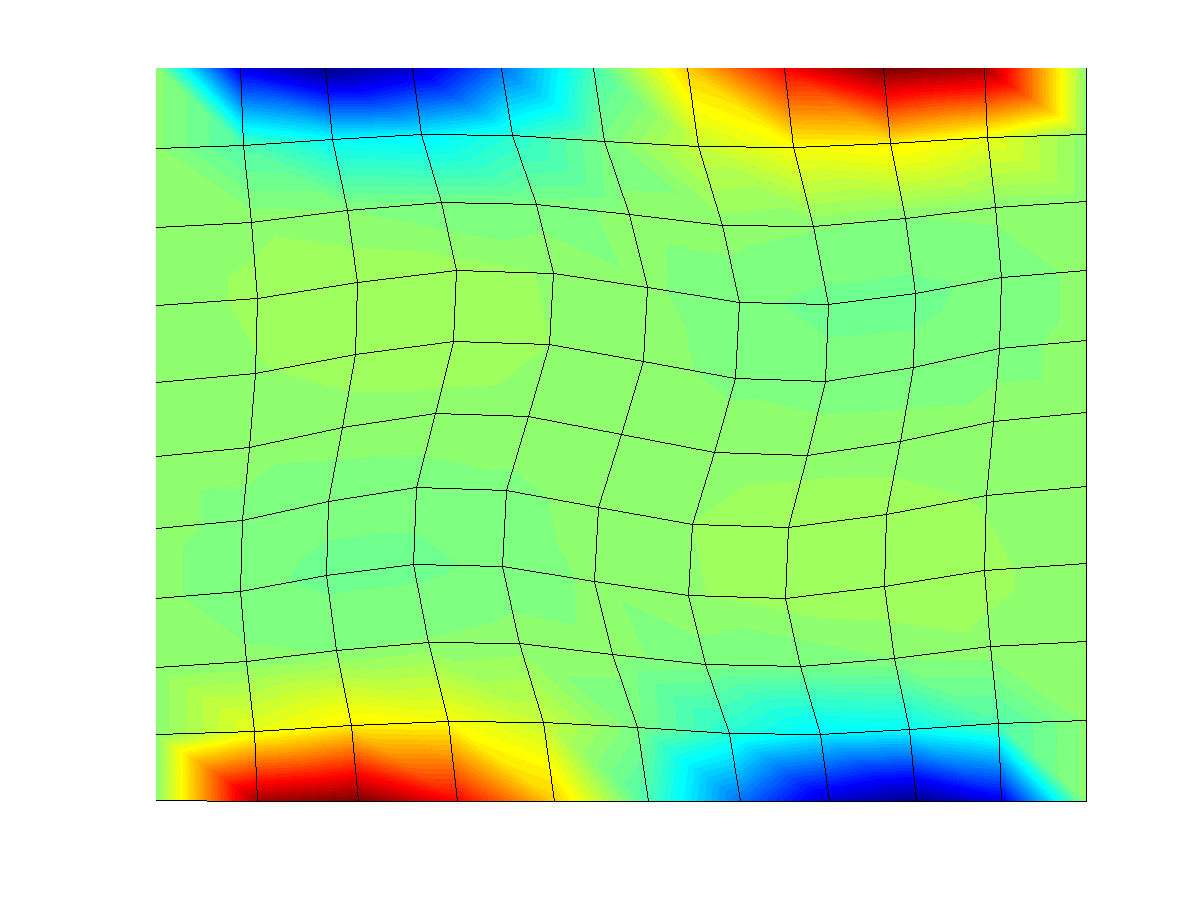}&
    \insertimage{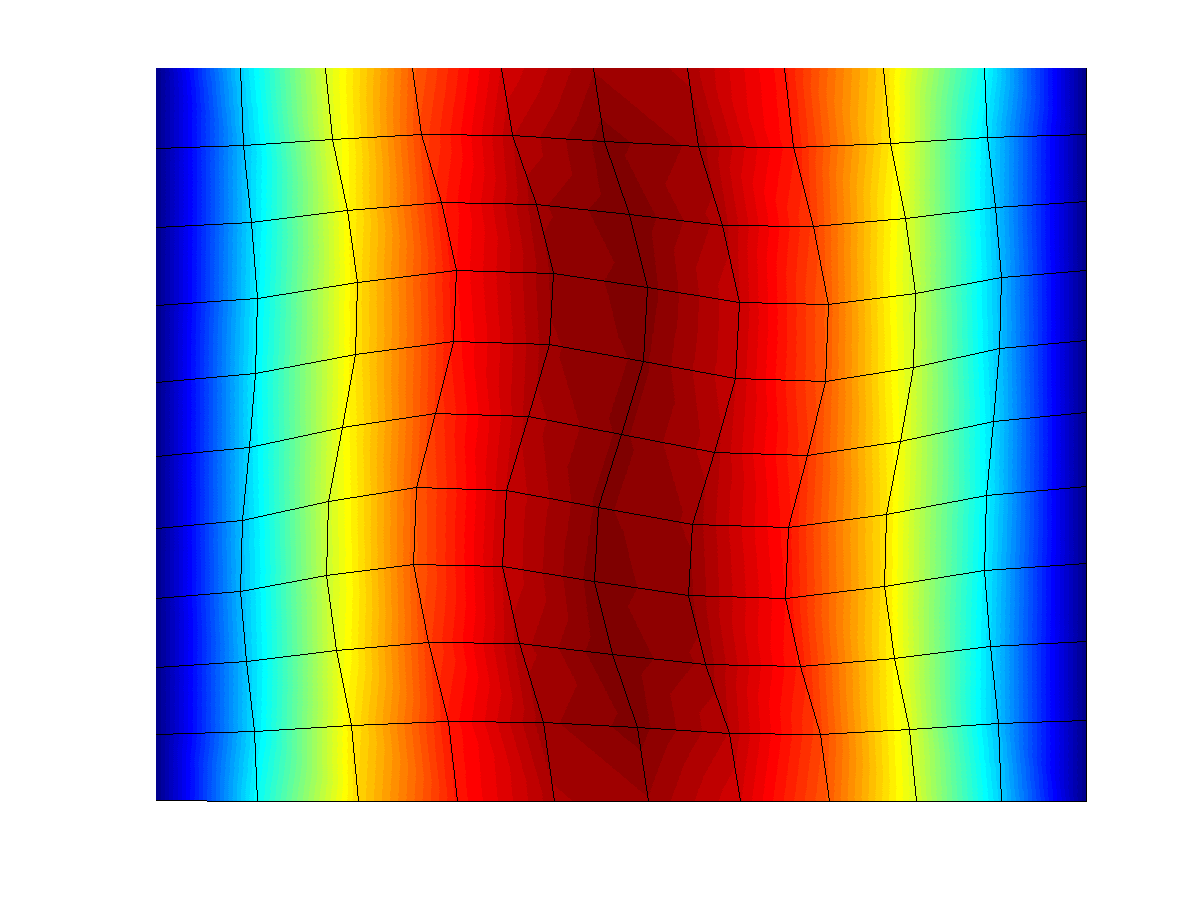}&
    \insertimage{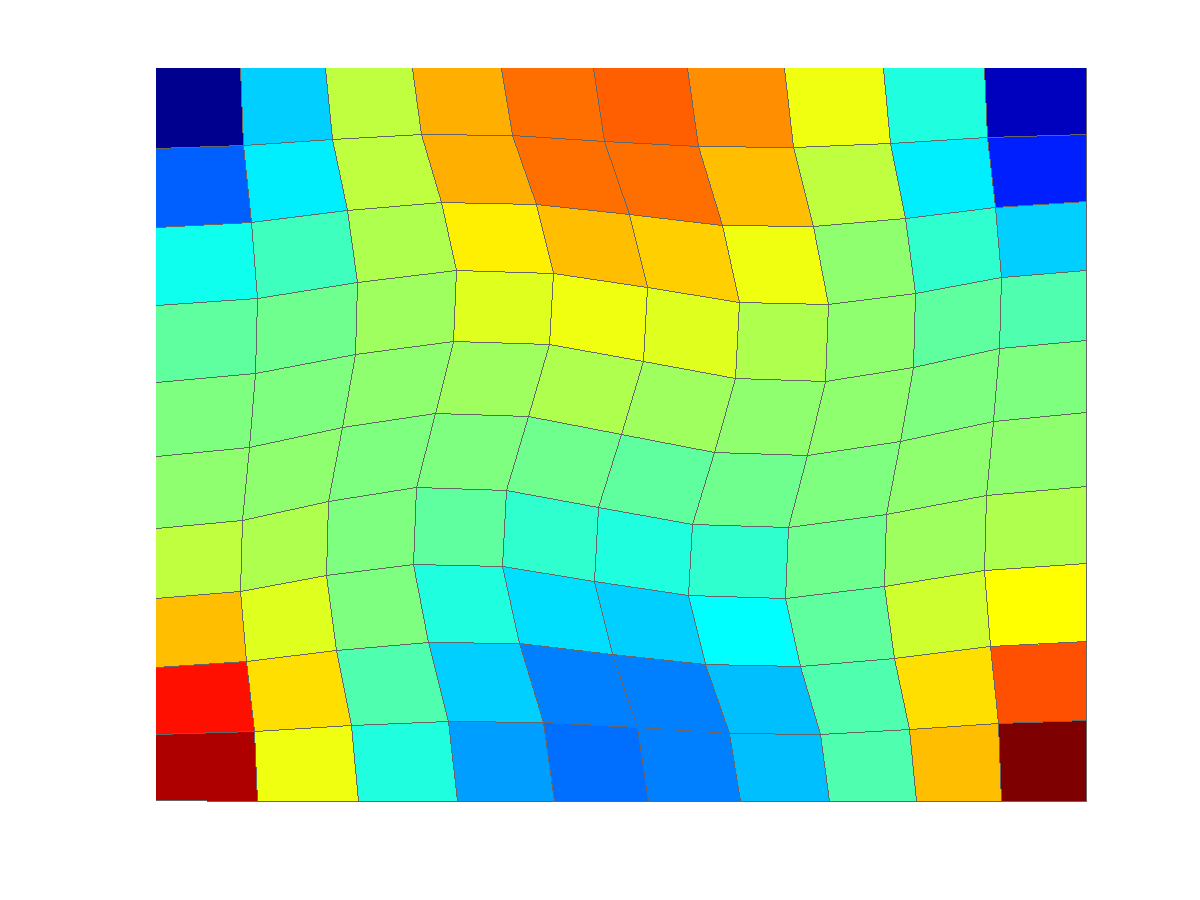}\\

    \insertcase{MPSA relax extra}&
    \insertimage{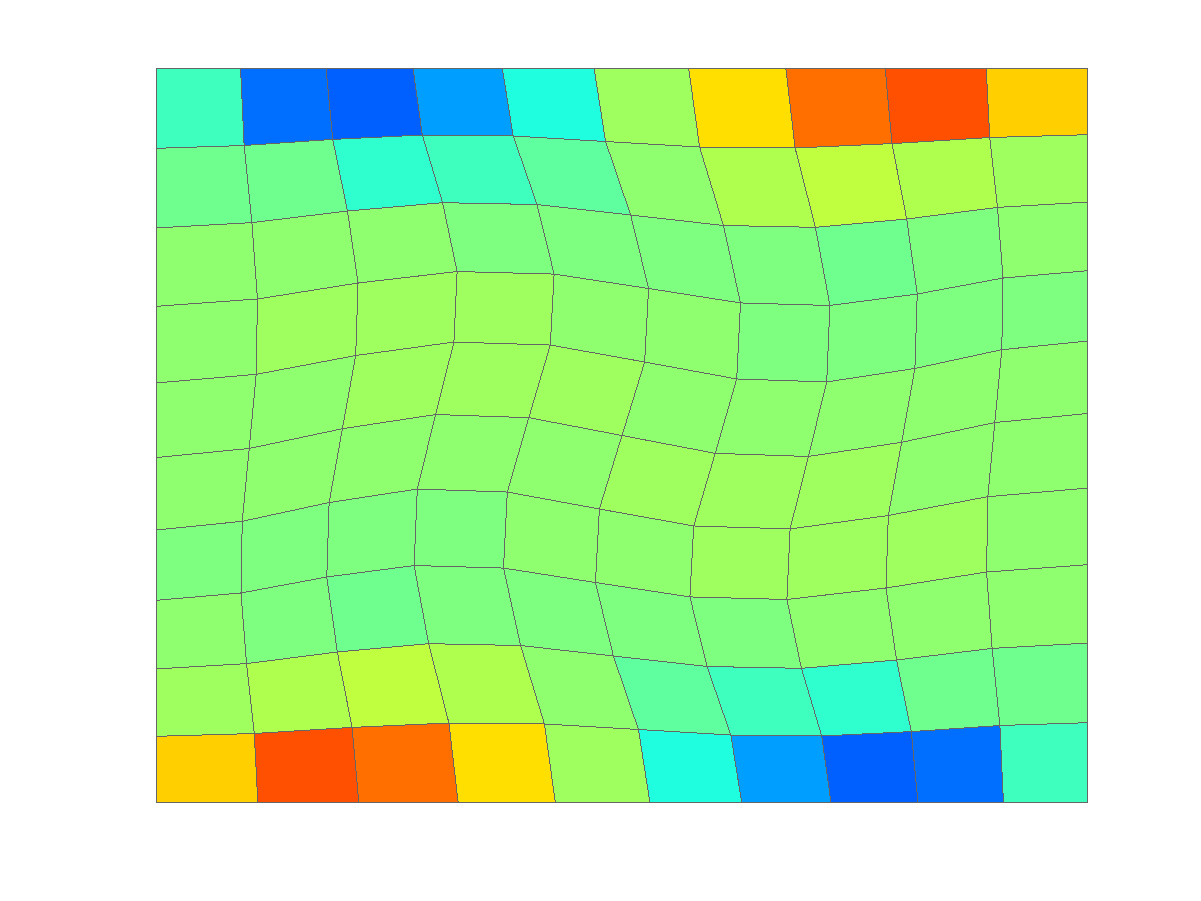}&
    \insertimage{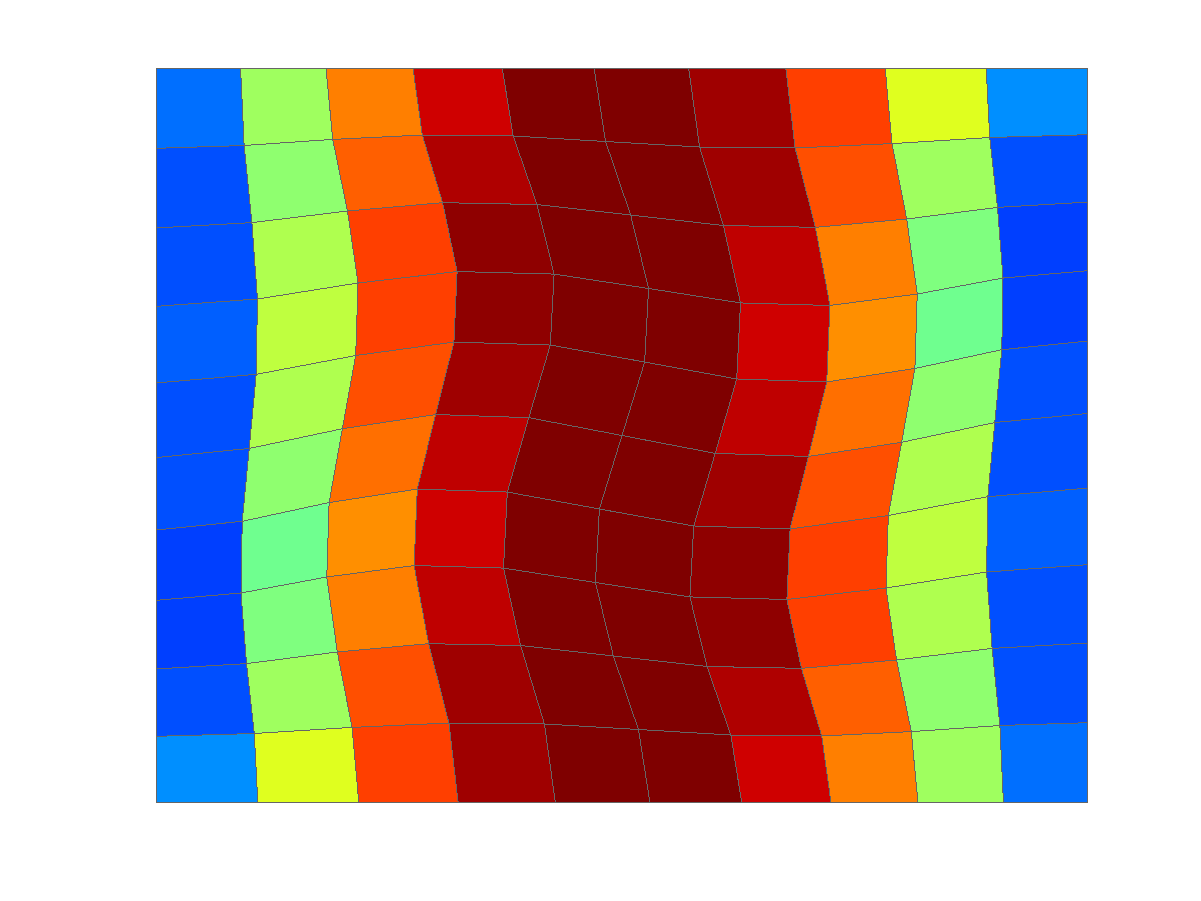}&
    \insertimage{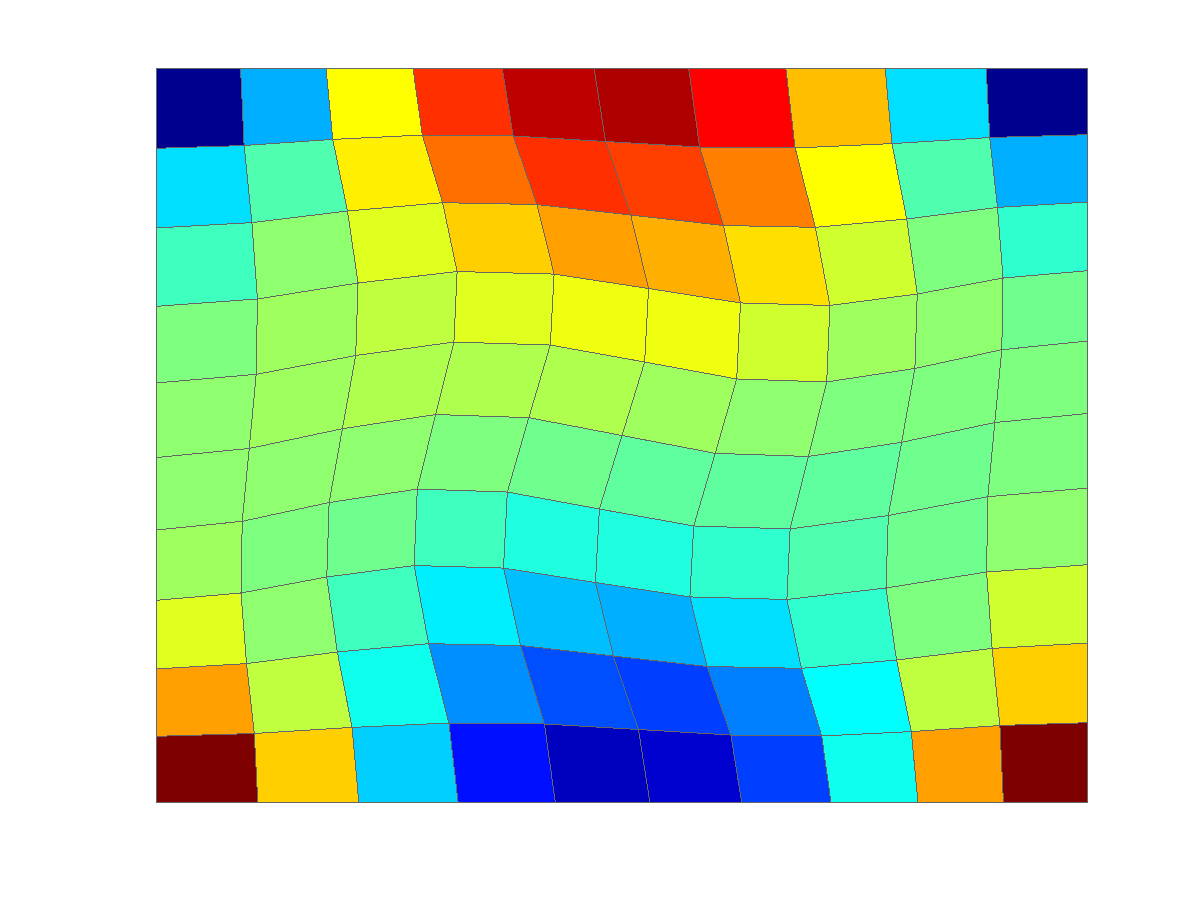}\\

    &
    \insertcolorbar{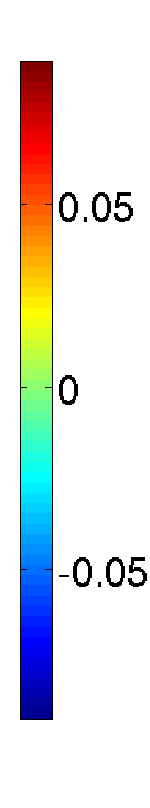}&
    \insertcolorbar{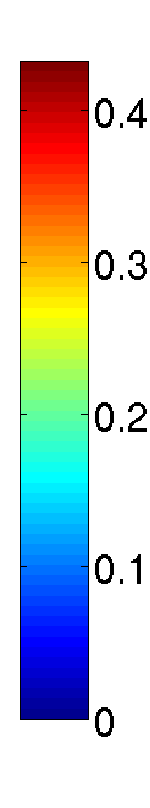}&
    \insertcolorbar{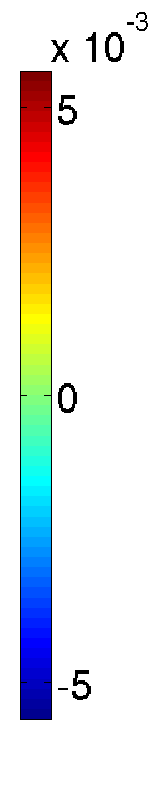}\\

  \end{tabular}
  \caption{Quadrilateral grid and $\nu=0.495$ (Case \ref{case:limincompquad}). VEM
    and VEM suffers from numerical locking. MPSA does not present any sign of
    locking. VEM relax extra is free from locking, as predicted by the theory}
\end{figure}
\endgroup

\clearpage

\begin{table}[h]
  \caption{Summary of the numerical tests}
  \begingroup
  \def\inserttc #1{\parbox[t][][t]{0.11\textwidth}{Case \ref{case:#1}}}
  \def\insertdes #1{\parbox[t][][t]{0.85\textwidth}{#1}}
  \def\arraystretch{1.2}
  \begin{center}
    \begin{tabular}[h]{l@{ : }c}
      \inserttc{fullcart} & \insertdes{Twisted grid}\\
      \inserttc{mixedgridplain} & \insertdes{Mixed grid with challenging features}\\
      \inserttc{stretchhexa} & \insertdes{Large aspect ratio with hexahedral grid}\\
      \inserttc{stretchtriangle} & \insertdes{Large aspect ratio with triangular grid}\\
      \inserttc{uniformVerticalRefinement} & \insertdes{Cartesian grid, with uniform refinement in the vertical direction}\\
      \inserttc{twoRegionsIsotropicRefinement} & \insertdes{Two regions, uniform refinement (both in $x$ and $y$) in one region}\\
      \inserttc{twoRegionsVerticalDirectionRefinement} & \insertdes{Two regions, refinement only in the $y$ direction in one region}\\
      \inserttc{refinedInterface} & \insertdes{Two regions with the same  Cartesian discretization, but with 20 extra nodes in each face at the interface \vspace*{1mm}}\\
      \inserttc{layerNoRefinement} & \insertdes{Vertical thin layer, no refinement inside the layer}\\
      \inserttc{layerRefinementCartesian} & \insertdes{Vertical thin layer with refinement inside the layer, Cartesian grid}\\
      \inserttc{layerRefinementTwisted} & \insertdes{Vertical thin layer with refinement inside the layer, twisted grid}\\
      \inserttc{limincomphexa} & \insertdes{$\nu=0.495$ with hexahedral grid}\\
      \inserttc{limincomptriangle} & \insertdes{$\nu=0.495$ with triangular grid}\\
      \inserttc{limincompquad} & \insertdes{$\nu=0.495$ with quadrilateral grid}
    \end{tabular}
  \end{center}
  \endgroup
  \label{tab:testcases}
\end{table}

\section{Conclusion}

In this paper, we have tested the behaviors of the VE and MPSA methods for linear
elasticity with respect to grid features and parameter values that are typical to
subsurface models. We can conclude that both methods are capable to handle in a
satisfactory manner the polyhedral grid structures that are standard in such
models. A priori, both methods have relative advantages. The MPSA method is
attractive from the physical point of view, due to the explicit treatment of the
force continuity at the cell interfaces. The MPSA offers a natural stable coupling
with poro-elasticity, see \cite{NorbottenStable2015}. From the implementation point
of view, the MPSA method is cell-centered and therefore shares the same grid
structure as the MPFA method, which is also often the preferred convergent method for
multi-phase flow. The VE method has the advantage of robustness. Obtained from a
variational principle, it is always symmetric definite positive. For simplexes, the
method reduces to the finite element method so that the large collection of
techniques developed for finite elements, such as preconditioning, super-convergence
techniques or patch recovery can be relatively easily applied to VEM. When we use the
approach presented in \cite{gain2014} as we do in this paper, the projection
operators are given explicitly and do not require extra local computations as they
usually do in a VE method, so that the local assembly has finally the same structure
as in the traditional finite element method. If one is interested in generating the
deformed grid obtained from a computed displacement field, another advantage of the
VE method is that that the deformed grid can be readily constructed as the method
yields nodal displacement. In comparison, the MPSA method requires an extra
post-processing step and the task of generating a deformed grid from displacements at
cell centers is not trivial. There is no explicit reconstruction of continuous
displacement field from values given at cells.

In geological models, the convergence properties of a method are not as important as
its performance on coarse and strongly irregular meshes.In a first series of test, we
have checked the convergence of the method for randomly perturbed quadrilateral
grids. Then, we study the behavior of the method on strongly distorted grid and grids
with high aspect ratio. At this stage, we reach the limit of both methods. For the
MPSA, we exceed the grid restriction of the method. The VE method is robust but
convergence is guaranteed with a uniform bound on the aspect ratio. At high aspect
ratio, it is therefore not surprising that we observe discrepancies of the
solutions. At the same time, in the examples we have tested, we observe that the
displacement field is not that strongly affected and could be used directly. The
pattern of the large oscillations in the divergence term can be understood and opens
for the possibility of a post-processing in the spirit of \cite{zienkiewicz92SPR},
which would enable us to used those values as well.

We have studied the behavior of both methods for grids containing two regions with
different refinements (cases \ref{case:twoRegions}). The first conclusion is that
both methods are robust with respect to the refinement ratio in averaged norms ($L^2$
norms). As the refinement ratio is increased, oscillations in the forces at the
interaction region appear for the VE method and, in the case of isotropic refinement
(case \ref{case:twoRegionsIsotropicRefinement}), the local values for the stress even
blow-up. We interpret this behavior by the highly non-linear nature of the virtual
basis in this case. The freedom we have in choosing the regularization term in VEM
could be used to dampen these unwanted effects but we do not explore this possibility
in this paper. The MPSA does not present the same level of oscillations in the force
term and yields more reliable values for the forces, which is in accordance with the
fact that the method is based on a force continuity principle. We have studied the
behavior of the methods in the case of a thin layer. The conclusion is that both
methods are robust with respect to the thickness of the layer. When the layer is
refined, the VE method introduces, as previously, oscillation in the forces at the
interface but also in the divergence term inside the layer. The error term in the
divergence remained confined to the layer in VEM while it is spread for MPSA.

Even if the rocks considered in a subsurface model are far from incompressible, the
coupling with fluid flow requires that the methods used for elasticity are robust
with respect to the incompressible limit and, in particular, not sensitive to
numerical locking. We have conducted tests for both methods. First, we confirm the
intuition that the MSPA and VEM methods have opposite responses to numerical locking
depending on the grid structure: In VEM, numerical locking will appear for grids with
relatively more cells than nodes (such as triangular grids) and, in MPSA, it will
appear for grids with relatively more nodes than cells (such as hexahedral grids). We
can get rid of numerical locking for VEM using established theory, as presented in
\cite{da2014mimetic} for the Stokes problem, by relaxing the divergence term and
adding extra degrees of freedom. For MPSA, by introducing an extra degree of freedom
at cell center, which correspond to pressure, it is also possible to derive a method
that is robust with respect to locking. All these methods have been tested and the
regularized methods fulfill the expectation we have concerning locking. Moreover, we
can conclude from our experiments that, using the standard versions, the MPSA method
seems more robust than the VE method with respect to locking. For example, MPSA can
handle quadrilateral grids where VEM fails.

\printbibliography

\end{document}